%% file: RR-8857.tex
\documentclass[twoside]{article}

\input{command}
\newcommand{\myclearpage}{\clearpage}
\newcommand{\myparagraph}[1]{\paragraph{#1.}}
\newcommand{\figloc}{htb}

\input{preamble}

\begin{document}

\RRNo{8857}
\makeRR

\tableofcontents

\input{body}

\bibliography{biblio}
\bibliographystyle{plain}

\end{document}

%% file: command.tex
\usepackage[english]{babel}
\usepackage[utf8]{inputenc} 
\usepackage{amsmath}
\usepackage{amsfonts}
\usepackage{graphicx}
\usepackage{color}
\usepackage{hyperref}
\usepackage{stmaryrd}

\allowdisplaybreaks

\graphicspath{{./}}

\newcommand{\twofigs}{0.45\textwidth}
\newcommand{\threefigs}{0.28\textwidth}

\newcommand{\tentative}{tentative}
\newcommand{\candidate}{candidate}

\newtheorem{Remark}{Remark}[section]

\newcommand{\eg}{{\em e.g.}}
\newcommand{\ie}{{\em i.e.}}
\newcommand{\apriori}{{\em a priori\/}}
\newcommand{\Apriori}{{\em A priori\/}}

\newcommand{\suchthat}{\,|\,}

\newcommand{\matR}{{\mathbb R}}

\newcommand{\tendsto}[1]{\rightarrow #1}
\newcommand{\tendstoz}{\tendsto{0}}
\newcommand{\half}{\frac{1}{2}}
\newcommand{\dps}{\displaystyle}
\newcommand{\derpar}[2]{\frac{\partial #1}{\partial #2}}
\newcommand{\jump}[1]{\llbracket#1\rrbracket}
\newcommand{\egaldef}{\stackrel{{\rm def}}{=}}
\newcommand{\alf}{\alpha}
\newcommand{\bet}{\beta}
\newcommand{\eps}{\varepsilon}
\newcommand{\fhi}{\varphi}
\newcommand{\fhip}{\fhi^+}
\newcommand{\fhim}{\fhi^-}

\newcommand{\nn}{\textbf{n}}
\newcommand{\nnp}{\nn^+}
\newcommand{\nnm}{\nn^-}
\newcommand{\ggamma}{\boldsymbol{\gamma}}
\newcommand{\nnu}{\nu}
\newcommand{\mmu}{\mu}

\newcommand{\vu}{\textbf{u}}
\newcommand{\vup}{\vu^+}
\newcommand{\vum}{\vu^-}
\newcommand{\ppp}{p^+}
\newcommand{\ppm}{p^-}
\newcommand{\pmeas}{\overline{p}}

\newcommand{\gamfrac}{\ggamma}
\newcommand{\gamfracp}{\gamfrac^+}
\newcommand{\gamfracm}{\gamfrac^-}

\newcommand{\pOmega}{\partial\Omega}
\newcommand{\pgamfrac}{\partial\gamfrac}

\newcommand{\ff}{f}
\newcommand{\ffT}{\ff_T}
\newcommand{\ffE}{\ff_E}

\newcommand{\ugf}{\vu_{\gamfrac}}
\newcommand{\fgf}{\ff_{\gamfrac}}
\newcommand{\pgf}{p_{\gamfrac}}
\newcommand{\ngf}{\nn_{\gamfrac}}

\newcommand{\K}{\textbf{K}}
\newcommand{\kh}{\kappa}

\newcommand{\grad}{\nabla}
\newcommand{\gradtau}{\grad_\tau}
\renewcommand{\div}{{\rm div}}
\newcommand{\divtau}{\div_\tau}

\newcommand{\alfE}{\alf_E}
\newcommand{\alfbar}{\overline{\alf}}
\newcommand{\alfbarE}{\alfbar_E}
\newcommand{\betE}{\bet_E}
\newcommand{\betbar}{\overline{\bet}}
\newcommand{\betbarE}{\betbar_E}
\newcommand{\epsalf}{\eps_\alf}
\newcommand{\epsbet}{\eps_\bet}
\newcommand{\xiE}{\xi_E}
\newcommand{\xibar}{\overline{\xi}}
\newcommand{\xibarE}{\xibar_E}
\newcommand{\eeps}{\boldsymbol{\eps}}
\newcommand{\epsof}[1]{\eeps^{#1}}

\newcommand{\epsk}{\epsof{k}}
\newcommand{\epsofcomp}[2]{\eps_{#2}^{#1}}
\newcommand{\epsmkofcomp}[1]{\epsofcomp{k-1}{#1}}
\newcommand{\epsmkl}{\epsmkofcomp{\ell}}
\newcommand{\epskofcomp}[1]{\epsofcomp{k}{#1}}
\newcommand{\epskl}{\epskofcomp{\ell}}
\newcommand{\epstof}[1]{\tilde{\eeps}^{#1}}
\newcommand{\epstk}{\epstof{k}}
\newcommand{\epstofcomp}[2]{\tilde{\eps}_{#2}^{#1}}
\newcommand{\epstkofcomp}[1]{\epstofcomp{k}{#1}}

\usepackage{mathrsfs} 

\newcommand{\calE}{\mathscr{E}}
\newcommand{\calF}{\mathscr{F}}
\newcommand{\calI}{\mathcal{I}}
\newcommand{\calJ}{\mathcal{J}}
\newcommand{\calL}{\mathcal{L}}
\newcommand{\calN}{\mathscr{N}}
\newcommand{\calQ}{\mathcal{Q}}
\newcommand{\calT}{\mathcal{T}}

\newcommand{\Th}{\calT_h}
\newcommand{\Thc}{\Th^c}
\newcommand{\Thm}{\Th^m}

\newcommand{\ET}{\calE_T}
\newcommand{\TE}{\calT_E}
\newcommand{\Eh}{\calE_h}
\newcommand{\Eih}{\Eh^i}
\newcommand{\Eeh}{\Eh^e}

\newcommand{\EN}{\calE_N}
\newcommand{\NE}{\calN_E}
\newcommand{\Fh}{\calF_h}
\newcommand{\Lag}{\calL}
\newcommand{\Nalf}{\calN_\alf}

\newcommand{\PT}{P_T}
\newcommand{\Pt}{P_t}
\newcommand{\PE}{P_E}
\newcommand{\Pe}{P_e}
\newcommand{\PN}{P_N}
\newcommand{\Pn}{P_n}
\newcommand{\UTE}{U_{T,E}}
\newcommand{\Ute}{U_{te}}
\newcommand{\UEN}{U_{E,N}}

\newcommand{\VEN}{V_{E,N}}
\newcommand{\Ven}{V_{en}}

\newcommand{\PmT}{\overline{P}_T}
\newcommand{\PmE}{\overline{P}_E}
\newcommand{\UmTE}{\overline{U}_{T,E}}

\newcommand{\LAT}{\lambda_T}
\newcommand{\LAt}{\lambda_t}
\newcommand{\LAE}{\lambda_E}
\newcommand{\LAe}{\lambda_e}
\newcommand{\LAN}{\lambda_N}
\newcommand{\LAn}{\lambda_n}
\newcommand{\LAbarN}{\overline{\lambda}_N}
\newcommand{\LAbarn}{\overline{\lambda}_{n}}
\newcommand{\MUte}{\mmu_{te}}
\newcommand{\MUTE}{\mu_{T,E}}
\newcommand{\MUen}{\mmu_{en}}
\newcommand{\MUEN}{\mu_{E,N}}
\newcommand{\NUen}{\nnu_{en}}
\newcommand{\NUEN}{\nu_{E,N}}

\newcommand{\Cost}{\calJ}
\newcommand{\J}[1]{\Cost_{#1}}
\newcommand{\JcFh}{\J{\currFh}}
\newcommand{\JtFh}{\J{\tentFh}}

\newcommand{\DJtFhepsalf}{\derpar{\JtFh}{\epsalf}}
\newcommand{\DJtFhepsbet}{\derpar{\JtFh}{\epsbet}}
\newcommand{\Xeps}{X(\epsalf,\epsbet)}
\newcommand{\Lambdaeps}{\Lambda(\epsalf,\epsbet)}
\newcommand{\DLagX}{\derpar{\Lag}{X}}
\newcommand{\DLageps}{\derpar{\Lag}{(\epsalf,\epsbet)}}
\newcommand{\Diff}{\delta}

\newcommand{\Falf}{\calF_\alf}
\newcommand{\Fbet}{\calF_\bet}
\newcommand{\F}{\calF}

\newcommand{\accr}[1]{\Delta #1}
\newcommand{\DFh}{\accr{\Fh}}
\newcommand{\DFalf}{\accr{\Falf}}
\newcommand{\DFbet}{\accr{\Fbet}}
\newcommand{\DF}{\accr{\F}}
\newcommand{\DNalf}{\accr{\Nalf}}

\newcommand{\currFh}{\Fh}
\newcommand{\currFalf}{\Falf}
\newcommand{\currFbet}{\Fbet}
\newcommand{\currF}{\F}
\newcommand{\currFc}{\currF^c}
\newcommand{\currNalf}{\Nalf}
\newcommand{\tentFh}{\DFh}
\newcommand{\tentFalf}{\DFalf}
\newcommand{\tentFbet}{\DFbet}
\newcommand{\tentF}{\DF}
\newcommand{\tentFc}{\tentF^c}
\newcommand{\tentNalf}{\DNalf}
\newcommand{\fullFh}{\currFh\cup\tentFh}
\newcommand{\fullFalf}{\currFalf\cup\tentFalf}
\newcommand{\fullFbet}{\currFbet\cup\tentFbet}
\newcommand{\fullNalf}{\currNalf\cup\tentNalf}

\newcommand{\currFof}[1]{\currF^{#1}}
\newcommand{\currFml}{\currFof{\ell-1}}
\newcommand{\currFmk}{\currFof{k-1}}
\newcommand{\currFk}{\currFof{k}}

\newcommand{\besttentFof}[1]{\tentF^{\star #1}}
\newcommand{\besttentFl}{\besttentFof{\ell}}
\newcommand{\besttentFmk}{\besttentFof{k-1}}
\newcommand{\besttentFk}{\besttentFof{k}}
\newcommand{\bestepsof}[1]{\epsof{\star #1}}
\newcommand{\bestepsmk}{\bestepsof{k-1}}
\newcommand{\bestepsk}{\bestepsof{k}}
\newcommand{\bestepsofcomp}[2]{\epsofcomp{\star #1}{#2}}
\newcommand{\bestepsmkofcomp}[1]{\bestepsofcomp{k-1}{#1}}
\newcommand{\bestepskofcomp}[1]{\bestepsofcomp{k}{#1}}

\newcommand{\bestcostof}[1]{\Cost^{\star #1}}
\newcommand{\bestcostmk}{\bestcostof{k-1}}
\newcommand{\bestcostk}{\bestcostof{k}}

\newcommand{\indic}{\calI}
\newcommand{\Ialf}{\indic_\alf}
\newcommand{\Ibet}{\indic_\bet}
\newcommand{\I}{\indic}
\newcommand{\bestIalfelem}{\Ialf^{\star,{\rm elem}}}
\newcommand{\bestIalfaggr}{\Ialf^{\star,{\rm aggr}}}
\newcommand{\bestIalfext}{\Ialf^{\star,{\rm ext}}}

\newcommand{\bestIelem}{\I^{\star,{\rm elem}}}
\newcommand{\bestIaggr}{\I^{\star,{\rm aggr}}}
\newcommand{\bestIext}{\I^{\star,{\rm ext}}}

\newcommand{\Ec}{E^c}

\newcommand{\Nx}{N_x}
\newcommand{\Ny}{N_y}
\newcommand{\Nm}{N_m}
\newcommand{\Nf}{N_f}
\newcommand{\Quad}[2]{\calQ^{#1,#2}}
\newcommand{\QuadNxy}{\Quad{\Nx}{\Ny}}
\newcommand{\Quadcomp}{\Quad{72}{72}}
\newcommand{\Quadmeas}{\Quad{\Nm}{\Nm}}
\newcommand{\Quadfrac}{\Quad{\Nf}{\Nf}}

\newcommand{\Dmeas}{\Delta_{\rm meas}}
\newcommand{\etaconv}{\eta_{\rm conv}}
\newcommand{\etastat}{\eta_{\rm stat}}
\newcommand{\kmax}{k_{\max}}
\newcommand{\thetaelem}{\theta_{\rm elem}}
\newcommand{\thetaext}{\theta_{\rm ext}}
\newcommand{\nmax}{n_{\max}}

%% file: preamble.tex
\usepackage[a4paper]{geometry}
\usepackage[T1]{fontenc} 
\usepackage{RR}

\usepackage{a4wide}

\RRdate{February 2016}

\RRauthor{
  Hend {Ben Ameur}%
  \thanks[lamsin]{%
    LAMSIN, ENIT, University of Tunis, BP 37, Le Belvédère,
    TN-1002 Tunis, Tunisia.
    {\tt hbenameur@yahoo.ca, cheikhfatmaenit@yahoo.fr}.}
  \and
  Guy Chavent%
  \thanks[serena]{Équipe Serena.
    {\tt \{Guy.Chavent,Francois.Clement,Jean.Roberts\}@inria.fr}.}
  \and
  Fatma Cheikh%
  \thanksref{lamsin}%
  \thanksref{serena}
  \and
  François Clément%
  \thanksref{serena}
  \and
  Vincent Martin%
  \thanks[lmac]{LMAC, UTC, BP 20529, FR-60205 Compiègne, France.
    {\tt Vincent.Martin@utc.fr}.}
  \and
  Jean E. Roberts%
  \thanksref{serena}
}
\authorhead{%
  H. Ben Ameur, G. Chavent, F. Cheikh, F. Clément, V. Martin, \&
  J. E. Roberts}

\RRtitle{%
  Indicateurs du premier ordre pour l'estimation de fractures discrètes dans
  des milieux poreux}
\RRetitle{%
  First-Order Indicators for the Estimation of Discrete Fractures in Porous
  Media}
\titlehead{%
  First-Order Indicators for the Estimation of Discrete Fractures in Porous
  Media}

\RRresume{\input{resume}}
\RRabstract{\input{abstract}}

\RRmotcle{\input{motsclefs}}
\RRkeyword{\input{keywords}}

\RRprojet{Serena}

\RCParis

%% file: resume.tex

Les écoulements dans les milieux poreux peuvent être radicalement modifiés
par la présence de failles ou de barrières géologiques.
De telles fractures peuvent être modélisées comme des interfaces qui
interagissent avec la matrice environnante.
Nous proposons une nouvelle technique pour l'estimation de l'emplacement et
des propriétés hydrogéologiques d'un petit nombre de grandes fractures dans
un milieu poreux à partir de mesures distribuées de pression ou de flux
données.
À chaque itération, l'algorithme construit une courte liste de candidats par
comparaison d'{\em indicateurs de fracture}.
Ces indicateurs quantifient au premier ordre la décroissance d'une fonction
d'écart aux données;
ils sont peut coûteux à calculer.
Le meilleur candidat est ensuite isolé par minimisation de la fonction
objectif pour chaque candidat.
Guidée de façon optimale par la reproduction des données, l'approche a le
grand avantage de ne pas nécessiter de remaillage, ni de dérivation de
forme.
La stabilité de l'algorithme est montrée sur une série d'exemples numériques
représentatifs de situations typiques.

%% file: abstract.tex

Faults and geological barriers can drastically affect the flow patterns in
porous media.
Such fractures can be modeled as interfaces that interact with the
surrounding matrix.
We propose a new technique for the estimation of the location and
hydrogeological properties of a small number of large fractures in a porous
medium from given distributed pressure or flow data.
At each iteration, the algorithm builds a short list of candidates by
comparing {\em fracture indicators}.
These indicators quantify at the first order the decrease of a data misfit
function; they are cheap to compute.
Then, the best candidate is picked up by minimization of the objective
function for each candidate.
Optimally driven by the fit to the data, the approach has the great advantage
of not requiring remeshing, nor shape derivation.
The stability of the algorithm is shown on a series of numerical examples
representative of typical situations.

%% file: motsclefs.tex
paramétrisation adaptative,
problème inverse,
écoulement dans un milieu poreux,
milieu poreux fracturé,
faille et barrière

%% file: keywords.tex
adaptive parameterization,
inverse problem,
flow in porous media,
fractured porous media,
fault and barrier

%% file: body.tex

\myclearpage
\input{intro}

\myclearpage
\input{frac_direct}
\myclearpage
\input{frac_inverse}
\myclearpage
\input{frac_indicators}
\myclearpage
\input{frac_algorithm}
\myclearpage
\input{frac_numerics}

\myclearpage
\input{concl}

%% file: intro.tex

\section{Introduction}

The accurate simulation of flow in porous media is important for many
applications such as petroleum reservoir management, monitoring and clean-up
of underground pollutants, and planning for underground nuclear-waste
disposal.
In porous media, large fractures are often
present~\cite{dietrich:helmig:book:2009,adler:thovert:book:frac:2013}, and
they can modify drastically the flow pattern: the fracture, sometimes called
a fault in this case, can be very permeable and thus channel rapidly the
fluid from and to the surrounding rock.
Another typical regime arises when the fracture is almost impervious to the
flow: it acts as a geological barrier.
Accurately simulating the interaction of flow between porous media and
fractures is therefore important.
For this reason, it is desirable to determine as precisely as possible the
location of any important fracture and its hydrogeological signature, {\ie}
its impact on the flow.

The goal of this article is to present a new methodology for estimating the
location and the hydrogeological properties of a small number of large
fractures in a porous medium.
Geological experiments~\cite{LeGoc:deDreuzy:2010} indeed show that the fluid
tends to choose its pathway along only some few of the many existing
fractures.
So we do not address the case of a large number of small fractures that could
be treated using a double continuum
model~\cite{Barenblat-etal-JAMM1960,WarrenRoot-SPEJ1963}, but search only for
a limited number of fractures, the ones having the major impact on the flow.

Imaging techniques for assessing underground media such as seismic
imaging~\cite{biondi:book:2006,bimpas:anditis:uzunoglu:book:2010} can
give important information about the location of fractures, but they are not
able to tell to what extent if any of the fracture influences the flow.
These techniques should be seen as complementary to the present work: they
could for instance be used as {\apriori} information for the inversion
algorithm designed here.

The numerical flow model for the fractured medium must take into account the
strong heterogeneity and the very different spatial scales involved.
Reduced models or co-dimension~1 models treat the fractures as interfaces in
the porous medium, and allow flow in the fracture as well as fluid exchange
between the fracture and the rock matrix.
Such interface models have been studied extensively in the engineering
literature,
{\eg} see~\cite{BacaAL-IJNMF-1984,Fard-Firoozabadi-SPE-2003,%
  Fard-Durlofsky-SPE2004,RJBH-AWR-2006,Kim-Deo-AIChE-2006,%
  Hoteit-Firoozabadi-AWR-2008},
as well as in the mathematical literature,
{\eg} see~\cite{MJR-SISC-2005,ABH-M2AN-2009,QuarteroniEtal,%
  DangeloScotti2011,nf-vm-jr-as:fracnonconf:12,Faille-etal-2012,KR-2014,%
  Masson-etal-2014,Masson-etal-2015},
to name just a few.
We use here the model developed in~\cite{MJR-SISC-2005}, where flow in the
fracture as well as in the rock is governed by Darcy's law.
The model was first presented only for a fault
regime~\cite{AJRS-DDproc-1999}, but was extended to treat both fault and
barrier regimes
in~\cite{JMR-Porquerolles-2002,MJR-SISC-2005,FailleEtAl-2002}.
In this model, fractures lie along the edges of the mesh and can be easily
opened or closed by adjusting the fracture parameters on edges, which makes
it convenient for the purpose of fracture determination.
This model has been extended so that one may use non-matching
grids~\cite{nf-vm-jr-as:fracnonconf:12,TuncEtal2011}, or disconnect the
fracture mesh from that of the domain~\cite{MS-JMAA-2010,DangeloScotti2011}.
It has also been extended to treat Forchheimer flow in the
fracture~\cite{FRS-COMG-2008,KR-2014}, and multiphase
flow~\cite{RJBH-AWR-2006}, but these extensions are not considered in this
paper.

We suppose that the permeability is known outside the fracture, and
concentrate on the determination of the position and intensities of a few
fractures that explain pressure and/or flow data.
We have considered both distributed pressure and distributed velocity
measurements throughout the domain.
These data are used to build a least-squares objective function which
measures the misfit between data and the corresponding computed values.
Depending on the available data, the determination of the position and
intensity of fractures which minimize this objective function can be poorly
conditioned or under determined.
We count on the parsimonious introduction of fractures, which is inherent to
our inversion algorithm, to have a beneficial regularizing effect,
see~\cite{engl:regularization:book:1996,chavent:book:2009} for instance.

The key ingredient in our approach is the notion of
indicator~\cite[section 3.7]{chavent:book:2009}.
This idea was originally developed in~\cite{BAChJa:ref:coarse:ind:2002}, for
the estimation of hydraulic transmissivities, and was then extended for the
estimation of vector valued distributed parameters
in~\cite{BAClWCh:jiipp:08,BAChClW:ipse:11}.
In this paper, we adapt this notion to construct {\em fracture indicators},
which are used to {\em locate} fractures.
Given the current model of the fractured porous medium one associates any
additional {\em {\candidate} fracture\/}, {\ie} any set of new contiguous
edges not belonging to any existing fracture, with a
{\em fracture indicator}.
This indicator measures, up to the first order, the rate of decrease of the
objective function which would be achieved by adding the {\candidate}
fracture to the current model.
But, the important feature of the method is that, once both direct and
adjoint systems have been solved for the {\em current} model, the fracture
indicator associated with {\em any {\candidate} fracture} can be computed at
a very small cost.

Ideally, one should compute fracture indicators for the extremely long list
of {\em all possible {\candidate} fractures}, and retain only a short list of
those with the largest indicators.
However, even with this cheap indicator, such an exhaustive search is way out
of reach.
Hence, one has to resort to heuristic exploration strategies to reduce the
length of this long seminal list.
One possible
strategy~\cite{BAChJa:ref:coarse:ind:2002,BAClWCh:jiipp:08,BAChClW:ipse:11}
is to select {\candidate} parameterizations in a long list made of predefined
families.
Here, this would have required to consider {\em predefined families}
of {\candidate} fractures.
Due to the lack of {\apriori} information for the choice of such families, we
propose instead a new {\em constructive} strategy based on
{\em elementary {\candidate} fractures}.
In this approach, the fracture indicators are first computed
for all items of the {\em long list\/} of elementary {\candidate} fractures.
Then, proceeding by aggregation and extension, a {\em short list\/} of
{\candidate} fractures with {\em large indicators\/} is produced.

The algorithm is iterative.
Each iteration provides an additional fracture that produces a better fit to
the data.
The sketch of the algorithm goes as follows.
For a given set of data measures, an initial model with no fracture is
chosen.
Then, each iteration is made of two steps: an indicator step, which provides
a short list of {\candidate} fractures as explained above, and an
optimization step, which selects the winning fracture in the short list of
{\candidate} fractures and determines the associated intensities.
In this optimization step, pre-optimal fracture intensities are computed for
all models obtained by adding in turn each {\candidate} fracture of the
short list to the current model.
Note that only a limited number of such minimizations are necessary, as only
a small number of {\candidate} fractures have been retained in the short
list during the indicator step.
Note also that these optimizations are performed on a space of small
dimension, that is equal to the number of fractures (say less than~10) times
the number of physical parameters (1 or~2).
At the end of the iteration, the winning {\candidate} fracture is added
to the collection of already retained fractures.
This new fracture may be either an extension of a previously found fracture
or an entirely new fracture.
Then, the iterative process may resume with the new estimated model.

With this approach, the underlying mesh is fixed during the procedure, and
the ``opening'' of a fracture requires only changing its intensity parameter
to a nonzero value.
Hence, one avoids remeshing, most of the assembly process, and the
computation of shape derivatives and topological gradients, or similar
techniques to ``open the fracture'' (in other words, to compute a gradient
with respect to a fracture that does not exist yet).
The one drawback of the method is that the search is limited to fractures
located on the edges of a fixed mesh.
The fracture indicator algorithm is presented and tested numerically for a
two-dimensional model;
however, it remains valid for a three-dimensional model.

\bigskip

The paper is organized as follows:
the direct model for the numerical simulation of interactions between a
porous medium and a discrete fracture is given in detail in
Section~\ref{s:direct}.
The inverse problem for the estimation of fracture parameters is presented in
Section~\ref{s:inverse}.
The fracture indicators are introduced in Section~\ref{s:indicators}, and the
fracture indicator algorithm is detailed in Section~\ref{s:algorithm}.
Finally, Section~\ref{s:numerics} is devoted to numerical results.

%% file: frac_direct.tex

\section{The direct model}
\label{s:direct}

We consider the simulation of single-phase flow in a porous medium~$\Omega$
containing a fracture.
For simplicity of exposition we describe both the continuous and discrete
flow models in a 2-dimensional setting, though we stress that the methodology
is equally valid in a 3-dimensional setting.
Also for simplicity, in the description of the models only, we suppose that
$\Omega$ contains one fracture~$\gamfrac$.

\subsection{The model problem}
\label{ss:direct:continuous}

\begin{figure}[\figloc]
  \begin{center}
    \input{figure1}
    \caption{Notations associated with a fracture~$\gamfrac$ in the
      domain~$\Omega$.}
    \label{fig:dom:gamma}
  \end{center}
\end{figure}
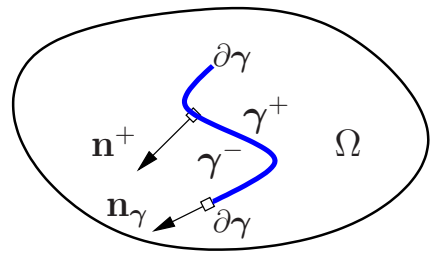

Let the porous medium $\Omega$ be a bounded domain in~$\matR^2$.
We suppose that the fracture $\gamfrac\subset\Omega$ is a regular,
non-self-intersecting, curve of finite length and of bounded curvature.
In this two-dimensional setting~$\pgamfrac$, the boundary of~$\gamfrac$,
consists of the endpoints of~$\gamfrac$ which we assume are distinct and do
not lie on the boundary of~$\Omega$.
Let~$\nnp$ be a unit normal vector field on $\gamfrac$, and let~$\nnm=-\nnp$.
We distinguish the two sides of the fracture, $\gamfracp$ and~$\gamfracm$,
such that $\nnp$ (resp.~$\nnm$) can be seen as the exterior normal
to $\Omega\setminus\gamfrac$ on~$\gamfracp$ (resp.~$\gamfracm$), see
Figure~\ref{fig:dom:gamma}.
For a sufficiently regular function~$\fhi$ in $\Omega\setminus\gamfrac$, we
consider two traces~$\fhip$ and~$\fhim$ on~$\gamfrac$, one for each side;
see~\cite{ABH-M2AN-2009} for a more rigorous definition.
We also introduce the tangential gradient~$\gradtau$ and the tangential
divergence~$\divtau$ operators along the fracture.

The fracture interface model~\cite{MJR-SISC-2005} may be written:
\begin{equation}
  \label{eq:fracmodel:direct:cont}
  \begin{array}{rcll}
    \vu & = & -\K \grad p &
    \mbox{in }\Omega \setminus \gamfrac,\\
    \div \, \vu & = & \ff &
    \mbox{in }\Omega \setminus \gamfrac,\\
    \ugf & = & -\alf \, \gradtau \pgf &
    \mbox{in } \gamfrac,\\
    \divtau \, \ugf & = & \fgf + \jump{\vu \cdot \nn} &
    \mbox{in } \gamfrac, \\
    - {\bet}\vup \cdot {\nnp} + \ppp & = & \pgf &
    \mbox{on } \gamfracp,\\
    - {\bet}\vum \cdot {\nnm} + \ppm & = & \pgf &
    \mbox{on } \gamfracm.
  \end{array}
\end{equation}
(With respect to the model~(4.1) in~\cite{MJR-SISC-2005}, the numerical
parameter~$\xi$ is taken equal to~$1$, since it has very little influence on
the flow in most cases, see~\cite{MJR-SISC-2005}.)

The first two equations in~\eqref{eq:fracmodel:direct:cont} are Darcy's law
and the law of mass conservation with a source term~$\ff$ in the domain
$\Omega\setminus\gamfrac$.
They relate the Darcy velocity~$\vu$ with the pressure~$p$ through the
permeability tensor~$\K$.
The next equation is Darcy's law, relating the tangential Darcy velocity~$\ugf$
to the fluid pressure in the fracture~$\pgf$.
The fourth equation is the law of mass conservation along the
fracture~$\gamfrac$, with a possible fracture source term~$\fgf$,
supplemented by an exchange term
$\jump{\vu\cdot\nn}=\vup\cdot\nnp+\vum\cdot\nnm$ accounting for fluid
exchange between~$\gamfrac$ and $\Omega\setminus\gamfrac$.
The last two equations that close the system represent an averaged Darcy law
across normal cross-sections of the physical fracture represented
by~$\gamfrac$.

The model~\eqref{eq:fracmodel:direct:cont} is supplemented with the boundary
conditions
\begin{equation}
  \label{prob:BC}
  \begin{array}{rclll}
    p & = & \pmeas & \mbox{on } \pOmega, & \\
    \ugf \cdot \ngf & = & 0 & \mbox{on }\pgamfrac, &
\end{array}
\end{equation}
where~$\pmeas$ is a given function, and~$\ngf$ is the unit normal vector
field on~$\pgamfrac$ pointing outwards from~$\gamfrac$ (and tangent
to~$\gamfrac$ at its endpoints), see Figure~\ref{fig:dom:gamma}.
The condition $\ugf\cdot\ngf=0$ is coherent with the hypothesis that the
endpoints~$\pgamfrac$ are inside the porous medium.

The coefficients~$\alf$ and~$\bet$ are hydrogeological parameters
characterizing the fracture:
$\alf$~represents an effective tangential permeability.
In the case of (conducting) geological faults, $\alf$~is positive and the
Darcy velocity is generally discontinuous across the fracture
($\jump{\vu\cdot\nn} \neq 0$).
In a two-dimension setting, $\alf$~is the product of a thickness and a
permeability divided by a length.
The coefficient~$\bet$ represents the inverse of an effective permeability
normal to the fracture.
So~$\bet$ models the resistivity of the fracture to flow {\em across\/} it.
In the case of geological barriers, $\bet$ is positive and the pressure is
generally discontinuous across the fracture ($\ppp \neq \ppm$).
In two dimensions, $\bet$~is the product of a thickness, a permeability, and
a length.

Our objective will be to determine both the location of the
fracture~$\gamfrac$ and the physical properties~$\alf$ and~$\bet$ of the
fracture.

We assume that the permeability tensor~$\K$ is symmetric and positive
definite, and bounded away from zero, and that the parameters~$\alf$
and~$\bet$ satisfy the hypotheses
\begin{equation}
  \label{eq:hyp:alpha:beta}
  \alf \geq \alf_{\min} > 0 \quad \mbox{ and } \quad \bet \geq 0.
\end{equation}
Under these hypotheses, the model~\eqref{eq:fracmodel:direct:cont}
and~\eqref{prob:BC} is well-posed, see~\cite{MJR-SISC-2005}.

\subsection{Discretization of the model problem}
\label{ss:direct:discret}

Let~$\Th$ be a finite element mesh of~$\Omega$, conforming with the
fracture~$\gamfrac$ in the sense that~$\gamfrac$ is a subset of the union of
the closures of the edges of the elements of~$\Th$.
We denote by $T,E,N$ respectively the elements ({\eg} triangles), edges, and
nodes of~$\Th$, by~$\ET$ the boundary edges of an element~$T$, by~$\TE$ the
elements having~$E$ as an edge, by~$\NE$ the endpoints of an edge~$E$, and
by~$\EN$ the edges having~$N$ as an endpoint.
Let~$\Eih$ be the set of internal edges of~$\Th$ and~$\Eeh$ the set of
boundary edges of~$\Th$, so that $\pOmega=\bigcup_{E\in\Eeh}\overline{E}$ and
$\Eeh\cap\Eih=\emptyset$.
We denote by $\Eh\egaldef\Eih\cup\Eeh$ the set of all edges of~$\Th$.

We say that the fracture~$\gamfrac$ is supported by the internal mesh edges
to mean that $\gamfrac\subset\bigcup_{E\in\Eih}\overline{E}$, and we assume
that the parameters~$\alf$ and~$\bet$ are discretized by piecewise constant
functions on these fracture edges.
To simplify the search for independent faults and barriers, we discretize the
fracture~$\gamfrac$ as
\begin{equation}
  \label{eq:frac:Fh}
  \currFh \egaldef (\currFalf, \currFbet),
\end{equation}
where we have separated {\em fault\/} edges in~$\currFalf$ from
{\em barrier\/} edges in~$\currFbet$,
\begin{equation}
  \label{eq:frac:sets:alf:bet}
  \currFalf \egaldef \{ E\in\Eih \suchthat \alfE > 0 \}
  \quad \mbox{and} \quad
  \currFbet \egaldef \{ E\in\Eih \suchthat \betE > 0 \}.
\end{equation}
We denote the collections of {\em fault} (or tangential) and {\em barrier}
(or normal) parameters by
\begin{equation}
  \label{eq:def:alf:bet}
  \alf \egaldef (\alfE)_{E \in \currFalf}
  \quad \mbox{and} \quad
  \bet \egaldef (\betE)_{E \in \currFbet}.
\end{equation}
The sets of fault and barrier edges need not be disjoint.
Let~$\currNalf$ be the set of all nodes of~$\currFalf$.

In order to simplify the presentation, we assume that the mesh is made up of
square elements of size~$h$, that the permeability tensor in~$\Omega$ is
constant and scalar $\K=K\mathbf{I}$ (with $K>0$), and that the
model~\eqref{eq:fracmodel:direct:cont} and~\eqref{prob:BC} is discretized
using a cell-centered finite volume method, see~\cite{FailleEtAl-2002}.
The discretization in the more general case of lowest order Raviart-Thomas
mixed hybrid finite elements with general meshes is presented
in~\cite{these:fatma:2016}.
These two discretizations in rectangular meshes are equivalent up to an
approximate quadrature formula, see~\cite{charob:91}.

Let~$\PT$, $\PE$, and~$\PN$ be the pressure unknowns on an element $T\in\Th$,
an edge $E\in\Eh$ and a node $N\in\currNalf$.
Let~$\UTE$ be the flow rate of fluid leaving an element $T\in\Th$ through an
edge $E\in\ET$, let~$\UEN$ be the flow rate of fluid leaving an edge
$E\in\currFalf$ through a node $N\in\NE$.
Using the notations $\kh\egaldef K^{-1}h$, $\ffT\egaldef\int_T\ff$,
$\ffE\egaldef\int_E\fgf$, and $\PmE\egaldef\frac{1}{|E|}\int_E\pmeas$, the
finite volume formulation of~\eqref{eq:fracmodel:direct:cont}
and~\eqref{prob:BC} reads
\begin{equation}
  \label{eq:fracmodel:direct:discr:Omega}
  \begin{array}{rcll}
    \dps \sum_{E\in\ET} \UTE & = & \ffT, &
    \mbox{for all } T\in\Th, \\
    \dps \frac{\kh}{2} \, \UTE & = & \PT - \PE, &
    \mbox{for all } E\notin\currFbet,\, T\in\TE, \\
    \dps \left( \frac{\kh}{2} + \betE \right) \, \UTE & = & \PT - \PE, &
    \mbox{for all } E\in\currFbet,\, T\in\TE, \\
    \dps \PE & = & \PmE, & \mbox{for all } E\in\Eeh, \\
    \dps - \sum_{T\in\TE} \UTE & = & 0, &
    \mbox{for all } E\notin\currFalf,\, E\in\Eih,
  \end{array}
\end{equation}
and
\begin{equation}
  \label{eq:fracmodel:direct:discr:gamma}
  \begin{array}{rcll}
    \dps - \sum_{T\in\TE} \UTE + \sum_{N\in\NE} \UEN & = & \ffE, &
    \mbox{for all } E\in\currFalf, \\
    \dps \frac{h}{2\alfE} \, \UEN & = & \PE - \PN, &
    \mbox{for all } E\in\currFalf,\, N\in\NE, \\
    \dps \sum_{\substack{E\in\currFalf \\ E\in\EN}} \UEN & = & 0, &
    \mbox{for all } N\in\currNalf,
  \end{array}
\end{equation}
where the last equation already contains, for endpoints of the fault~$\Falf$,
the no-flow condition $\ugf\cdot\ngf=0$ of~\eqref{prob:BC}.

In Section~\ref{s:indicators}, the construction of fracture indicators will
require to pass to the limit in the above equations when $\alfE\tendstoz$ on
some edges of~$\currFalf$.
But a quick look at the second equation
in~\eqref{eq:fracmodel:direct:discr:gamma} let us foresee some problem in
doing so.
So we replace on all fault edges the flow rate unknown~$\UEN$ by a new
unknown~$\VEN$ defined by (remember that $\alfE>0$ in this section)
\begin{equation}
  \label{eq:chg:var:UEN}
  \alfE \VEN = \UEN,
  \quad
  \mbox{for all } E\in\currFalf,\, N\in\NE.
\end{equation}
With these unknowns, \eqref{eq:fracmodel:direct:discr:gamma} is replaced by
\begin{equation}
  \label{eq:fracmodel:direct:discr:gamma:chgvar}
  \begin{array}{rcll}
    \dps - \sum_{T\in\TE} \UTE + \alfE \sum_{N\in\NE} \VEN & = & \ffE, &
    \mbox{for all } E\in\currFalf, \\
    \dps \frac{h}{2} \, \VEN & = & \PE - \PN, &
    \mbox{for all } E\in\currFalf,\, N\in\NE,\\
    \dps \sum_{\substack{E\in\currFalf \\ E\in\EN}} \alfE \, \VEN & = & 0, &
    \mbox{for all } N\in\currNalf.
  \end{array}
\end{equation}
The second equation in~\eqref{eq:fracmodel:direct:discr:gamma:chgvar} shows
that the new unknown~$\VEN$ is the (discrete) pressure gradient between the
middle of edge~$E$ and its end node~$N$.

\bigskip

For the sake of simplicity, we suppose, in the sequel, the fracture source
term to be zero ($\ffE=0$
in~\eqref{eq:fracmodel:direct:discr:gamma:chgvar}).
To sum up, the (discrete) {\em direct model\/} for the simulation of flow in
a porous domain~$\Omega$ containing a discrete fracture $(\Falf,\Fbet)$ with
fault and barrier parameters $(\alf,\bet)$ is made up of
equations~\eqref{eq:fracmodel:direct:discr:Omega}
and~\eqref{eq:fracmodel:direct:discr:gamma:chgvar}.
Thus the {\em direct problem\/} associates to the fracture, {\ie} the pair
$\Fh=(\Falf,\Fbet)$ together with the fracture parameters, the solution of
the direct model.

%% file: figure1.tex
\begin{picture}(0,0)%
\includegraphics{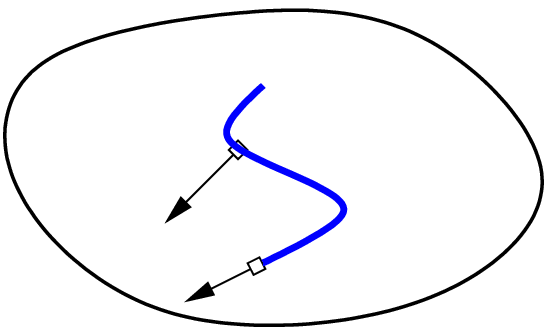}%
\end{picture}%
\setlength{\unitlength}{4144sp}%
\begingroup\makeatletter\ifx\SetFigFont\undefined%
\gdef\SetFigFont#1#2#3#4#5{%
  \reset@font\fontsize{#1}{#2pt}%
  \fontfamily{#3}\fontseries{#4}\fontshape{#5}%
  \selectfont}%
\fi\endgroup%
\begin{picture}(2501,1485)(868,-1089)
\put(2071, 29){\makebox(0,0)[lb]{\smash{{\SetFigFont{12}{14.4}{\familydefault}{\mddefault}{\updefault}{$\pgamfrac$}%
}}}}
\put(2791,-511){\makebox(0,0)[lb]{\smash{{\SetFigFont{14}{16.8}{\familydefault}{\mddefault}{\updefault}{$\Omega$}%
}}}}
\put(2251,-331){\makebox(0,0)[lb]{\smash{{\SetFigFont{14}{16.8}{\familydefault}{\mddefault}{\updefault}{$\gamfracp$}%
}}}}
\put(1441,-871){\makebox(0,0)[lb]{\smash{{\SetFigFont{14}{16.8}{\familydefault}{\mddefault}{\updefault}{$\ngf$}%
}}}}
\put(1351,-511){\makebox(0,0)[lb]{\smash{{\SetFigFont{14}{16.8}{\familydefault}{\mddefault}{\updefault}{$\nnp$}%
}}}}
\put(2071,-961){\makebox(0,0)[lb]{\smash{{\SetFigFont{12}{14.4}{\familydefault}{\mddefault}{\updefault}{$\pgamfrac$}%
}}}}
\put(1981,-601){\makebox(0,0)[lb]{\smash{{\SetFigFont{14}{16.8}{\familydefault}{\mddefault}{\updefault}{$\gamfracm$}%
}}}}
\end{picture}%

%% file: frac_inverse.tex

\section{The inverse problem}
\label{s:inverse}

Focusing now on the {\em inverse problem}, we assume that pressure and/or
flow measurements are available at some points of the domain~$\Omega$, and
consider the problem of estimating the location and parameters of the
corresponding fracture.
For the sake of simplicity, we limit ourselves in this section to the simple
case where a pressure measurement~$\PmT$ is available for each element~$T$
in~$\Th$.
More realistic measurements will be considered for numerical experiments in
Section~\ref{s:numerics}.

\subsection{The least-squares formulation}
\label{ss:least:squares}

In the framework of an iterative process for solving the inverse problem, let
the {\em current fracture\/} be described by its edges
$\Fh=(\currFalf,\currFbet)$ and parameters $(\alf,\bet)$.
We associate the current fracture with the least-square data misfit
\begin{equation}
  \label{eq:def:cost:fun}
  \JcFh (\alf, \bet) \egaldef \half \sum_{T\in\Th} (\PT - \PmT)^2,
\end{equation}
where~$\PT$ is given by the
{\em current model\/}~\eqref{eq:fracmodel:direct:discr:Omega}
and~\eqref{eq:fracmodel:direct:discr:gamma:chgvar} set on the current
fracture~$\Fh$.

If this misfit is above the uncertainty level of the data measurements, one
is led to consider adding new {\em {\candidate} edges\/}
$\tentFh=(\tentFalf,\tentFbet)$ to the current fracture~$\currFh$.
We suppose that current and {\candidate} fractures have no common edge,
\begin{equation}
  \label{eq:FcapDF=0}
  \currFalf \cap \tentFalf = \emptyset
  \quad\mbox{and}\quad
  \currFbet \cap \tentFbet = \emptyset.
\end{equation}
Let~$\tentNalf$ be the set of new nodes of the {\candidate}
fault~$\tentFalf$, so that
\begin{equation}
  \label{eq:NcapDN=0}
  \currNalf \cap \tentNalf = \emptyset.
\end{equation}

We use {\em nominal parameters\/}~$\alfbar$ and~$\betbar$ to describe typical
parameter profiles expected along the new {\candidate} fault and barrier
edges.
In the absence of such information, nominal parameters are chosen to be
constant.
We model the introduction of the {\candidate} fracture by new
{\em intensity parameters\/}~$\epsalf$ and~$\epsbet$, which let the
fracture parameters~$\alf$ and~$\bet$ grow from zero to their nominal value
and beyond,
\begin{eqnarray}
  \label{param:alfa:vanish:prop}
  \alfE = \epsalf \, \alfbarE, & &
  \mbox{for all } E\in\tentFalf,\, \epsalf\geq 0, \\
  \label{param:beta:vanish:prop}
  \betE =\epsbet \, \betbarE, & &
  \mbox{for all } E\in\tentFbet,\, \epsbet\geq 0.
\end{eqnarray}

Because of~\eqref{param:alfa:vanish:prop}, the last equation
of~\eqref{eq:fracmodel:direct:discr:gamma:chgvar} can be divided by~$\epsalf$
on~$\tentNalf$, thus the model~\eqref{eq:fracmodel:direct:discr:Omega}
and~\eqref{eq:fracmodel:direct:discr:gamma:chgvar} set on both current and
{\candidate} fractures~$\currFh$ and~$\tentFh$ simplifies to the
{\em {\tentative} model}
\begin{equation}
  \label{eq:fracmodel:direct:discr:tentative}
  \begin{array}{rcll}
    \dps \sum_{E\in\ET} \UTE & = & \ffT, &
    \mbox{for all } T\in\Th, \\
    \dps \frac{\kh}{2} \, \UTE & = & \PT - \PE, &
    \mbox{for all } E\notin\fullFbet,\, T\in\TE, \\
    \dps \left( \frac{\kh}{2} + \betE \right) \, \UTE & = & \PT - \PE, &
    \mbox{for all } E\in\fullFbet,\, T\in\TE, \\
    \dps \PE & = & \PmE, &
    \mbox{for all } E\in\Eeh, \\
    \dps - \sum_{T\in\TE} \UTE & = & 0, &
    \mbox{for all } E \notin\fullFalf,\, E\in\Eih, \\
    \dps - \sum_{T\in\TE} \UTE + \alfE \sum_{N\in\NE} \VEN & = &0, &
    \mbox{for all } E\in\fullFalf, \\
    \dps \frac{h}{2} \, \VEN & = & \PE - \PN, &
    \mbox{for all } E\in\fullFalf,\, N\in\NE, \\
    \dps \sum_{\substack{E\in\fullFalf \\ E\in\EN}} \alfE \, \VEN & = & 0, &
    \mbox{for all } N\in\Nalf,\\
    \dps \sum_{\substack{E\in\tentFalf \\ E\in\EN}} \alfbarE \, \VEN & = & 0, &
    \mbox{for all } N\in\tentNalf.
  \end{array}
\end{equation}

\goodbreak

When $\epsalf=\epsbet=0$, one easily checks that the {\tentative}
model~\eqref{eq:fracmodel:direct:discr:tentative} reduces to the current
model~\eqref{eq:fracmodel:direct:discr:Omega}
and~\eqref{eq:fracmodel:direct:discr:gamma:chgvar} complemented with the
{\em limit fault model\/} set on the {\candidate} fault,
\begin{equation}
  \label{eq:fracmodel:direct:discr:gamma:chgvar:alfa=0}
  \begin{array}{rcll}
    \dps \frac{h}{2} \, \VEN & = & \PE - \PN &
    \mbox{for all } E\in\tentFalf,\, N\in\NE, \\
    \dps \sum_{\substack{E\in\tentFalf \\ E\in\EN}} \alfbarE \, \VEN & = & 0 &
    \mbox{for all } N\in\tentNalf.
  \end{array}
\end{equation}
Hence, we see that, for any $\epsalf,\epsbet\geq 0$, the {\tentative}
model~\eqref{eq:fracmodel:direct:discr:tentative} determines
uniquely~$\PT$, $\PE$, and~$\UTE$ on all elements and edges, and~$\PN$
and~$\VEN$ on all nodes and edges of the current {\em and\/} {\candidate}
faults.
So, we can extend the definition~\eqref{eq:def:cost:fun} of the least-squares
data misfit to a function of the intensity parameters~$\epsalf,\epsbet\geq 0$
on the {\candidate} fracture by
\begin{equation}
  \label{eq:def:cost:fun:eps}
  \JtFh (\epsalf, \epsbet) \egaldef \half \sum_{T\in\Th} (\PT - \PmT)^2,
\end{equation}
where~$\PT$ is now given by the {\tentative}
model~\eqref{param:alfa:vanish:prop}, \eqref{param:beta:vanish:prop},
and~\eqref{eq:fracmodel:direct:discr:tentative}.
The purpose of the indicator algorithm developed in the next sections will
be to select, in a computationally efficient way, the {\candidate}
fractures~$\tentFh$ for which there exists $\epsalf,\epsbet>0$ such that
\begin{equation}
  \label{tentative:decrease:objective:function}
  \JtFh (\epsalf, \epsbet) < \JtFh (0, 0)
  = \JcFh (\alf, \bet),
\end{equation}
with a preferably large decrease of the data misfit.

\begin{Remark}
  \label{rem:alpha-to-zero:general}
  Replacing the hypothesis $\epsalf=\epsbet=0$
  in~\eqref{param:alfa:vanish:prop} and~\eqref{param:beta:vanish:prop} by the
  weaker assumption that $\alf,\bet\tendstoz$ for all {\candidate} edges
  of~$\tentFh$, makes the last equation
  in~\eqref{eq:fracmodel:direct:discr:gamma:chgvar} disappear on nodes
  of~$\tentNalf$: unknowns $(\PN)_{N\in\tentNalf}$ and
  $(\VEN)_{E\in\tentFalf,N\in\NE}$ become undetermined for the limit
  equations, but all other unknowns remain uniquely determined.
  
  Actually, in Section~\ref{s:indicators}, the key ingredient for the
  definition of fracture indicators is the fading of fracture parameters on
  the {\candidate} fracture.
  And, choosing a specific direction in the space of fracture edge parameters
  through~\eqref{param:alfa:vanish:prop} and~\eqref{param:beta:vanish:prop},
  makes the solution of the {\tentative}
  model~\eqref{eq:fracmodel:direct:discr:tentative} behave continuously
  for vanishing intensity parameters.
\end{Remark}

\subsection{Gradient with respect to the intensity parameters}
\label{ss:gradient}

We compute now the gradient of the objective function with respect
to~$\epsalf$ and~$\epsbet$ by the adjoint state method, for instance
see~\cite[Section 2.4]{chavent:book:2009}.
For the sake of simplicity, the Dirichlet boundary condition (fourth
equation in~\eqref{eq:fracmodel:direct:discr:tentative}) is not written
explicitly, it is understood that whenever~$E$ is a boundary edge, $\PE$
means~$\PmE$ (and hence its differential vanishes).
With this convention, the vector of edge pressure unknowns becomes
$(\PE)_{E\in\Eih}$.

The Lagrangian associated with the minimization of the objective
function~\eqref{eq:def:cost:fun:eps} under the constraint of the state
equations~\eqref{eq:fracmodel:direct:discr:tentative} for the current and
{\candidate} fractures~$\fullFh$ is then
\begin{eqnarray}
  \label{eq:def:lagr}
  \Lag(X, \Lambda; \epsalf, \epsbet) & \egaldef &
  \half \sum_{T\in\Th} (\PT - \PmT)^2 \\
  \nonumber &&
  + \sum_{T\in\Th}
  \left( -\sum_{E\in\ET} \UTE + \ffT \right) \, \LAT \\
  \nonumber &&
  + \sum_{E\notin\fullFbet} \sum_{T\in\TE}
  \left( \frac{\kh}{2} \, \UTE - \PT + \PE \right) \, \MUTE \\
  \nonumber &&
  + \sum_{E\in\fullFbet} \sum_{T\in\TE} \left(
    \left( \frac{\kh}{2} + \betE \right) \, \UTE - \PT + \PE
  \right) \, \MUTE \\
  \nonumber &&
  + \sum_{\substack{E\notin\fullFalf \\ E\in\Eih}} \sum_{T\in\TE} \UTE \, \LAE \\
  \nonumber &&
  + \sum_{E\in\fullFalf} \left(
    \sum_{T\in\TE} \UTE - \alfE \sum_{N\in\NE} \VEN \right) \, \LAE \\
  \nonumber &&
  + \sum_{E\in\fullFalf} \sum_{N\in\NE}
  \left( \frac{h}{2} \, \VEN - \PE + \PN \right) \, \MUEN \\
  \nonumber &&
  + \sum_{N\in\Nalf}
  \sum_{\substack{E\in\fullFalf \\ E\in\EN}} \alfE \, \VEN \, \LAN \\
  \nonumber &&
  + \sum_{N\in\tentNalf}
  \sum_{\substack{E\in\tentFalf \\ E\in\EN}} \alfbarE \, \VEN \, \LAbarN,
\end{eqnarray}
where~$\epsalf$ and~$\epsbet$ are the intensity parameters,
$X\egaldef(\Pt,\Pe,\Pn,\Ute,\Ven)$ is the state vector, and
$\Lambda\egaldef(\LAt,\LAe,\LAn,\LAbarn,\MUte,\MUen)$ is the adjoint state
vector, with the notations
\begin{equation}
  \label{eq:def:unknowns}
  \begin{array}{ll}
    \Pt \egaldef (\PT)_{T\in\Th}, &
    \Ute \egaldef (\UTE)_{T\in\Th,\, E\in\ET}, \\
    \Pe \egaldef (\PE)_{E\in\Eih}, &
    \Ven \egaldef (\VEN)_{E\in\fullFalf,\, N\in\NE}, \\
    \Pn \egaldef (\PN)_{N\in\fullNalf}, \\
    \LAt \egaldef (\LAT)_{T\in\Th}, &
    \MUte \egaldef (\MUTE)_{T\in\Th,\, E\in\ET}, \\
    \LAe \egaldef (\LAE)_{E\in\Eih}, &
    \MUen \egaldef (\MUEN)_{E\in\fullFalf,\, N\in\NE}, \\
    \LAn \egaldef (\LAN)_{N\in\currNalf}, &
    \LAbarn \egaldef (\LAbarN)_{N\in\tentNalf}.
  \end{array}
\end{equation}

Let~$\Xeps$ denote the state vector corresponding to the solution of the
{\tentative} model~\eqref{eq:fracmodel:direct:discr:tentative} with the
fracture parameters~$\alf$ and~$\bet$ defined
by~\eqref{param:alfa:vanish:prop} and~\eqref{param:beta:vanish:prop} on the
{\candidate} fracture~$\tentFalf$.
Let~$\Lambda$ be a given adjoint state vector.
Then, for all $\epsalf,\epsbet\geq 0$, we have the identity
\begin{equation}
  \label{eq:lagrangian:equals:cost}
  \JtFh (\epsalf, \epsbet) = \Lag (\Xeps, \Lambda; \epsalf, \epsbet).
\end{equation}
Differentiation with respect to parameters~$\epsalf$ and~$\epsbet$ gives
(remember that~$\Lambda$ is given)
\begin{equation}
  \label{cost:fun:differentiation:step1}
  \Diff\JtFh (\epsalf, \epsbet) =
  \Diff \Lag (\Xeps, \Lambda; \epsalf, \epsbet),
\end{equation}
hence, for all $\epsalf,\epsbet\geq 0$,
\begin{equation}
  \label{cost:fun:differentiation:step2}
  \Diff\JtFh (\epsalf, \epsbet) =
  \DLagX (\Xeps, \Lambda; \epsalf, \epsbet) \, \Diff X
  + \DLageps (\Xeps, \Lambda; \epsalf, \epsbet) \, \Diff (\epsalf, \epsbet).
\end{equation}
This identity holds for any adjoint state vector~$\Lambda$.
Thus, nothing can prevent us from choosing for $\Lambda$ precisely the unique
solution~$\Lambdaeps$ of the {\em adjoint equations\/} such that
\begin{equation}
  \label{eq:fracmodel:adjoint:discrete:lagrangian}
  \DLagX (\Xeps, \Lambda; \epsalf, \epsbet) \, \Diff X = 0,
  \quad \forall \Diff X.
\end{equation}
Then, \eqref{cost:fun:differentiation:step2} simplifies to the desired
{\em gradient equations}, for all $\epsalf,\epsbet\geq 0$,
\begin{equation}
  \label{cost:fun:differentiation:step3}
  \Diff\JtFh (\epsalf, \epsbet) =
  \DLageps (\Xeps, \Lambdaeps; \epsalf, \epsbet) \, \Diff (\epsalf, \epsbet),
\end{equation}
and the gradient of~$\JtFh$ is simply obtained by differentiation of the
Lagrangian~\eqref{eq:def:lagr} with respect to the intensity
parameters~$\epsalf$ and~$\epsbet$.
This gives, for all $\epsalf,\epsbet\geq 0$,
\begin{equation}
  \label{eq:grad}
  \begin{array}{rcl}
    \dps \DJtFhepsalf (\epsalf, \epsbet) & = &
    \dps - \sum_{E\in\tentFalf} \alfbarE \sum_{N\in\NE} \VEN \, \LAE
    + \sum_{E\in\tentFalf} \alfbarE
    \sum_{\substack{N\in\NE \\ N\in\currNalf}} \VEN \, \LAN, \\
    \dps \DJtFhepsbet (\epsalf, \epsbet) & = &
    \dps \sum_{E\in\tentFbet} \betbarE \sum_{T\in\TE} \UTE \, \MUTE,
  \end{array}
\end{equation}
where~$\VEN$ and~$\UTE$ are solutions of the state
equations~\eqref{eq:fracmodel:direct:discr:tentative}, and~$\LAE$, $\LAN$,
and $\MUTE$ are solution of the adjoint
equations~\eqref{eq:fracmodel:adjoint:discrete:lagrangian}, which we explicit
now.

Differentiation of the Lagrangian~\eqref{eq:def:lagr} with respect to the
state vector~$X$ gives
\begin{eqnarray}
  \label{eq:lagr:diff}
  \DLagX (X, \Lambda; \epsalf, \epsbet) \, \Diff X & = &
  \sum_{T\in\Th} (\PT - \PmT) \, \Diff \PT \\
  \nonumber &&
  - \sum_{T\in\Th} \sum_{E\in\ET} \Diff \UTE \, \LAT \\
  \nonumber &&
  + \sum_{E\notin\fullFbet} \sum_{T\in\TE}
  \left( \frac{\kh}{2} \, \Diff \UTE - \Diff \PT + \Diff \PE \right)
  \, \MUTE \\
  \nonumber &&
  + \sum_{E\in\fullFbet} \sum_{T\in\TE} \left(
    \left( \frac{\kh}{2} + \betE \right) \, \Diff \UTE
    - \Diff \PT + \Diff \PE
  \right) \, \MUTE \\
  \nonumber &&
  + \sum_{\substack{E\notin\fullFalf \\ E\in\Eih}}
  \sum_{T\in\TE} \Diff \UTE \, \LAE \\
  \nonumber &&
  + \sum_{E\in\fullFalf} \left(
    \sum_{T\in\TE} \Diff \UTE - \alfE \sum_{N\in\NE} \Diff \VEN
  \right) \, \LAE \\
  \nonumber &&
  + \sum_{E\in\fullFalf} \sum_{N\in\NE} \left(
    \frac{h}{2} \, \Diff \VEN - \Diff \PE + \Diff \PN
  \right) \, \MUEN \\
  \nonumber &&
  + \sum_{N\in\Nalf} \sum_{\substack{E\in\fullFalf \\ E\in\EN}}
  \alfE \, \Diff \VEN \, \LAN \\
  \nonumber &&
  + \sum_{N\in\tentNalf} \sum_{\substack{E\in\tentFalf \\ E\in\EN}}
  \alfbarE \, \Diff \VEN \, \LAbarN,
\end{eqnarray}
where, according to our convention, $\Diff\PE=0$ on boundary edges.
The adjoint state equation~\eqref{eq:fracmodel:adjoint:discrete:lagrangian}
is satisfied if and only if all coefficients of~$\Diff\PT$, $\Diff\PE$,
$\Diff\PN$, $\Diff\UTE$, and~$\Diff\VEN$ are equal to zero.
This gives the following {\em {\tentative} adjoint model\/} for the
determination of the adjoint variables~$\LAT$, $\LAE$, $\LAN$, $\LAbarN$,
$\MUTE$, and~$\MUEN$,
\begin{equation}
  \label{eq:fracmodel:adj:discr:Omega}
  \begin{array}{rcll}
    \dps \sum_{E\in\ET} \MUTE & = & \PT - \PmT, &
    \mbox{for all } T\in\Th, \\
    \dps \frac{\kh}{2} \, \MUTE & = & \LAT - \LAE, &
    \mbox{for all } E\notin\fullFbet,\, T\in\TE, \\
    \dps \left( \frac{\kh}{2} + \betE \right) \, \MUTE & = & \LAT - \LAE, &
    \mbox{for all } E\in\fullFbet,\, T\in\TE, \\
    \dps \LAE & = & 0, &
    \mbox{for all } E\in\Eeh, \\
    \dps - \sum_{T\in\TE} \MUTE & = & 0, &
    \mbox{for all } E\notin\fullFalf,\, E\in\Eih, \\
  \end{array}
\end{equation}
\begin{equation}
  \label{eq:fracmodel:adj:discr:gamma}
  \begin{array}{rcll}
    \hspace*{-4em}
    \dps - \sum_{T\in\TE} \MUTE + \sum_{N\in\NE} \MUEN & = & 0, &
    \mbox{for all } E\in\fullFalf, \\
    \dps \frac{h}{2} \, \MUEN & = & \alfE \, (\LAE - \LAN), &
    \mbox{for all } E\in\fullFalf,\, N\in\Nalf\cap\NE, \\
    \dps \frac{h}{2} \, \MUEN & = & \alfE \, \LAE - \alfbarE \, \LAbarN, &
    \mbox{for all } E\in\tentFalf,\, N\in\tentNalf\cap\NE, \\
    \dps \sum_{\substack{E\in\fullFalf \\ E\in\EN}} \MUEN & = & 0, &
    \mbox{for all } N\in\fullNalf.
  \end{array}
\end{equation}

We rewrite now these adjoint equations in a form similar to that of
the {\tentative} model~\eqref{eq:fracmodel:direct:discr:tentative}.
First, define~$\LAN$ for~$N\in\tentNalf$ and~$\NUEN$ for all
current and {\candidate} fault edges by (compare with the seventh equation
in~\eqref{eq:fracmodel:direct:discr:tentative})
\begin{equation}
  \label{eq:discrete:adjoint:gamma:newvar}
  \begin{array}{rcll}
    \LAN & = &
    \dps \frac{\sum_{\substack{E\in\tentFalf \\ E\in\EN}} \alfbarE \, \LAE}%
    {\sum_{\substack{E\in\tentFalf \\ E\in\EN}} \alfbarE}, &
    \mbox{for all } N\in\tentNalf, \\
    \dps \frac{h}{2} \, \NUEN & = & \LAE - \LAN, &
    \mbox{for all } E\in\fullFalf,\, N\in\NE.
  \end{array}
\end{equation}
Then, from the last two equations in~\eqref{eq:fracmodel:adj:discr:gamma},
one checks that the adjoint variables~$\LAbarN$ and~$\MUEN$ satisfy (compare
with~\eqref{eq:chg:var:UEN})
\begin{equation}
  \label{eq:discrete:adjoint:gamma:chgtvar}
  \begin{array}{rcll}
    \LAbarN & = & \epsalf \, \LAN &
    \mbox{for all } N\in\tentNalf, \\
    \dps \MUEN & = & \alfE \, \NUEN &
    \mbox{for all } E\in\fullFalf,\, N\in\NE.
  \end{array}
\end{equation}
With the definition~\eqref{eq:discrete:adjoint:gamma:newvar},
\eqref{eq:fracmodel:adj:discr:gamma} becomes (compare with the last four
equations in~\eqref{eq:fracmodel:direct:discr:tentative})
\begin{equation}
  \label{eq:fracmodel:adj:discr:gamma:chgtvar}
  \begin{array}{rcll}
    \dps - \sum_{T\in\TE} \MUTE + \alfE \sum_{N\in\NE} \NUEN & = & 0, &
    \mbox{for all } E\in\fullFalf, \\
    \dps \frac{h}{2} \, \NUEN & = & \LAE - \LAN, &
    \mbox{for all } E\in\fullFalf,\, N\in\NE, \\
    \dps \sum_{\substack{E\in\fullFalf \\ E\in\EN}} \alfE \, \NUEN & = & 0, &
    \mbox{for all } N\in\currNalf, \\
    \dps \sum_{\substack{E\in\tentFalf \\ E\in\EN}} \alfbarE \, \NUEN & = & 0, &
    \mbox{for all } N\in\tentNalf.
  \end{array}
\end{equation}

\begin{Remark}
  \label{rem:limit:adjoint:model}
  Unsurprisingly, the tentative adjoint
  model~\eqref{eq:fracmodel:adj:discr:Omega}
  and~\eqref{eq:fracmodel:adj:discr:gamma},
  or~\eqref{eq:fracmodel:adj:discr:Omega}
  and~\eqref{eq:fracmodel:adj:discr:gamma:chgtvar}, makes sense for
  $\epsalf,\epsbet\geq 0$, as the tentative direct
  model~\eqref{eq:fracmodel:direct:discr:tentative} does.
  When $\epsalf=\epsbet=0$, the adjoint model reduces to the\/ {\em current}
  adjoint model, that will be defined
  in~\eqref{eq:fracmodel:adj:discr:current}, complemented with a limit fault
  adjoint model similar
  to~\eqref{eq:fracmodel:direct:discr:gamma:chgvar:alfa=0} and whose solution
  is given by~\eqref{eq:discrete:adjoint:gamma:newvar} on the {\candidate}
  fault.
\end{Remark}

The introduction of~$\LAN$ and~$\NUEN$ gives a new formula for
$\partial\JtFh/\partial\epsalf$.
Indeed, the multiplication of the last equation
in~\eqref{eq:fracmodel:direct:discr:tentative} by~$\LAN$ and the summation
over all nodes of~$\tentNalf$ give
\begin{equation}
  \label{eq:grad:missing:term}
  \sum_{N\in\tentNalf} \sum_{\substack{E\in\tentFalf \\ E\in\EN}}
  \alfbarE \, \VEN \, \LAN = 0.
\end{equation}
Then, by adding this zero term to the right-hand side of the first gradient
equation in~\eqref{eq:grad}, one obtains a new formula for the derivative
of~$\JtFh$ with respect to~$\epsalf$, for all $\epsalf,\epsbet\geq 0$,
\begin{equation}
  \label{eq:grad:oldadjvar}
  \DJtFhepsalf (\epsalf, \epsbet) =
  - \frac{h}{2} \,
  \sum_{E\in\tentFalf} \alfbarE \sum_{N\in\NE} \VEN \, \NUEN.
\end{equation}
This formula is more elegant than the one in~\eqref{eq:grad}.
It requires to solve the adjoint limit fault
model~\eqref{eq:discrete:adjoint:gamma:newvar} for~$\LAN$ and~$\NUEN$, but
only on the {\candidate} fault.

\begin{Remark}
  \label{rem:measures}
  When we do not have a pressure measurement for all elements of the mesh,
  the right-hand side of the first equation
  in~\eqref{eq:fracmodel:adj:discr:Omega} vanishes for elements with no
  measurement.
  In the same way, when a Darcy velocity measurement~$\UmTE$ is provided for
  element~$T\in\Th$ and for edge~$E\in\ET$, the additional source term
  $\UTE-\UmTE$ appears on the right-hand side in either the second or the
  third equation in~\eqref{eq:fracmodel:adj:discr:Omega}, following that~$E$
  is a barrier edge or not.
\end{Remark}

\bigskip

To conclude this section, we recall that the {\tentative} direct and adjoint
systems~\eqref{eq:fracmodel:direct:discr:tentative},
\eqref{eq:fracmodel:adj:discr:Omega},
and~\eqref{eq:fracmodel:adj:discr:gamma:chgtvar} define uniquely the direct
and adjoint variables~$X$ and~$\Lambda$ for all $\epsalf,\epsbet\geq 0$,
{\em including $\PN$, $\LAN$, $\VEN$, and~$\NUEN$ on the {\candidate} fault}.
Hence, the gradient~\eqref{eq:grad} or~\eqref{eq:grad:oldadjvar} of the
objective function~$\JtFh(\epsalf,\epsbet)$ is also defined when the
intensity parameters~$\epsalf$ and~$\epsbet$ vanish on the {\candidate}
fracture.

%% file: frac_indicators.tex

\section{Fracture indicators}
\label{s:indicators}

We address in this section the definition of {\em fault\/} and
{\em barrier indicators\/} for the selection of additional {\candidate}
faults and barriers which are likely to produce a better fit to the data.

Let the current fracture, described by its
locations~$\currFh=(\currFalf,\currFbet)$ and parameters~$\alf$
and~$\bet$, be given by {\apriori} knowledge or from previous computations.
$\currFh$ may be void, for instance at the beginning of the fracture
determination procedure when no {\apriori} information is available.
Let $X=(\Pt,\Pe,\Pn,\Ute,\Ven)$ be the solution of the current
model~\eqref{eq:fracmodel:direct:discr:Omega}
and~\eqref{eq:fracmodel:direct:discr:gamma:chgvar}, which reads
\begin{equation}
  \label{eq:fracmodel:direct:discr:current}
  \begin{array}{rcll}
    \dps \sum_{E\in\ET} \UTE & = & \ffT, &
    \mbox{for all } T\in\Th, \\
    \dps \frac{\kh}{2} \, \UTE & = & \PT - \PE, &
    \mbox{for all } E\notin\currFbet,\, T\in\TE, \\
    \dps \left( \frac{\kh}{2} + \betE \right) \, \UTE & = & \PT - \PE, &
    \mbox{for all } E\in\currFbet,\, T\in\TE, \\
    \dps \PE & = & \PmE, &
    \mbox{for all } E\in\Eeh, \\
    \dps - \sum_{T\in\TE} \UTE & = & 0, &
    \mbox{for all } E\notin\currFalf,\, E\in\Eih, \\
    \dps - \sum_{T\in\TE} \UTE + \alfE \sum_{N\in\NE} \VEN & = & 0, &
    \mbox{for all } E\in\currFalf, \\
    \dps \frac{h}{2} \, \VEN & = & \PE - \PN, &
    \mbox{for all } E\in\currFalf,\, N\in\NE, \\
    \dps \sum_{\substack{E\in\currFalf \\ E\in\EN}} \alfE \, \VEN & = & 0, &
    \mbox{for all } N\in\currNalf,
  \end{array}
\end{equation}
and let $\Lambda\egaldef(\LAt,\LAe,\LAn,\MUte,\NUen)$ be the solution of the
{\em current adjoint model\/} obtained
from~\eqref{eq:fracmodel:adj:discr:Omega}
and~\eqref{eq:fracmodel:adj:discr:gamma:chgtvar} for a vanishing {\candidate}
fracture with zero intensity parameters,
\begin{equation}
  \label{eq:fracmodel:adj:discr:current}
  \begin{array}{rcll}
    \dps \sum_{E\in\ET} \MUTE & = & \PT - \PmT, &
    \mbox{for all } T\in\Th, \\
    \dps \frac{\kh}{2} \, \MUTE & = & \LAT - \LAE, &
    \mbox{for all } E\notin\currFbet,\, T\in\TE, \\
    \dps \left( \frac{\kh}{2} + \betE \right) \, \MUTE & = & \LAT - \LAE, &
    \mbox{for all } E\in\currFbet,\, T\in\TE, \\
    \dps \LAE & = & 0, &
    \mbox{for all } E\in\Eeh, \\
    \dps - \sum_{T\in\TE} \MUTE & = & 0, &
    \mbox{for all } E\notin\currFalf,\, E\in\Eih, \\
    \dps - \sum_{T\in\TE} \MUTE + \alfE \sum_{N\in\NE} \NUEN & = & 0, &
    \mbox{for all } E\in\currFalf, \\
    \dps \frac{h}{2} \, \NUEN & = & \LAE - \LAN, &
    \mbox{for all } E\in\currFalf,\, E\in\EN, \\
    \dps \sum_{\substack{E\in\currFalf \\ E\in\EN}} \alfE \NUEN & = & 0, &
    \mbox{for all } N\in\currNalf.
  \end{array}
\end{equation}

If the current fracture~$\currFh$ does not give a small enough value to the
objective function~$\J{\currFh}(\alf,\bet)$ defined
in~\eqref{eq:def:cost:fun}, one considers testing the effect of adding a set
$\tentFalf\subset\Eih$ of new fault edges and/or a set $\tentFbet\subset\Eih$
of new barrier edges.
As in Section~\ref{s:inverse}, we suppose that these {\candidate} fractures
have no common edge with the current fracture, {\ie}
satisfy~\eqref{eq:FcapDF=0} and~\eqref{eq:NcapDN=0}.

In order to select {\candidate} fractures~$\tentFh=(\tentFalf,\tentFbet)$
which are likely to decrease the objective function in the sense
of~\eqref{tentative:decrease:objective:function}, we make a first-order
development of~$\JtFh$ defined in~\eqref{eq:def:cost:fun:eps},
\begin{eqnarray}
  \label{eq:indic:first:order}
  \JtFh (\epsalf, \epsbet)
  - \underbrace{\JtFh (0, 0)}_{\JcFh (\alf, \bet)}
  & \simeq &
  \DJtFhepsalf (0, 0) \, \epsalf + \DJtFhepsbet (0, 0) \, \epsbet.
\end{eqnarray}
Then, following~\cite[Section~3.7]{chavent:book:2009}, we define the
{\em fault indicator\/}~$\Ialf$ and the {\em barrier indicator\/}~$\Ibet$
associated respectively with the {\candidate} fault~$\tentFalf$ and the
{\candidate} barrier~$\tentFbet$ by
\begin{equation}
  \label{eq:indic}
  \Ialf (\tentFalf) \egaldef \DJtFhepsalf (0, 0)
  \quad\mbox{and}\quad
  \Ibet (\tentFbet) \egaldef \DJtFhepsbet (0, 0).
\end{equation}
When these indicators are negative and of large absolute value, $\tentFalf$
and/or~$\tentFbet$ will be good {\candidate} fault or barrier, as they should
provide a strong decrease of the cost function (at least at first order) when
the intensity parameters~$\epsalf$ or~$\epsbet$ increase.
Of course, these fracture indicators only give first-order information, but
they are inexpensive to compute, as we see now.

Equation~\eqref{eq:grad:oldadjvar} gives for the fault indicator the formula
\begin{equation}
  \label{eq:indic:alf}
  \Ialf (\tentFalf) =
  - \frac{h}{2} \,
  \sum_{E\in\tentFalf} \alfbarE \sum_{N\in\NE} \VEN \, \NUEN,
\end{equation}
where $(\VEN,\NUEN)_{E\in\tentFalf,\,N\in\NE}$ are given by the direct limit
fault model~\eqref{eq:fracmodel:direct:discr:gamma:chgvar:alfa=0} and its
adjoint counterpart set on the {\candidate} fault edges of~$\tentFalf$ only,
whose solutions are given by
\begin{equation}
  \label{eq:limit:fault:model:solution}
  \begin{array}{lll}
    \PN =
    \dps \frac{\sum_{\substack{E\in\tentFalf \\ E\in\EN}} \alfbarE \, \PE}%
    {\sum_{\substack{E\in\tentFalf \\ E\in\EN}} \alfbarE}, &
    \LAN =
    \dps \frac{\sum_{\substack{E\in\tentFalf \\ E\in\EN}} \alfbarE \, \LAE}%
    {\sum_{\substack{E\in\tentFalf \\ E\in\EN}} \alfbarE}, &
    \mbox{for all } N\in\tentNalf, \\
    \dps \frac{h}{2} \, \VEN = \PE - \PN, &
    \dps \frac{h}{2} \, \NUEN = \LAE - \LAN, &
    \mbox{for all } E\in\tentFalf,\, N\in\NE.
  \end{array}
\end{equation}
Here, $(\PE,\LAE)_{E\in\tentFalf}$ and
$(\PN)_{N\in\currNalf\cap\NE,\,E\in\tentFalf}$ are simply extracted from the
solutions~$X$ and~$\Lambda$ of the current direct and adjoint
models~\eqref{eq:fracmodel:direct:discr:current}
and~\eqref{eq:fracmodel:adj:discr:current}.

Similarly, the second equation in~\eqref{eq:grad} gives for the barrier
indicator the formula
\begin{equation}
  \label{eq:indic:bet}
  \Ibet (\tentFbet) =
  \sum_{E\in\tentFbet} \betbarE \sum_{T\in\TE} \UTE \, \MUTE,
\end{equation}
where $(\UTE,\MUTE)_{E\in\tentFbet,\,T\in\TE}$ are simply extracted from the
solutions~$X$ and~$\Lambda$ of the current direct and adjoint
models~\eqref{eq:fracmodel:direct:discr:current}
and~\eqref{eq:fracmodel:adj:discr:current}.

\begin{Remark}
  \label{rem:physical:indicator}
  The fault indicator formula~\eqref{eq:indic:alf} depends on the direct and
  adjoint tangential Darcy velocities along the {\candidate} fault.
  From Section~\ref{ss:direct:continuous}, this is physically sound, since
  the fault parameter~$\alf$ is linked to the flow\/ {\em along} the
  {\candidate} fault.
  In the same way, the barrier indicator formula~\eqref{eq:indic:bet} depends
  on the direct and adjoint Darcy velocities normal to the {\candidate}
  barrier, which is natural, since~$\bet$ is linked to the flow\/
  {\em across} to the {\candidate} barrier.
\end{Remark}

\begin{Remark}
  \label{rem:null:fault:indicator:isolated:edge}
  If~$\tentFalf$ is made of a single edge~$E$ with no endpoint on the current
  fault~$\currFalf$, one checks easily from the last equation
  of~\eqref{eq:fracmodel:direct:discr:tentative} that
  $\JtFh(\epsalf,0)=\JcFh(\alf,\bet)$ for all $\epsalf\geq 0$, {\ie}
  opening a fault on a single edge does not change the pressure and flow
  patterns.
  Notice also that in this case,
  \eqref{eq:limit:fault:model:solution}~implies immediately that $\VEN=0$ at
  the endpoints of~$E$, and the corresponding fault indicator is zero (which
  is coherent).
  Thus, it will be necessary to use only {\candidate} faults constituted of\/
  {\em at least two connected edges}, or having a\/
  {\em common node with the current fault}.
\end{Remark}

\bigskip

To conclude this section, we remark that once the current direct and adjoint
models~\eqref{eq:fracmodel:direct:discr:current}
and~\eqref{eq:fracmodel:adj:discr:current} have been solved, the computation
of the indicators~$\Ialf(\tentFalf)$ or~$\Ibet(\tentFbet)$ only requires the
summation of known quantities along the {\candidate} fractures~$\tentFalf$
or~$\tentFbet$ (plus the resolution of~\eqref{eq:limit:fault:model:solution}
in the case of faults).
Hence, the computation of the indicators for many {\candidate} faults or
barriers is inexpensive since it requires the resolution of
{\em only one set of direct and adjoint model\/} with the current
fracture~$\currFh$.
This is the basis of the indicator algorithm developed in
Section~\ref{s:algorithm}.

%% file: frac_algorithm.tex

\section{Algorithm for the estimation of fractures}
\label{s:algorithm}

We have now all the ingredients for the presentation of an iterative
algorithm, based on fracture indicators, designed for the estimation of
discrete fractures in a porous medium from given pressure and/or flow
measurements.

We suppose in this section that pressure data are available in all elements
of the computational mesh, but the case of coarser measurements, or of flow
data can be treated as well by modifying the objective
function~\eqref{eq:def:cost:fun:eps} and the adjoint state
model~\eqref{eq:fracmodel:adj:discr:current} according to
Remark~\ref{rem:measures}.

The simultaneous determination of faults and barriers is a delicate matter:
the number of parameters is multiplied by two, and this can lead to
underdetermination unless a large number of measurements are available.
Hence, we limit ourselves to an algorithm for estimating either faults or
barriers, but not both at the same time.
We present the algorithm under a generic form, hence throughout this section,
we drop the subscripts~${}_\alf$ and~${}_\bet$: $\currF$,
$\tentF$, and~$\eps$ will refer to either~$\currFalf$, $\tentFalf$,
and~$\epsalf$, or~$\currFbet$, $\tentFbet$, and~$\epsbet$, depending on the
nature of the sought fracture.
When necessary, the fracture parameter of interest~$\alf$ or~$\bet$ will be
denoted~$\xi$.

Fractures have been assumed so far to be any set of internal edges.
Again, this can lead to underdetermination, which is usually associated with
a poor conditioning of the inverse problem and with the presence of many
meaningless local minima for the objective function.
Moreover, we are only interested in locating fractures that have a
macroscopic impact on the flow.
Therefore, we make the assumption that fractures are defined on a coarser
mesh than the computational mesh.
Nevertheless, the amount of possible fractures for such a coarse
representation can still be quite large, so in order to avoid
overparameterization, we use the fracture indicators defined in
Section~\ref{s:indicators} to introduce fractures one at a time,
and we use~\eqref{param:alfa:vanish:prop} and~\eqref{param:beta:vanish:prop}
to parameterize each fracture with only one parameter~$\eps$.

\subsection{The fracture mesh}
\label{ss:fracture:mesh}

To exclude very small fractures, typically made of only a few edges, we
decouple the {\em computational mesh\/}~$\Th$ from the
{\em fracture research mesh\/}~$\Thc$.
The computational mesh should be fine enough to perform accurate
calculations, and the fracture mesh should correspond to the expected size of
the fracture.
To avoid interpolation issues, we simply choose the fine mesh~$\Th$ to be a
refinement of the coarse mesh~$\Thc$.
Hence, coarse edges are made up of several computational edges.

From now on, we suppose that both current and {\candidate} fractures are made
of internal coarse edges, and we denote by~$\currFc$ (resp.~$\tentFc$)
the set of coarse edges belonging to the current fracture (resp. a
{\candidate} fracture).
In term of computational edges, these fractures are given by
\begin{equation}
  \label{frac:fine:from:coarse}
  \begin{array}{c}
    \currF =
    \{ E\in\Eih \suchthat \exists \Ec\in\currFc, E\subset\Ec \}, \\
    \tentF =
    \{ E\in\Eih \suchthat \exists \Ec\in\tentFc, E\subset\Ec \}.
  \end{array}
\end{equation}
It is understood that in the sequel all current and {\candidate}
fractures~$\currF$ and~$\tentF$ are of the
form~\eqref{frac:fine:from:coarse} for some~$\currFc$ and~$\tentFc$,
{\ie} actually made of coarse edges.

\subsection{A collection of fractures}
\label{ss:collection:fractures}

At each iteration~$\ell\geq 1$, a new fracture~$\besttentFl$ is added to the
current collection of fractures~$\currFml$.
Hence, after $k-1$ iterations, the current collection of fractures is of the
form
\begin{equation}
  \label{structure:currFk-1}
  \currFmk =
  \currFof{0} \cup \besttentFof{1} \cup \ldots \cup \besttentFmk,
\end{equation}
where~$\currFof{0}$ is the (usually void) set of given {\apriori}
fractures.

The setting of Sections~\ref{s:inverse} and~\ref{s:indicators} is replicated
for each fracture.
Let~$\epsmkl$ be the intensity parameter associated with the
fracture~$\besttentFl$.
Equations~\eqref{param:alfa:vanish:prop} and~\eqref{param:beta:vanish:prop}
become (different fractures have no common edges)
\begin{equation}
  \label{current:frac:param:from:intensity:param}
  \xiE = \epsmkl \xibarE,
  \quad \forall E\in\besttentFl,\, \forall \ell=1,\ldots,k-1.
\end{equation}
When {\apriori} information is available, values~$\xi^0$ of the fracture
parameter on the initial set of fractures~$\currFof{0}$ are also provided;
they are kept fixed.
Then, the data misfit of equation~\eqref{eq:def:cost:fun:eps} becomes a
function of the vector of intensity parameters,
\begin{equation}
  \label{eq:def:cost:fun:intensity:param}
  \J{\currFmk} (\epsmkofcomp{1}, \dots, \epsmkofcomp{k-1})
  \egaldef \half \sum_{T\in\Th} (\PT - \PmT)^2,
\end{equation}
where~$\PT$ is given by the current
model~\eqref{eq:fracmodel:direct:discr:current}
with the current collection of fractures~$\currFmk$.
The {\em optimal\/} intensities~$\bestepsmk$ are the minimizer of the
objective function satisfying
\begin{equation}
  \label{best:fracture:intensity:parameters:k-1}
  \bestepsmk
  = (\bestepsmkofcomp{1}, \ldots, \bestepsmkofcomp{k-1})
  \egaldef \arg \min_{\epsmkl>0}
  \J{\currFmk} (\epsmkofcomp{1}, \dots, \epsmkofcomp{k-1}).
\end{equation}

As the current collection of fractures is known by its locations~$\currFmk$
and intensities~$\bestepsmk$, one can solve the associated current direct and
adjoint models~\eqref{eq:fracmodel:direct:discr:current}
and~\eqref{eq:fracmodel:adj:discr:current}, complemented
with~\eqref{current:frac:param:from:intensity:param} for the optimal
intensities~$\bestepsmk$.
This determines the current direct and adjoint variables
$X^{k-1}=(\Pt,\Pe,\Pn,\Ute,\Ven)$ and
$\Lambda^{k-1}=(\LAt,\LAe,\LAn,\MUte,\NUen)$.

To sum up, after $k-1$ iterations, the following quantities are available:
the collection of $k-1$ estimated fracture locations~$\besttentFof{1}$,
\ldots, $\besttentFmk$, the associated optimal
intensities~$\bestepsmk=(\bestepsmkofcomp{1},\ldots,\bestepsmkofcomp{k-1})$,
and the corresponding vectors of current direct and adjoint
variables~$X^{k-1}$ and~$\Lambda^{k-1}$.

The case $k=1$ corresponds to the initialization of the algorithm.
It is specific: there is no estimated fracture ({\ie} no intensity), $X^0$
and~$\Lambda^0$ are obtained as solution
of~\eqref{eq:fracmodel:direct:discr:current}
and~\eqref{eq:fracmodel:adj:discr:current} with the set of given {\apriori}
fractures~$\currFof{0}$.

\subsection{The next iteration}
\label{ss:next:iteration}

When the fit to the data obtained
in~\eqref{best:fracture:intensity:parameters:k-1} after $k-1$ iterations is
not satisfactory, one needs to perform an additional iteration.
This means selecting a new fracture~$\besttentFk$ enabling the data misfit to
decrease, and determining the resulting optimal fracture intensities
$\bestepsk=
(\bestepskofcomp{1},\ldots,\bestepskofcomp{k-1},\bestepskofcomp{k})$
for the new collection of fractures $\currFk=\currFmk\cup\besttentFk$.
This is achieved in three steps.

\myparagraph{Indicator step}
A {\em long list of {\candidate} fractures} is chosen.
As in Sections~\ref{s:inverse} and~\ref{s:indicators}, the {\candidate}
fractures must have no common edge with the current collection of
fractures~$\currFmk$.
The fracture indicators~$\I(\tentF)$ for all {\candidate} fractures~$\tentF$
of the long list are computed by~\eqref{eq:indic:alf}
and~\eqref{eq:limit:fault:model:solution} for faults, or~\eqref{eq:indic:bet}
for barriers.
These computations are very fast as all terms are already available from the
current direct and adjoint variables~$X^{k-1}$ and~$\Lambda^{k-1}$ provided
by the previous iteration.

Then, a {\em short list of {\candidate} fractures} is built.
The {\candidate} fractures of the short list are associated with (negative)
fracture indicators of large magnitude.
Thus, according to the first-order information carried by the indicators,
they are the most likely to produce a significant enhancement of the fit to
the data.

Strategies to choose the long list and to build the short list are discussed
in Section~\ref{ss:build:cand:frac}.

\myparagraph{Optimization step}
For {\em each {\candidate} fracture~$\tentF$ of the short list}, solve the
minimization problem
\begin{equation}
  \label{best:fracture:intensity:parameters:k}
  \epstk (\tentF) = (\epstkofcomp{1}, \ldots, \epstkofcomp{k})
  = \arg \min_{\epskl>0}
  \J{\currFmk\cup\tentF} (\epskofcomp{1}, \ldots, \epskofcomp{k}),
\end{equation}
where the objective function~$\J{\currFmk\cup\tentF}$ is similar
to~\eqref{eq:def:cost:fun:intensity:param}, here $\PT$ is still given by the
current model~\eqref{eq:fracmodel:direct:discr:current}, but set on
$\currFmk\cup\tentF$ with the~$k$ fracture intensities
$(\epskofcomp{1},\ldots,\epskofcomp{k})$.
Thanks to the use of a nominal parameter~$\xibar$, the natural initial guess
for the minimization is
$(\bestepsmkofcomp{1},\ldots,\bestepsmkofcomp{k-1},1)$.

The computation of the gradient of~$\J{\currFmk\cup\tentF}$ is required at
each iteration of the minimization algorithm.
Hence, the need for the resolution of both direct and adjoint
models~\eqref{eq:fracmodel:direct:discr:current}
and~\eqref{eq:fracmodel:adj:discr:current} set on $\currFmk\cup\tentF$ with
the current value of $(\epskofcomp{1},\ldots,\epskofcomp{k})$ at each
minimization iteration.
Thus, each minimization is much more computationally intensive than the
calculation of the fracture indicators by~\eqref{eq:indic:alf}
and~\eqref{eq:limit:fault:model:solution}, or~\eqref{eq:indic:bet}.
This is why a short list of a small number of {\candidate} fractures is
built.
The optimization step is usually more expensive than the indicator step, by
several orders of magnitude.

\myparagraph{Update step}
The winner~$\besttentFk$ is the {\candidate} fracture that gives the smallest
minimum value to the objective function.
The new current collection of fractures is $\currFk=\currFmk\cup\besttentFk$,
and the new optimal fracture intensity vector is
$\bestepsof{k}=\epstk(\besttentFk)$.

Finally, the new vectors of direct and adjoint variables~$X^k$
and~$\Lambda^k$ are determined by solving direct and adjoint
models~\eqref{eq:fracmodel:direct:discr:current}
and~\eqref{eq:fracmodel:adj:discr:current} set on~$\currFk$ with
intensities~$\bestepsk$.

\subsection{Stopping the algorithm}
\label{ss:algo:stopping}

The algorithm stops either when the maximal number of fractures is reached,
when adding a new fracture does not significantly improve the objective
function, or when the data misfit is of the same magnitude as the uncertainty
on the given data measures.
The constants driving the stopping criteria are:
the absolute uncertainty on the data measures~$\Dmeas\geq 0$,
the relative tolerance for convergence~$\etaconv>0$,
the relative tolerance for stationary sequence~$\etastat>0$,
and the maximum number of admissible fractures~$\kmax$.

When the algorithm stops, the result is given in terms of the fracture
parameter: $\xi^k$ is restored from the {\apriori} given values~$\xi^0$ on
the initial set of fractures~$\currFof{0}$, and from the last estimated
intensities~$\bestepsk$
through~\eqref{current:frac:param:from:intensity:param} (expressed for~$k$
instead of~$k-1$).

\subsection{Algorithm}
\label{ss:algo}

Here comes the precise description of the algorithm.
Remember that the fault and barrier indices~${}_\alf$ and~${}_\bet$ are
omitted: $\currF$, $\tentF$, and~$\eps$ stand either for~$\currFalf$,
$\tentFalf$, and~$\epsalf$, or for~$\currFbet$, $\tentFbet$, and~$\epsbet$.
Moreover, $\xi$ stands either for~$\alf$ or~$\bet$.

{\Apriori} information, or maybe a previous computation, provides the initial
guess for fractures (location~$\currFof{0}$ and fracture parameter~$\xi^0$,
which may be void), and a value~$\xibarE$ of the nominal parameter in each
edge~$E$ of the mesh.

\begin{description}
\item[Initialization.] \mbox{}
  \begin{enumerate}
  \item
    Compute the initial solutions~$X^0$ and~$\Lambda^0$ of the direct and
    adjoint models~\eqref{eq:fracmodel:direct:discr:current}
    and~\eqref{eq:fracmodel:adj:discr:current} set on the initial
    fractures~$\currFof{0}$ with fracture parameter~$\xi^0$.
  \item
    Compute the initial objective function
    $\bestcostof{0}=\J{\currFof{0}}(\xi^0)$ with~\eqref{eq:def:cost:fun}.
  \end{enumerate}
\item[Iterations.]
  For $k\geq 1$, do:
  \begin{enumerate}
  \item
    \textbf{Indicator step.}
    Build the short list of {\candidate} fractures according to the chosen
    strategy (see Section~\ref{ss:build:cand:frac}).
    This uses fracture indicators~\eqref{eq:indic:alf}
    and~\eqref{eq:limit:fault:model:solution}, or~\eqref{eq:indic:bet}, that
    depends on~$X^{k-1}$ and~$\Lambda^{k-1}$.
  \item
    If the short list is empty,
    then the algorithm stops,
    the result is~$\currFmk$ (with $k-1$ estimated fractures),
    and~$\xi^{k-1}$ restored from~$\xi^0$ and intensities~$\bestepsmk$
    through~\eqref{current:frac:param:from:intensity:param}.
  \item
    \textbf{Optimization step.}
    For each {\candidate} fracture~$\tentF$ in the short list,
    solve~\eqref{best:fracture:intensity:parameters:k}:
    \begin{itemize}
    \item
      minimize the objective function $\J{\currFmk\cup\tentF}$ with
      respect to the~$k$ fracture intensity parameters
      $\epsk=(\epskofcomp{1},\ldots,\epskofcomp{k})$;
    \item
      call $\epstk(\tentF)$ the minimizer.
    \end{itemize}
  \item
    \textbf{Update step.}
    Retain in the short list the {\candidate} fracture~$\besttentFk$ that
    gives the smallest value to
    $\J{\currFmk\cup\tentF}(\epstk(\tentF))$, and set
    \begin{equation*}
      \currFk = \currFmk \cup \besttentFk,
      \quad
      \bestepsk = \epstk (\besttentFk),
      \quad
      \bestcostk = \J{\currFk} (\bestepsk).
    \end{equation*}
    Then, compute the solutions~$X^k$ and~$\Lambda^k$ of the current direct
    and adjoint models~\eqref{eq:fracmodel:direct:discr:current}
    and~\eqref{eq:fracmodel:adj:discr:current} set on the current collection
    of fractures~$\currFk$ with intensities~$\bestepsk$.
  \item
    If $\bestcostmk-\bestcostk\leq\etastat\,\bestcostof{0}$,
    then the algorithm is stationary,
    the result is~$\currFmk$ (with $k-1$ estimated fractures),
    and~$\xi^{k-1}$ restored from~$\xi^0$ and intensities~$\bestepsmk$.
  \item
    If $\bestcostk-\Dmeas\leq\etaconv\,\bestcostof{0}$,
    then the algorithm has converged,
    the result is~$\currFk$ (with~$k$ estimated fractures), and~$\xi^k$
    restored from~$\xi^0$ and intensities~$\bestepsk$.
  \item
    If $k\geq\kmax$,
    then the algorithm has reached the maximum number of expected fractures,
    the result is~$\currFk$ (with~$k$ estimated fractures), and~$\xi^k$
    restored from~$\xi^0$ and intensities~$\bestepsk$.
  \end{enumerate}
\end{description}

\subsection{Strategy to build a short list of {\candidate} fractures}
\label{ss:build:cand:frac}

The general idea is to choose a {\em long list\/} of {\candidate} fractures,
trying to be as exhaustive as possible, or to incorporate {\apriori}
knowledge, and then to use the inexpensive first-order indicators to select a
{\em short list\/} of {\candidate} fractures that are likely to provide a
large decrease of the objective function.

For the design of {\em refinement indicators\/}
in~\cite{BAChJa:ref:coarse:ind:2002,BAClWCh:jiipp:08,BAChClW:ipse:11},
the unknowns were of a completely different nature: the parameters were
distributed all over the domain, with a value in each element of the mesh.
The algorithm sought for the location of parameter discontinuities
represented by {\em cutting\/} curves.
The long list of {\candidate} cuttings was chosen from a small number of
predefined families of curves, such as parallel lines, or circles.
Here, some fractures have to be created at specific locations on the edges of
the mesh.
It would have been possible to use the same kind of choices for the long
list of {\candidate} fractures.
Instead, we propose a constructive strategy.

The key idea is to choose the initial long list as a collection of
{\em elementary {\candidate} fractures\/} that allow for a cartography of the
values of the fracture indicators all over the domain.
Elementary fractures are made of a few edges to allow for exhaustiveness, and
they are meant to grow into larger fractures having a macroscopic impact on
the flow and associated with (negative) fracture indicators of larger
magnitude.

The strategy depends on the choice of a coarse fracture research mesh~$\Thc$
of which the computational mesh~$\Th$ is a refinement.
It is also parameterized by constants driving the selection criteria:
two ratios~$\thetaelem$ and~$\thetaext$ (between~0 and~1) for indicator-based
selection of {\candidate} fractures, and the maximum admissible number of
{\candidate} fractures~$\nmax$.

The three stages of the strategy go as follows.
First, according to Remark~\ref{rem:null:fault:indicator:isolated:edge}, the
fracture indicators are computed for a long list of elementary fractures made
of all pairs of contiguous internal coarse edges with at most one node in
common with the current collection of fractures, complemented with all single
internal coarse edges with one node on the current collection of fractures.
Let~$\bestIelem$ be the best ({\ie} minimum) indicator for all elementary
fractures (remember that {\em good\/} indicators are {\em very negative}).
If this minimum is positive, then the short list of {\candidate} fractures is
empty (and the fracture indicator algorithm stops).
Otherwise, we select all elementary fractures whose indicator is lower or
equal to $\thetaelem\,\bestIelem$.

Then, we aggregate together all selected elementary fractures that share a
common coarse edge to form the longest possible {\candidate} fractures.
The minimum indicator for all aggregates of elementary fractures~$\bestIaggr$
is usually lower than~$\bestIelem$, and there is no need for a selection here
since the number of aggregates is smaller than the number of selected pairs.

Finally, we allow for an extension of all aggregates by one coarse edge at
one or two of its endpoints.
Note that there could be three endpoints, or more.
Note also that it can increase considerably the number of {\candidate}
fractures.
Let~$\bestIext$ be the best ({\ie} minimum) indicator for all extended
aggregates.
Again, this minimum is lower or equal to~$\bestIaggr$.
We select all extended aggregates whose indicator is lower or equal to
$\thetaext\,\bestIext$.
At the end, in order to control the cost of the optimization step of the
algorithm, we truncate the short list of {\candidate} fractures to
the~$\nmax$ selected extended aggregates associated with the lowest
indicators.

%% file: frac_numerics.tex

\section{Numerical results}
\label{s:numerics}

We present now numerical experiments conducted with a Matlab$^\circledR$
implementation of the fracture indicator-based algorithm described in
Section~\ref{s:algorithm} for the estimation of discrete fractures in a
porous medium from given pressure measurements.

First of all, as it must always be the case, the implementation of the
closed-form formulas~\eqref{eq:grad}, \eqref{eq:grad:oldadjvar},
\eqref{eq:indic:alf}, and~\eqref{eq:indic:bet} for the gradient of the
objective function with respect to the intensity parameters, and for the
fracture indicators, have been successfully checked by comparison with a
finite-difference approximation based on the sole computation of the
objective function.

The test-cases are simple, but represent typical situations of interest:
faults parallel to the main direction of the flow, and barriers normal to the
flow.
The synthetic pressure and flow data are obtained by using the same
simulation program: the tests range from the so-called ``inverse crime'' for
which the very same model is used for the inversion and simulation of data,
to situations where random noise of increasing level is added, and also to
the much more difficult but interesting situation where the target fractures
are not carried by the fracture research mesh.

We do not use any {\apriori} information on the target fractures: the initial
estimation is no fracture ($\currFof{0}=\emptyset$), and nominal
parameters are taken constant $\alfbar=\betbar=1$;
thus $\epsalf=\alf$ and $\epsbet=\bet$.
The constants for tuning the behavior of the algorithm are chosen as follows:
$\Dmeas$ is proportional to the noise level, $\etaconv=\etastat=10^{-2}$, and
$\kmax=8$ for the stopping criteria, and $\thetaelem=0.8$, $\thetaext=0.9$,
and $\nmax=10$ for the building of the short list of {\candidate} fractures.

\myparagraph{Three levels of mesh}
To allow for the decimation of the measurements, in addition to the (fine)
computational mesh~$\Th$ and to the (coarse) fracture research mesh~$\Thc$,
we also use a specific mesh for the synthetic measurements, denoted
by~$\Thm$ (the sum in~\eqref{eq:def:cost:fun:intensity:param} is now meant
for $T\in\Thm$).

Let $\Omega$ be a rectangular domain, and~$\QuadNxy$ denote the regular
rectangular grid of $\Nx\times\Ny$ cells for positive integers~$\Nx$
and~$\Ny$.
In the sequel, the domain~$\Omega$ is the unit square, the computational mesh
is $\Th=\Quadcomp$, the measurement mesh is $\Thm=\Quadmeas$ with~$\Nm$
ranging from~72 down to~8, and the fracture mesh defined in
Section~\ref{ss:fracture:mesh} is $\Thc=\Quadfrac$ with $\Nf=12$ or~9.
Note that~$\Th$ is a refinement of~$\Thc$ since~$\Nf$ is a divisor of~72.

\myparagraph{Test-cases}
The porous medium is homogeneous with a permeability~$K=1$.
Impervious Neumann conditions are imposed on top and bottom, and Dirichlet
conditions impose a pressure drop from right ($\pmeas=1$) to left
($\pmeas=0$).
Without any fracture, the solution of the direct model is a linear pressure
(see Figure~\ref{fig:one:fault}b) and a uniform Darcy velocity.

A fault normal to the flow, or a barrier tangential to the flow, has a null
hydrogeological signature: the solution of the direct model is still the same
uniform flow from right to left.
Therefore, we consider either tangential faults, or normal barriers, for a
maximum effect on the solution.
In the following examples, the fractures have length~0.5 (half the size of
the domain), and they are located either in the middle of the domain (single
fracture), or centered at one quarter and at three quarters of the domain
(two fractures).
The fault parameter~$\alf$ (conductivity along the fracture) ranges
from~0.2 up to~200.
The barrier parameter~$\bet$ (resistivity across to the fracture) ranges
from~0.02 up to~200.

\myclearpage
\subsection{Illustration of the algorithm for the estimation of faults}
\label{ss:illustration:algo}

\begin{figure}[\figloc]
  \begin{center}
    (a)
    \includegraphics[width=\twofigs]{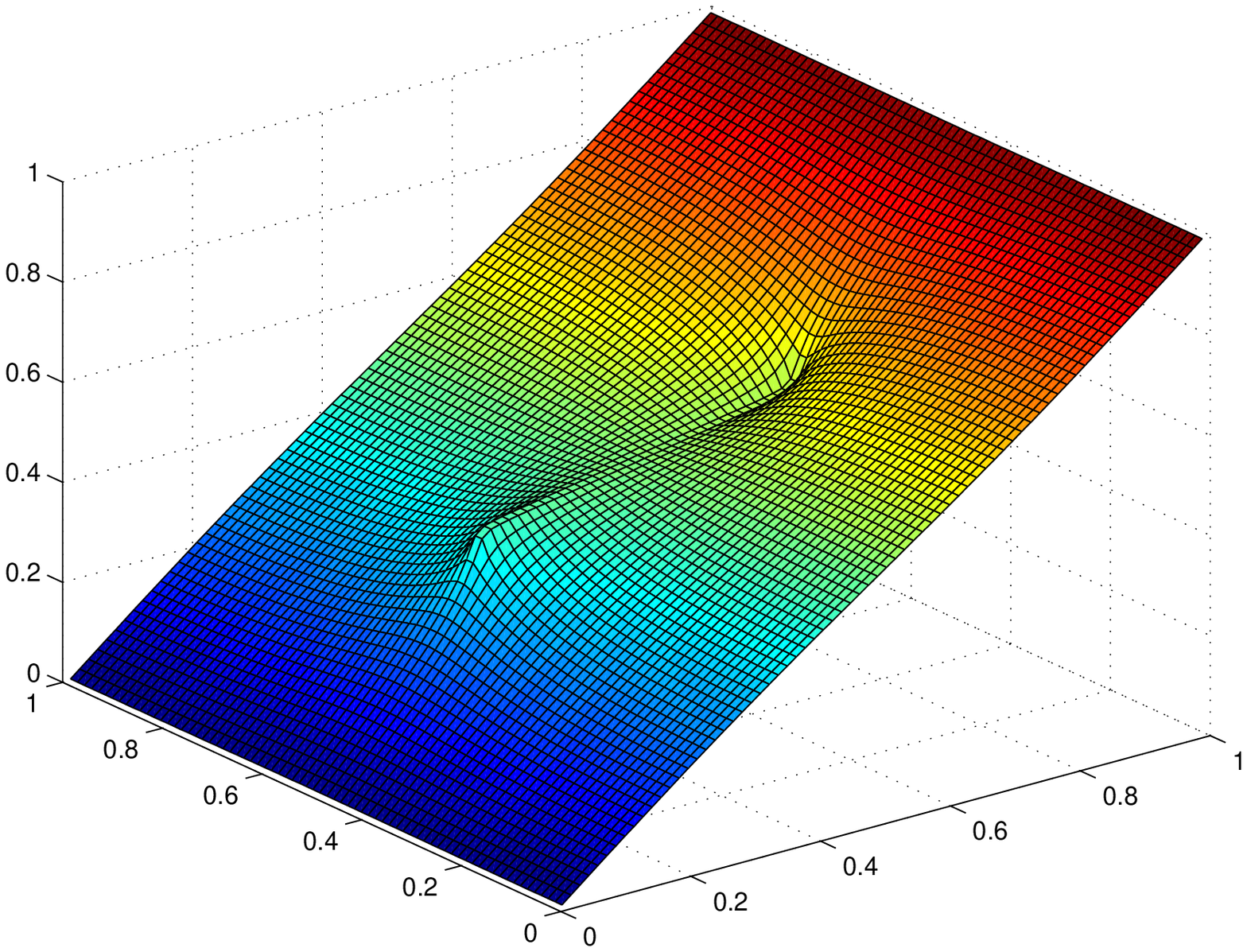}
    (b)
    \includegraphics[width=\twofigs]{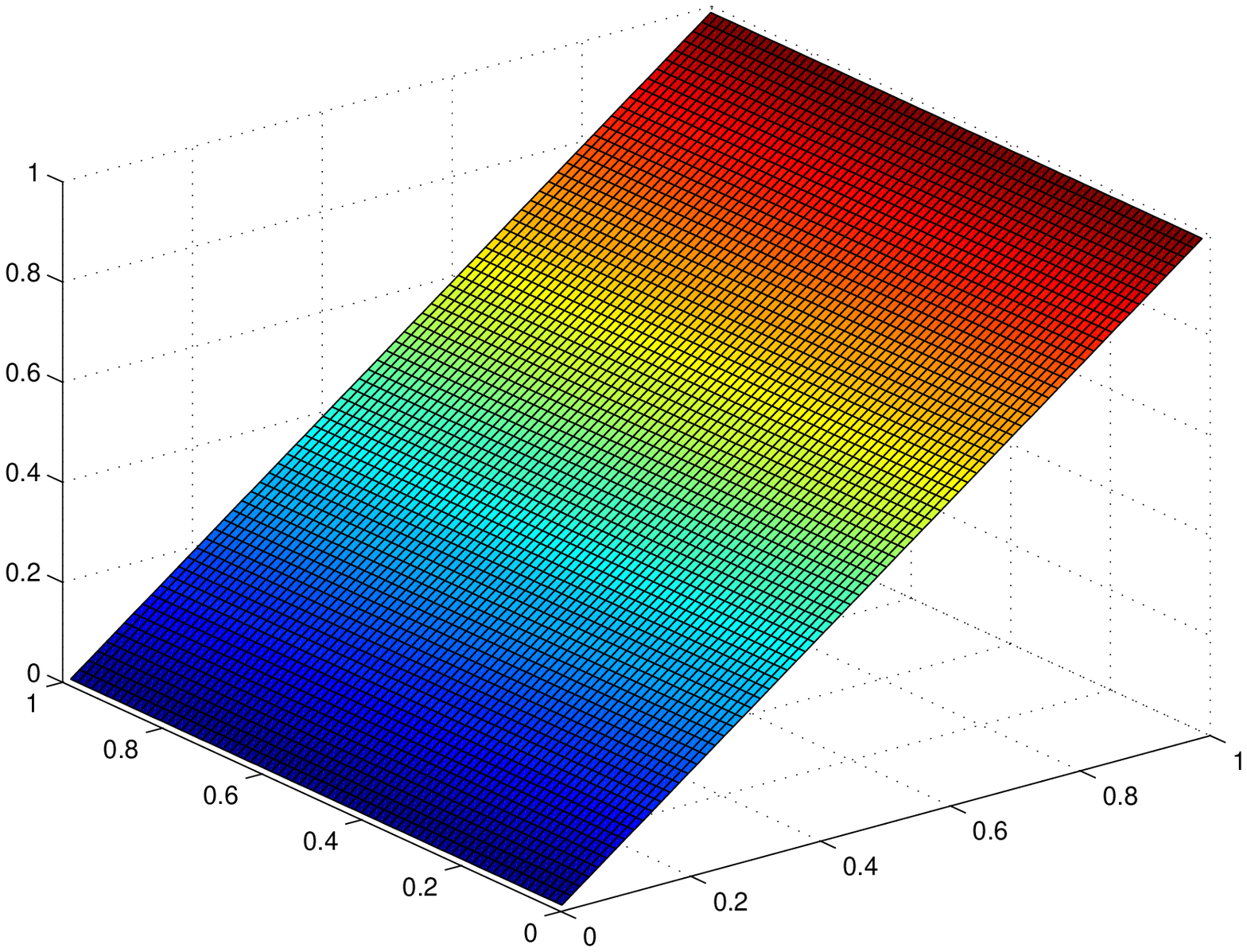} \\
    (c)
    \includegraphics[width=\twofigs]{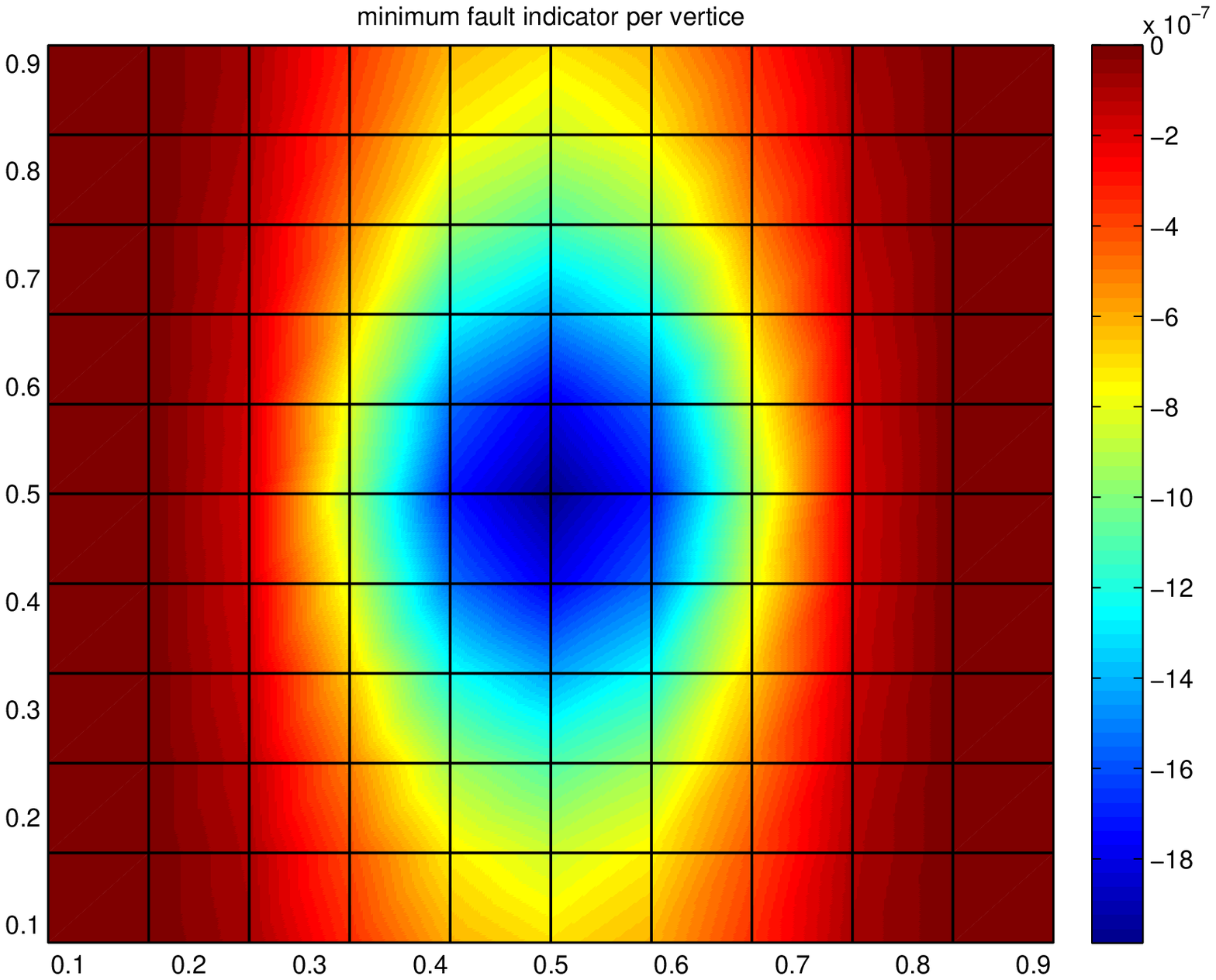}
    (d)
    \includegraphics[width=\twofigs]{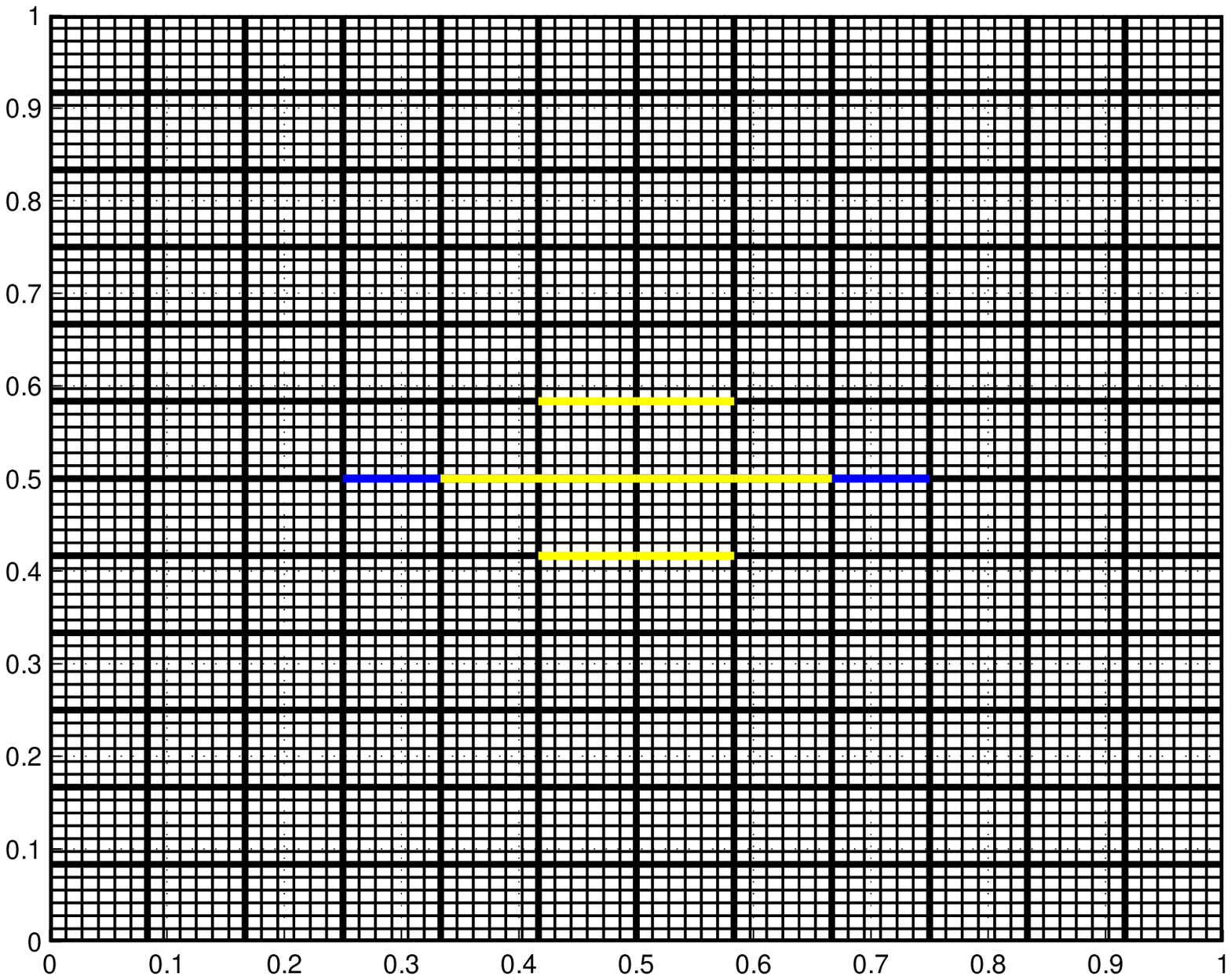} \\
    (e)
    \includegraphics[width=\twofigs]{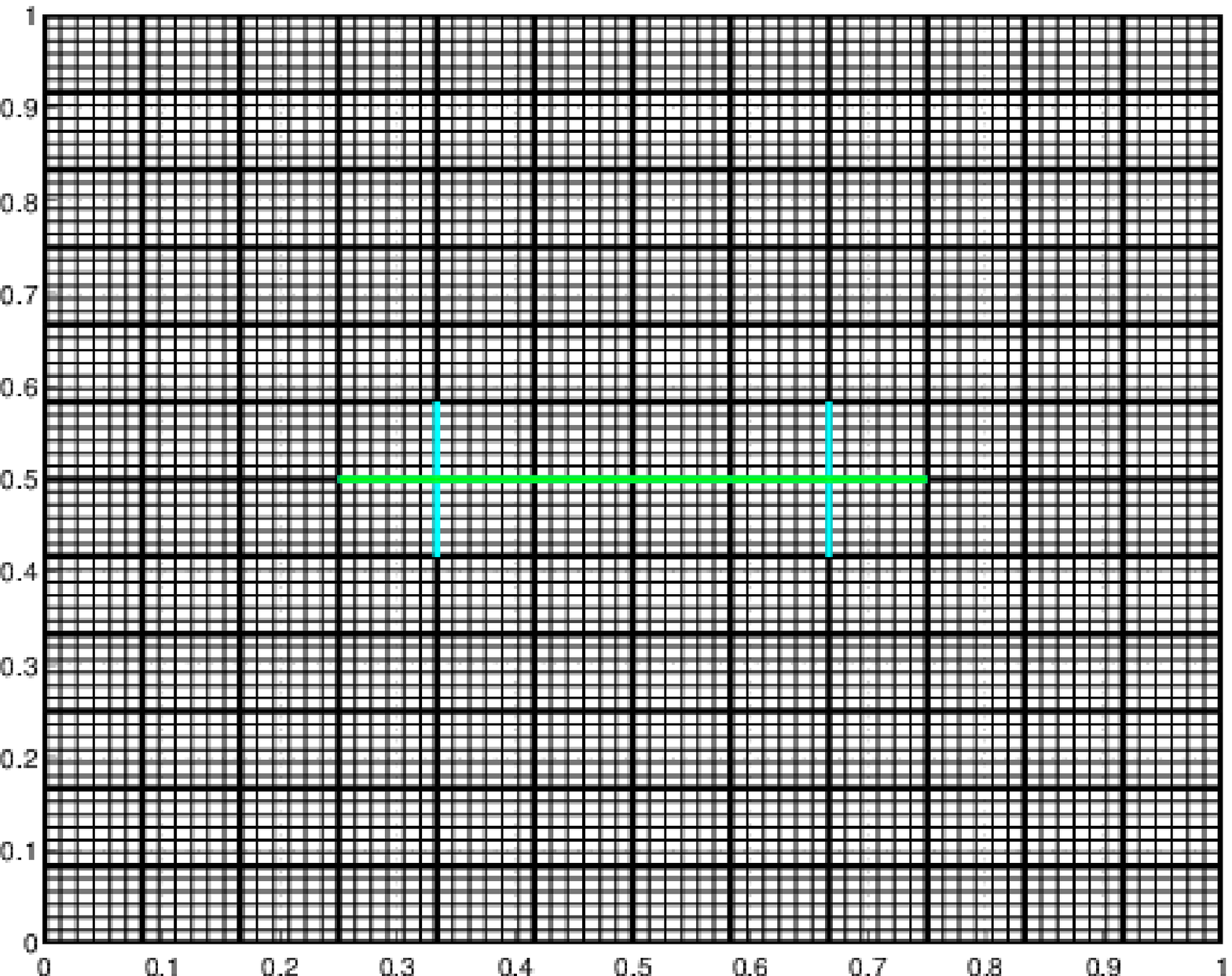}
    (f)
    \includegraphics[width=\twofigs]{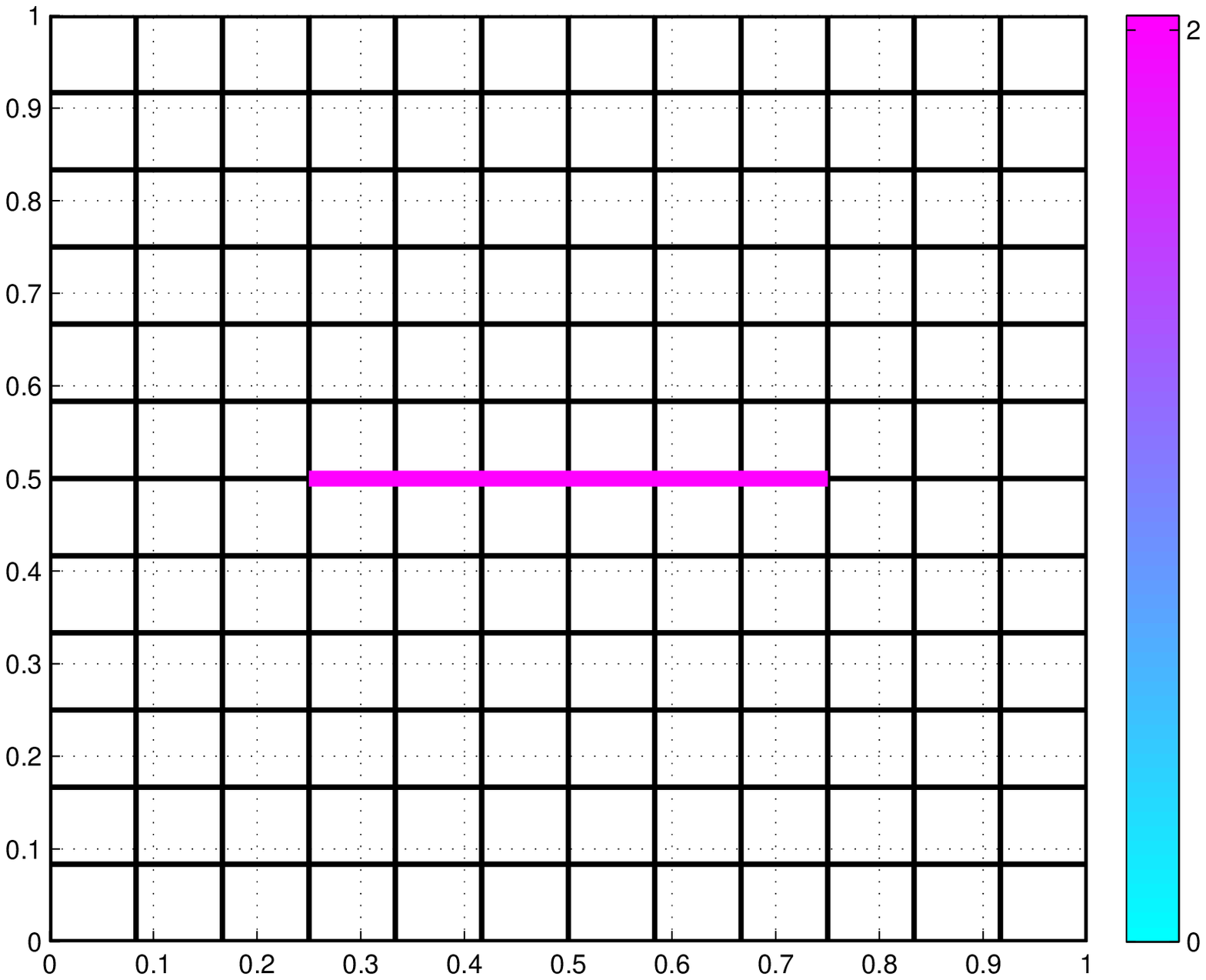}
    \caption{Case of a single fault tangential to the flow.
      \newline
      (a) Data pressure, solution of the direct
      model~\eqref{eq:fracmodel:direct:discr:current} with a fault at $y=0.5$
      ($\alf=2$).
      \newline
      (b) Initial pressure, solution of the direct model with no fracture.
      \newline
      (c) Distribution of indicators for the long list of elementary
      {\candidate} faults:
      the lowest indicator for all~6 pairs of edges centered in each interior
      coarse node is represented.
      \newline
      (d) In dark blue, the target fault;
      above in yellow, the selected pairs of edges (three of them are
      superimposed above the target fault).
      \newline
      (e) In light blue, the {\candidate} faults of the short list after
      aggregation and extension;
      above in green, the best {\candidate} fault.
      \newline
      (f) Best result after minimization for all {\candidate} faults of the
      short list (the lowest permeabilities are in light blue, and the
      highest in pink).
    }
    \label{fig:one:fault}
  \end{center}
\end{figure}

Let us first consider the simple case of a fault tangential to the flow
located in the middle of the domain with $\alf=2$.
The synthetic data are represented in Figure~\ref{fig:one:fault}a.
The pressure is measured in all elements of the computational mesh ({\ie}
$\Nm=72$).
The fracture research mesh corresponds to~$\Nf=12$, {\ie} each coarse edge is
made of 6 fine edges.

The algorithm is initialized with no fracture ($\currFof{0}=\emptyset$).
The initial direct and adjoint variables~$X^0$ and~$\Lambda^0$ are computed
as solution to~\eqref{eq:fracmodel:direct:discr:current}
and~\eqref{eq:fracmodel:adj:discr:current} with no fracture.
The initial pressure is linear (Figure~\ref{fig:one:fault}b).

The first iteration starts with the computation of the fault indicators
through~\eqref{eq:limit:fault:model:solution} and~\eqref{eq:indic:alf}
from~$X^0$ and~$\Lambda^0$ for all elementary {\candidate} faults of the long
list made of the~726 pairs of contiguous interior coarse edges
of~$\Thc=\Quad{12}{12}$.
The cartography of these indicators (Figure~\ref{fig:one:fault}c) shows a
nice blue spot in the middle of the domain, right where the target fault is
located (remember that good indicators are large negative numbers).
Only five pairs of coarse edges have an indicator lower than
$\thetaelem\bestIalfelem$;
they are parallel and close to the target fault
(Figure~\ref{fig:one:fault}d).

The aggregation stage collapses the 3 pairs of edges superimposed with the
target fault, reducing the number of {\candidate} faults down to~3
(Figure~\ref{fig:one:fault}d).
There is a significant improvement of the indicator as
$\bestIalfaggr\simeq 1.8\,\bestIalfelem$.

The extension stage adds one coarse edge at zero, one, or two of the
endpoints of the aggregated pairs of edges.
This multiplies the number of {\candidate} faults by 16 (6 ways to add one
edge, 9 ways to add two edges, plus the initial one).
The best indicator increases again, by~14\%;
it corresponds to the target fault (in green in
Figure~\ref{fig:one:fault}e).
The short list of {\candidate} faults is made of the~7 extended aggregated
pairs having an indicator lower than $\thetaext\bestIalfext$ (in light blue
in Figure~\ref{fig:one:fault}e).

After the minimization of the objective function with respect to the fault
parameter (remember that $\epsalf=\alf$) for all~7 {\candidate} faults of the
short list, the best data misfit is divided by a huge factor of~$10^{16}$
from its initial value.
This corresponds to the {\candidate} fault with the best indicator, which
happens to be the target fault.
And at the minimum, the target value $\alf=2$ is recovered
(Figure~\ref{fig:one:fault}f).
The algorithm stops on the convergence criterion.

\bigskip

\begin{figure}[\figloc]
  \begin{center}
    (a)
    \includegraphics[width=\twofigs]{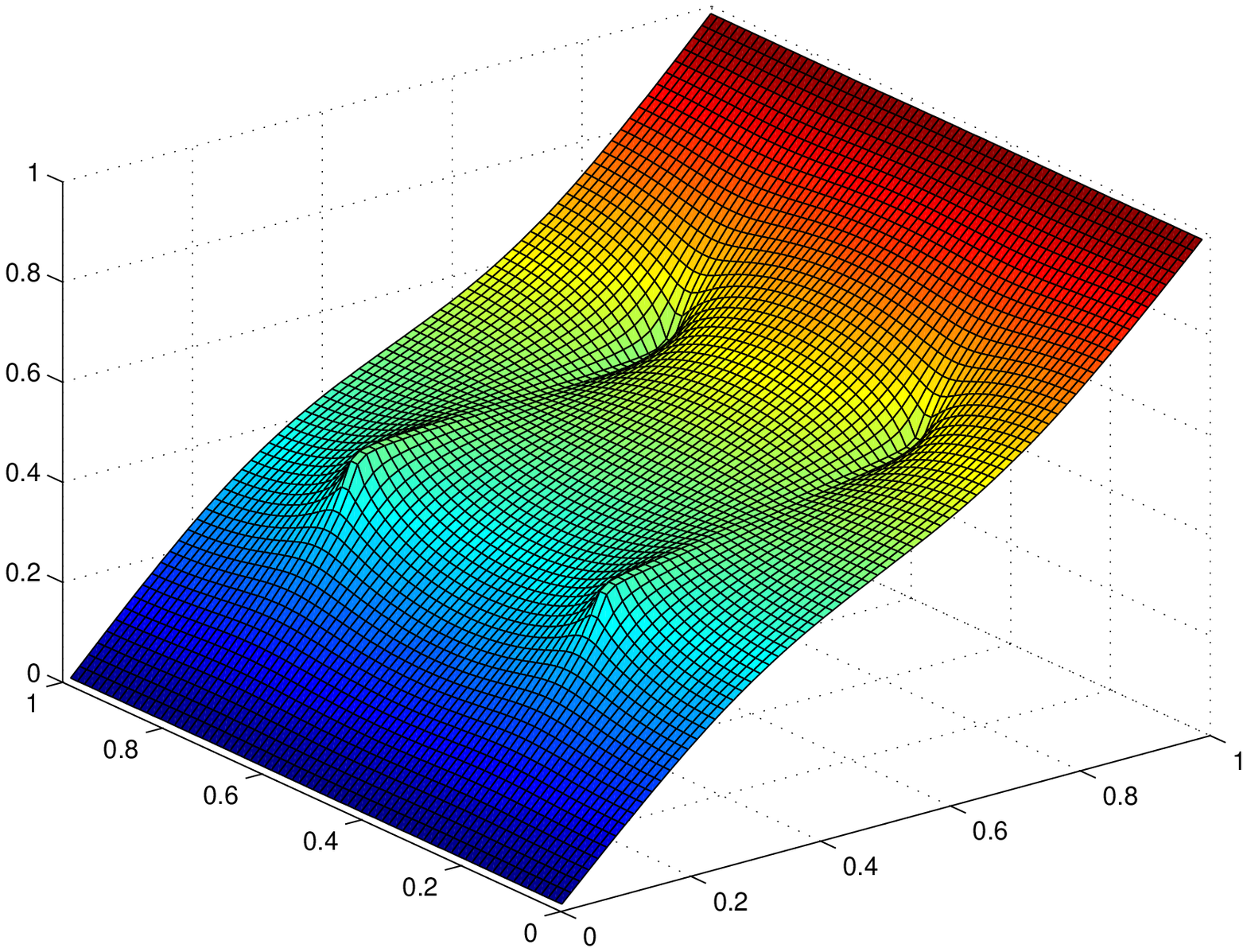}
    (d)
    \includegraphics[width=\twofigs]{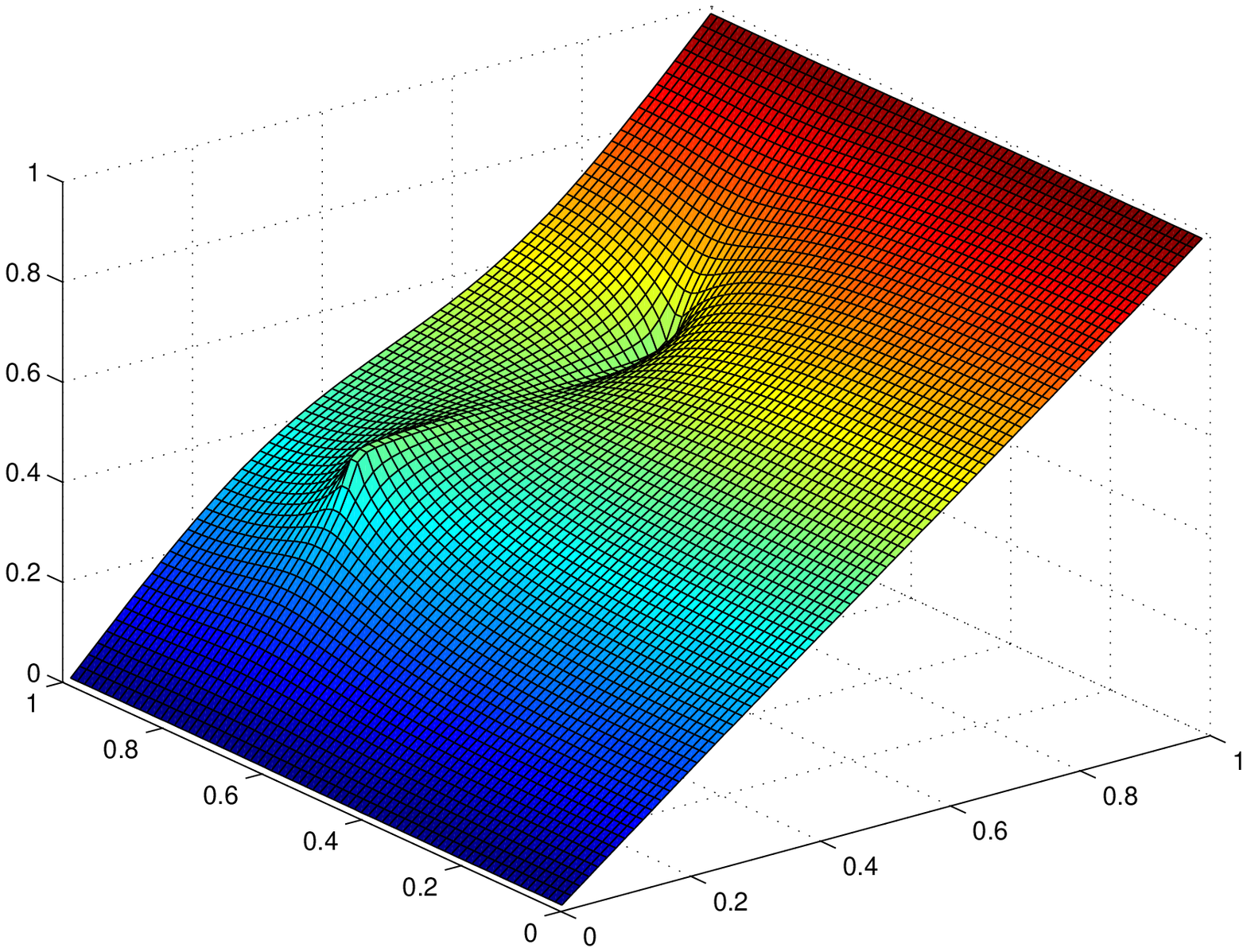} \\
    (b)
    \includegraphics[width=\twofigs]{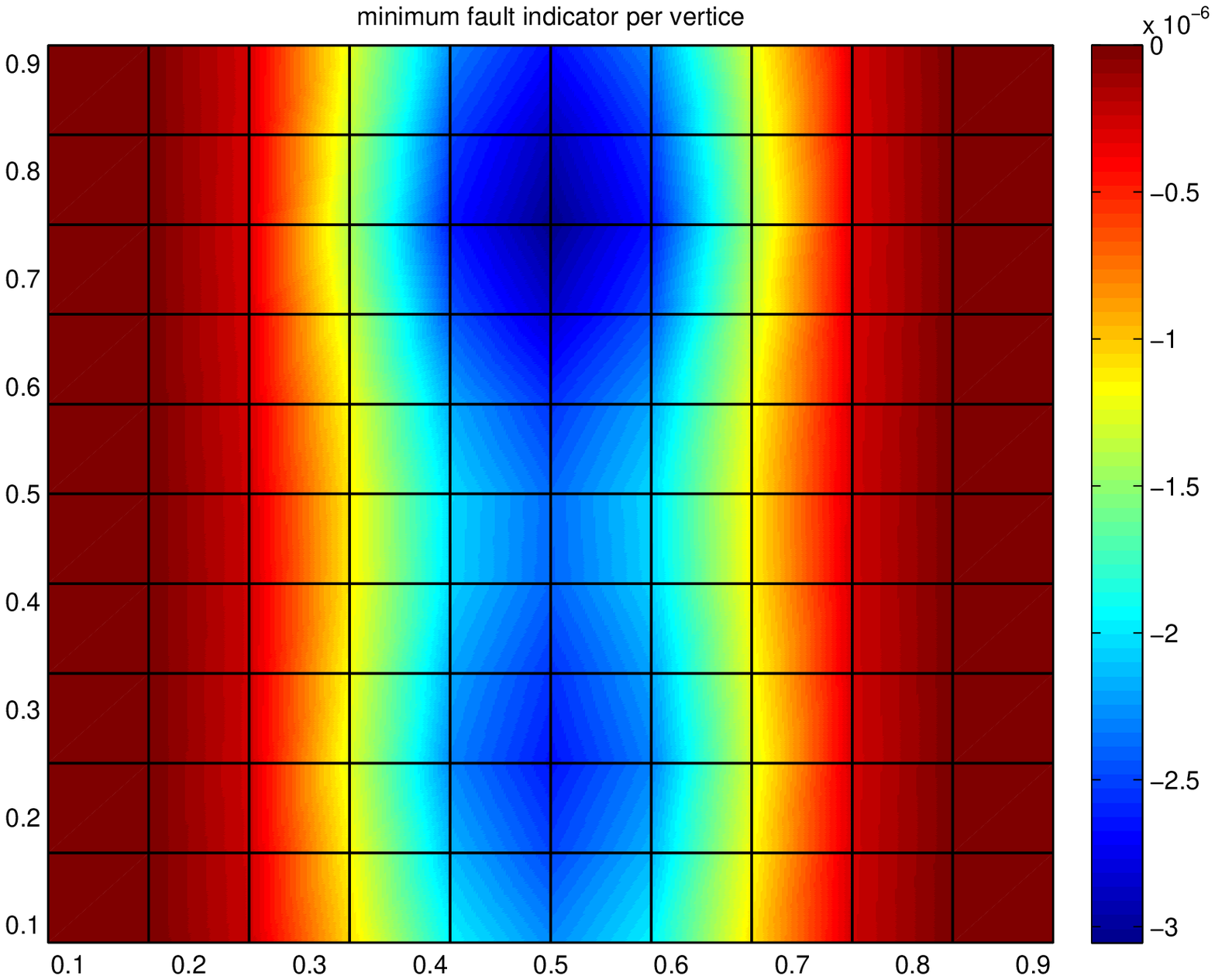}
    (e)
    \includegraphics[width=\twofigs]{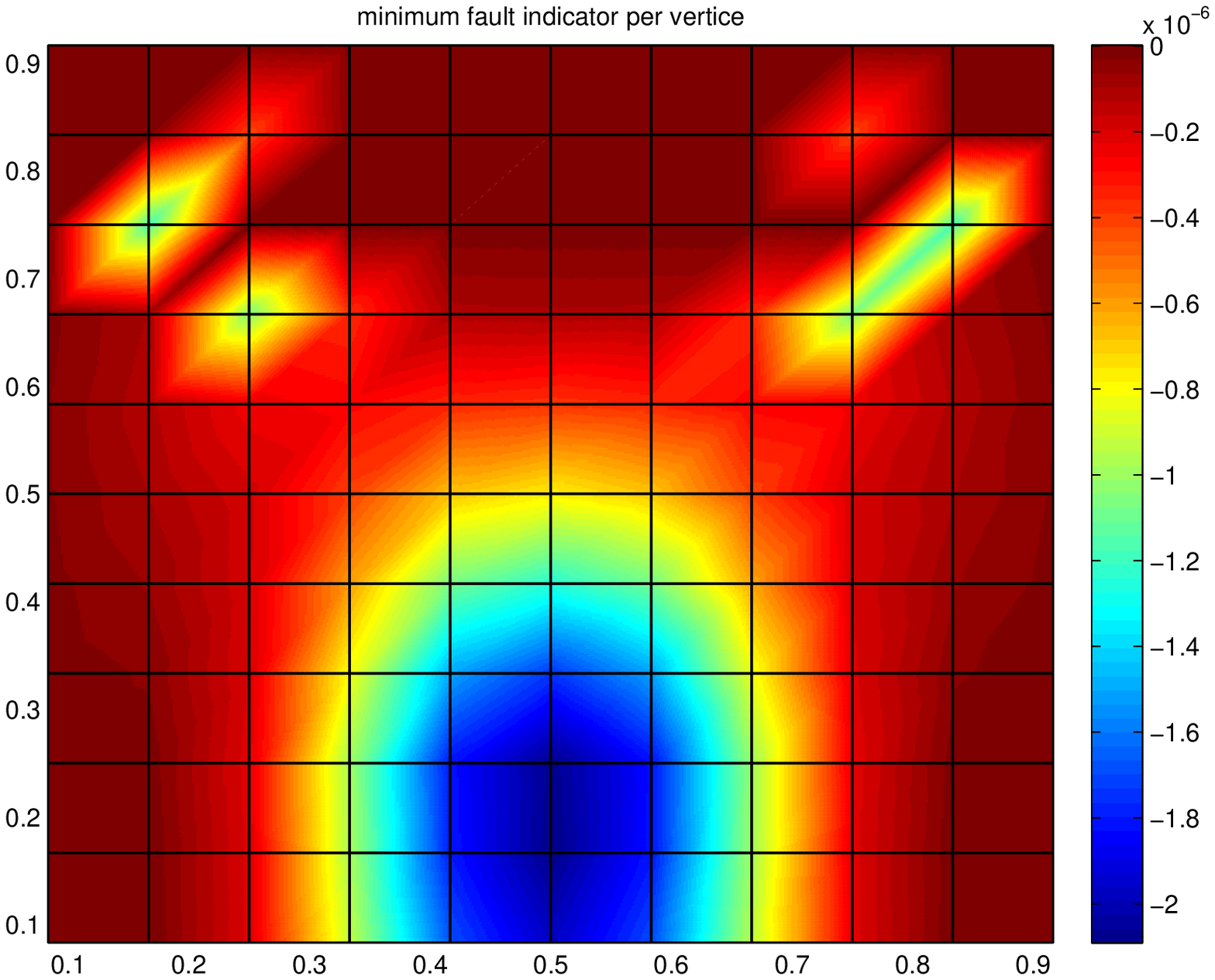} \\
    (c)
    \includegraphics[width=\twofigs]{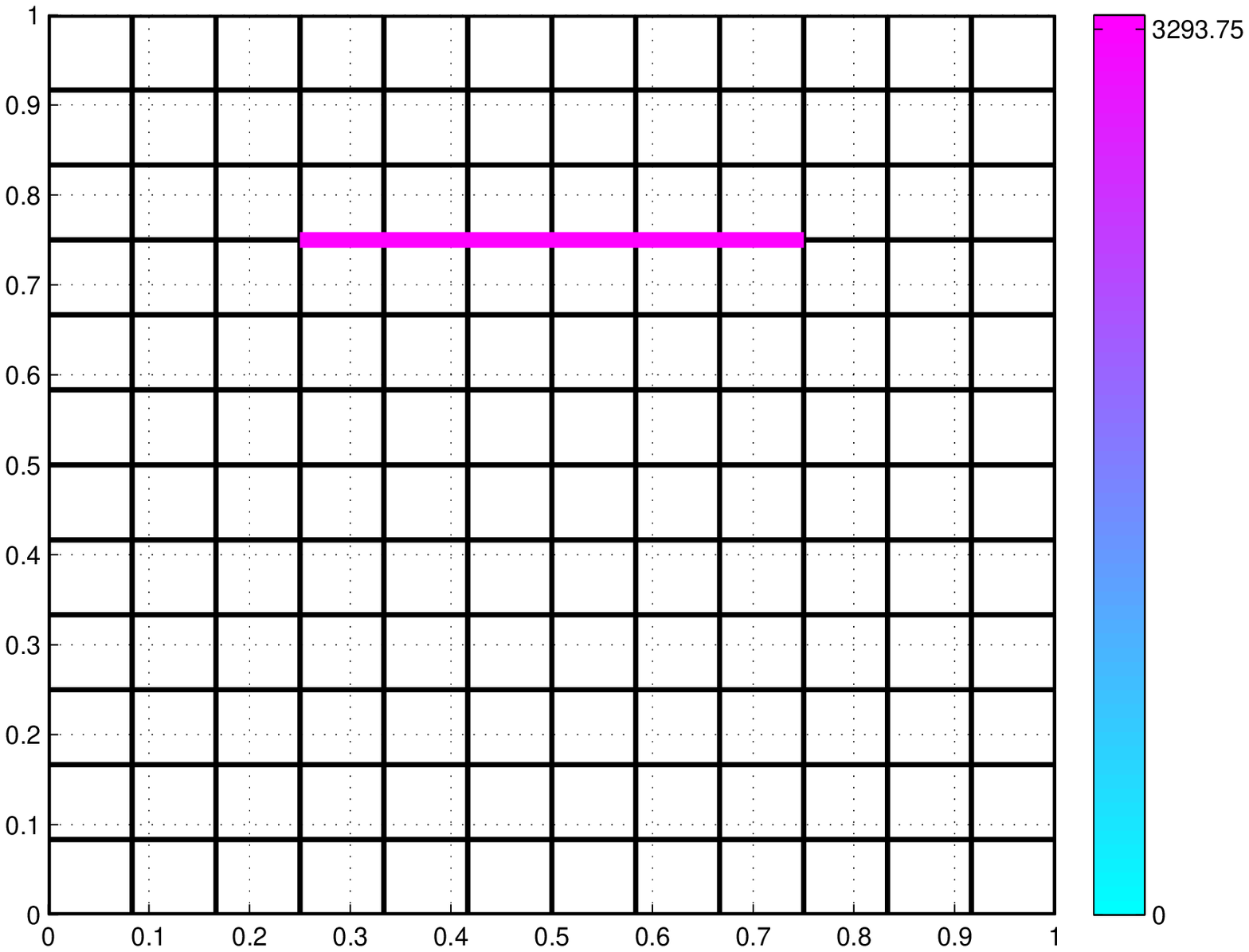}
    (f)
    \includegraphics[width=\twofigs]{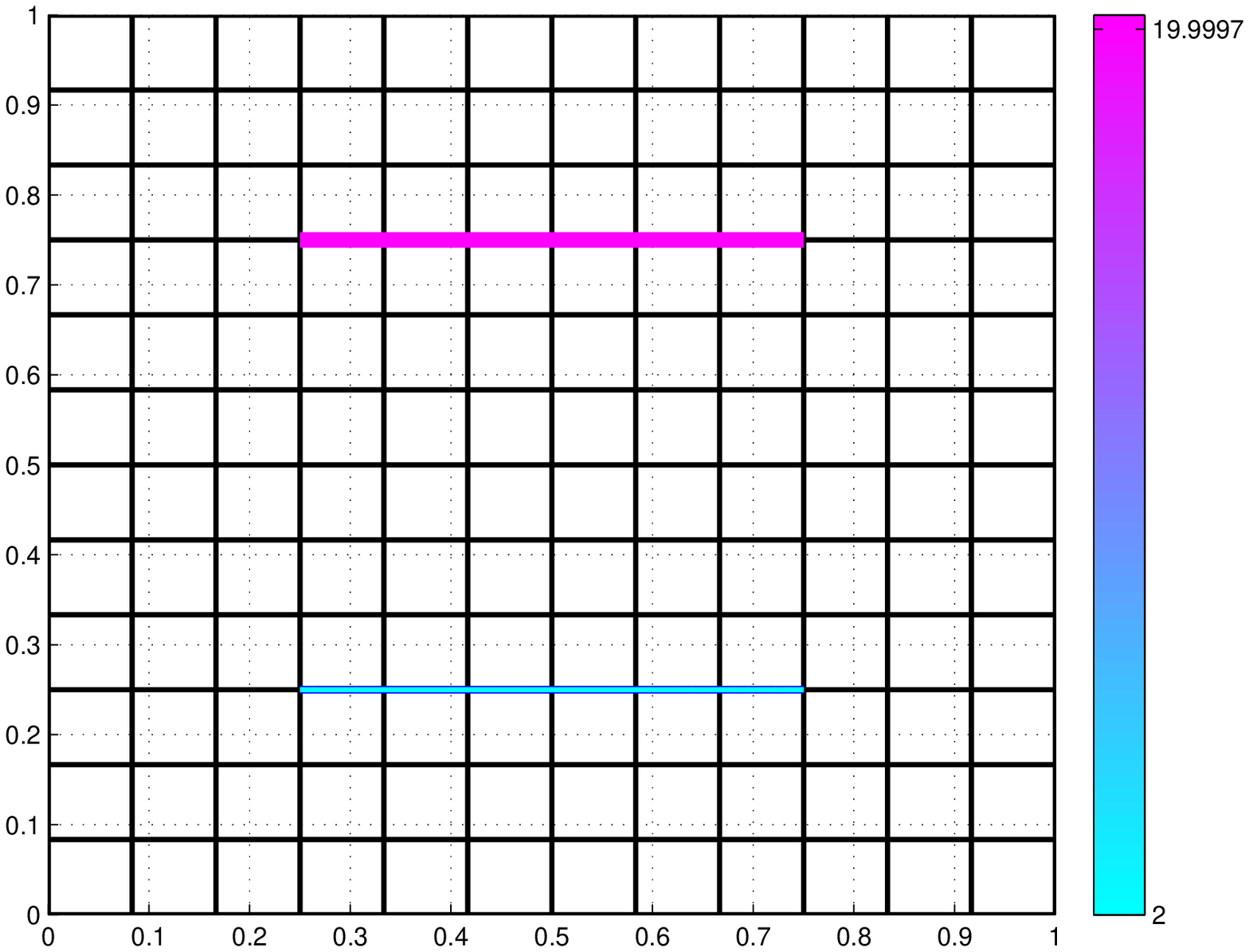}
    \caption{Case of two faults tangential to the flow.
      \newline
      (a) Data pressure, solution of the direct
      model~\eqref{eq:fracmodel:direct:discr:current} with two faults at
      $y=0.25$ ($\alf=2$) and $y=0.75$ ($\alf=20$).
      \newline
      (b) Distribution of indicators for the long list of elementary
      {\candidate} faults at the first iteration (see
      Figure~\ref{fig:one:fault}c).
      \newline
      (c) Best result after minimization at the first iteration for all
      {\candidate} faults of the short list.
      \newline
      (d) Pressure, solution of the direct model with the estimated fault
      of~(c).
      \newline
      (e) Distribution of indicators for the long list of elementary
      {\candidate} faults at the second iteration.
      \newline
      (f) Best result after minimization at the second iteration for all
      {\candidate} faults of the short list.
    }
    \label{fig:two:faults}
  \end{center}
\end{figure}

The picture is even more interesting with two faults.
Let us now consider the case of two faults tangential to the flow centered at
one quarter and at three quarters of the domain with $\alf=2$ (bottom) and~20
(top).
The synthetic data are represented in Figure~\ref{fig:two:faults}a.
The measurement and fracture meshes remain unchanged ($\Nm=72$ and $\Nf=12$).

At the first iteration, the cartography of indicators for all elementary
{\candidate} faults of the long list shows in Figure~\ref{fig:two:faults}b
two blue spots of highly negative values centered at the correct location of
the target faults.
Moreover, the best indicators are located around the upper fault having the
highest permeability (in dark blue).
As in the case of a single fault, the aggregation and extension stages
multiply the best indicator by a factor of~2, and the best indicator
corresponds to the most permeable target fault (the upper one).
Again, the best minimizer among the~9 {\candidate} faults of the short list
corresponds to the best indicator (the upper target fault).
The optimization step decreases the objective function by a factor of~3 only,
and at the minimum, a highly overestimated value of about~3300 is recovered
for~$\alf$ (Figure~\ref{fig:two:faults}c).
Yet, the pressure of Figure~\ref{fig:two:faults}d recovers the main
perturbation due to the upper target fault (compare with
Figure~\ref{fig:two:faults}a).
As the stopping criteria are not satisfied, another iteration is performed.

At the second iteration, the cartography of indicators presents now in
Figure~\ref{fig:two:faults}e a strong dark blue spot located at the center of
the least permeable target fault.
Clearly, there is nothing left to do around the upper target fault (the area
is mostly dark red).
The indicators are smaller in magnitude than in the first iteration by a
factor of about one third.
The optimization step exhibits the best {\candidate} fault of the
short list (among 10 of them) as the best minimizer.
Again, this best minimizer is exactly the missing lower target fault, the
data misfit is divided by a huge factor of~$2\,10^{13}$ from its initial
value, and the target values $\alf=20$ and $\alf=2$ are perfectly recovered
(Figure~\ref{fig:two:faults}f).
The algorithm stops again on the convergence criterion.

\myclearpage
\subsection{Beyond the inverse crime}
\label{ss:num:noise}

We investigate now less favorable situations with less measurements, and
where the synthetic data are no longer in the range of the model used for
inversion.
The comprehensive campaign of numerical tests is detailed
in~\cite{these:fatma:2016}.

\begin{figure}[\figloc]
  \begin{center}
    (a)
    \includegraphics[width=\twofigs]{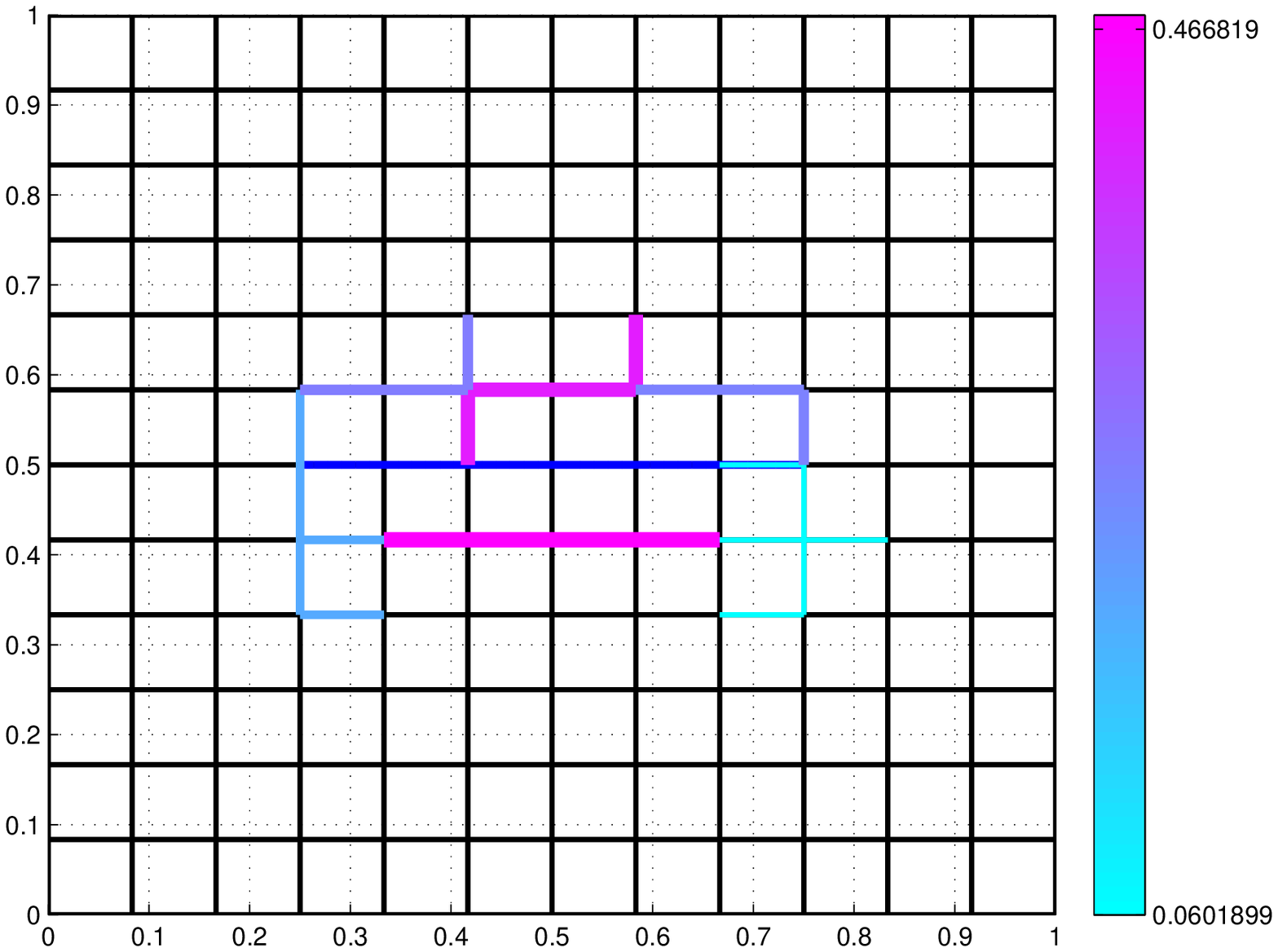}
    (b)
    \includegraphics[width=\twofigs]{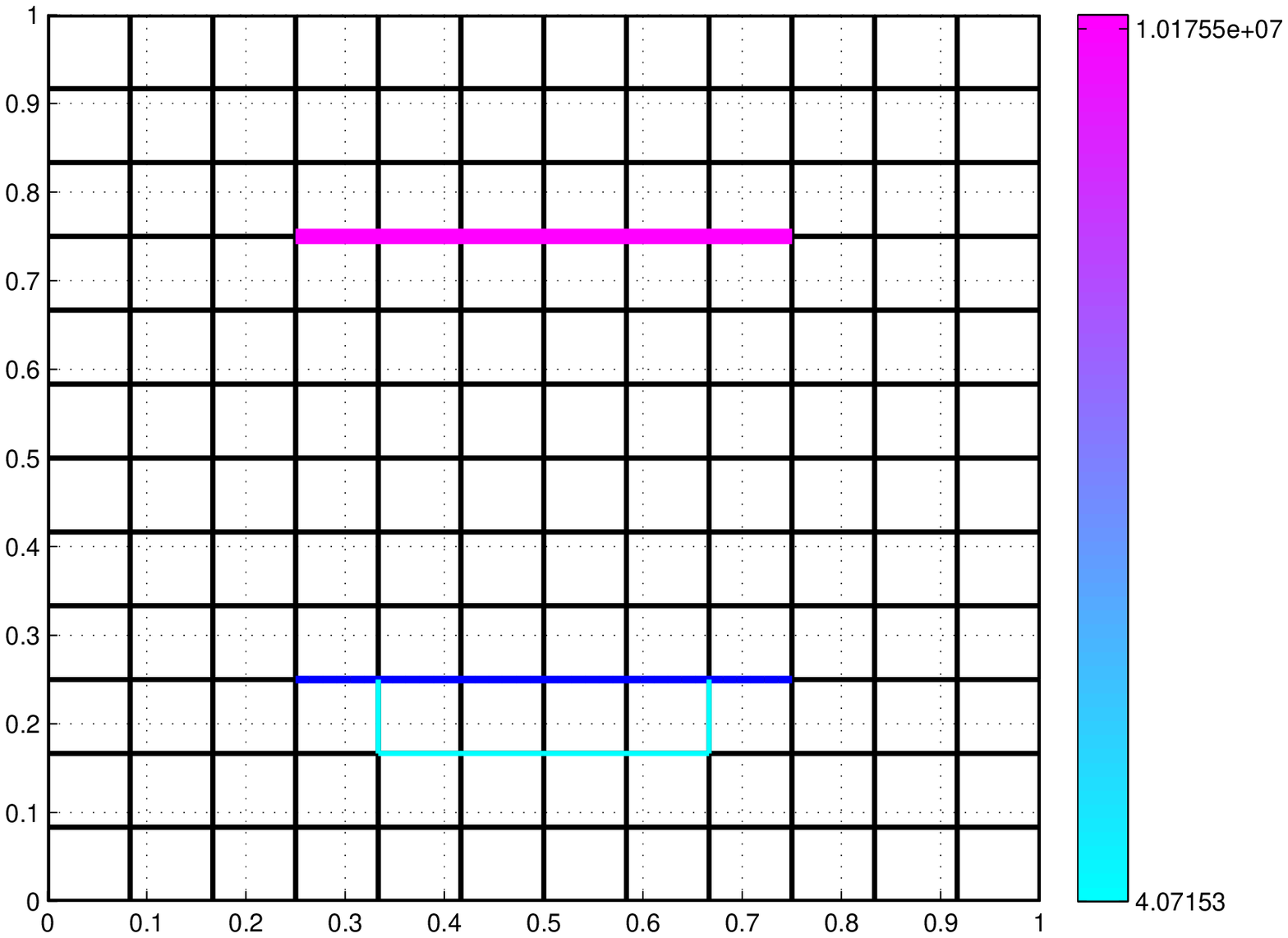} \\
    (c)
    \includegraphics[width=\twofigs]{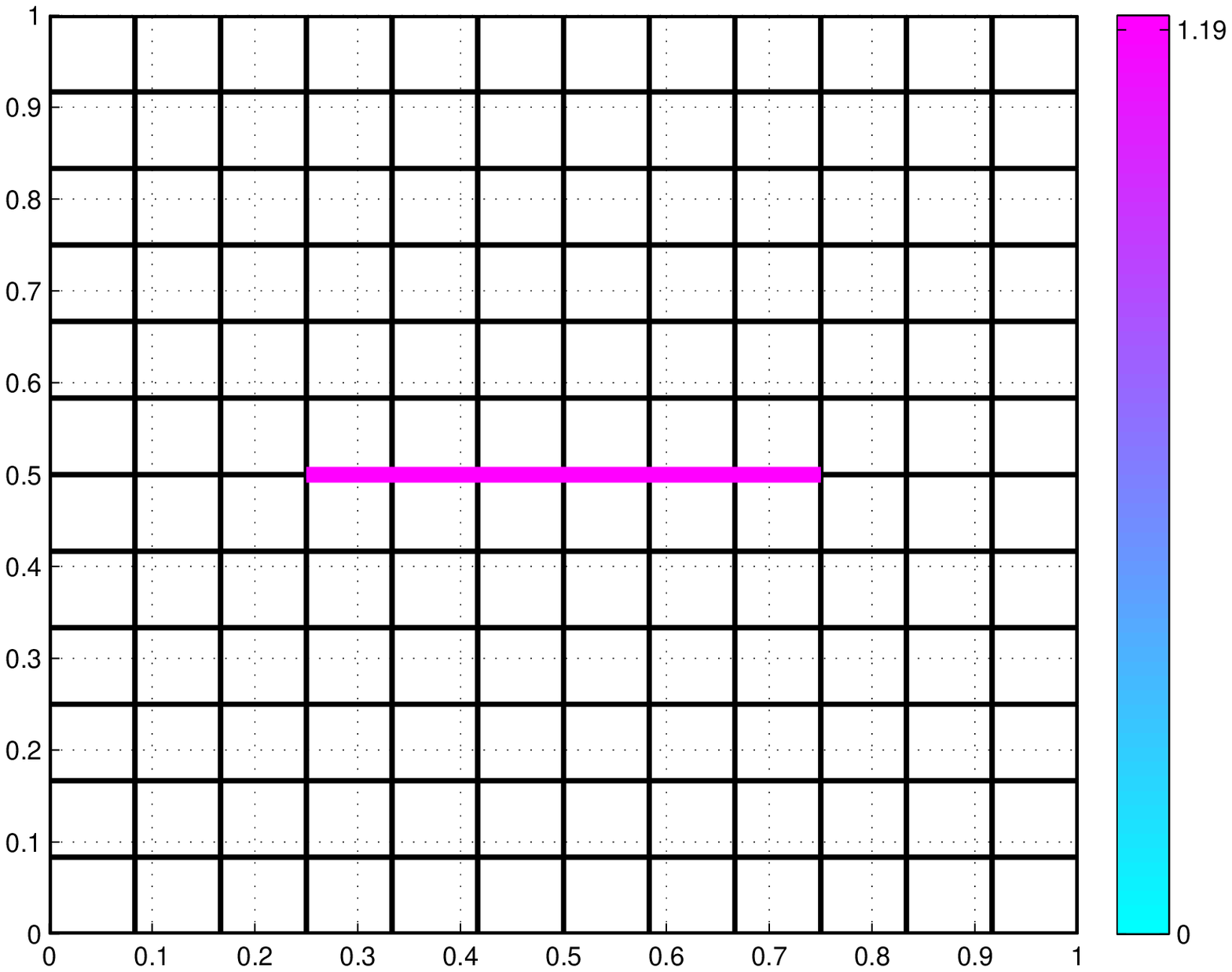}
    (d)
    \includegraphics[width=\twofigs]{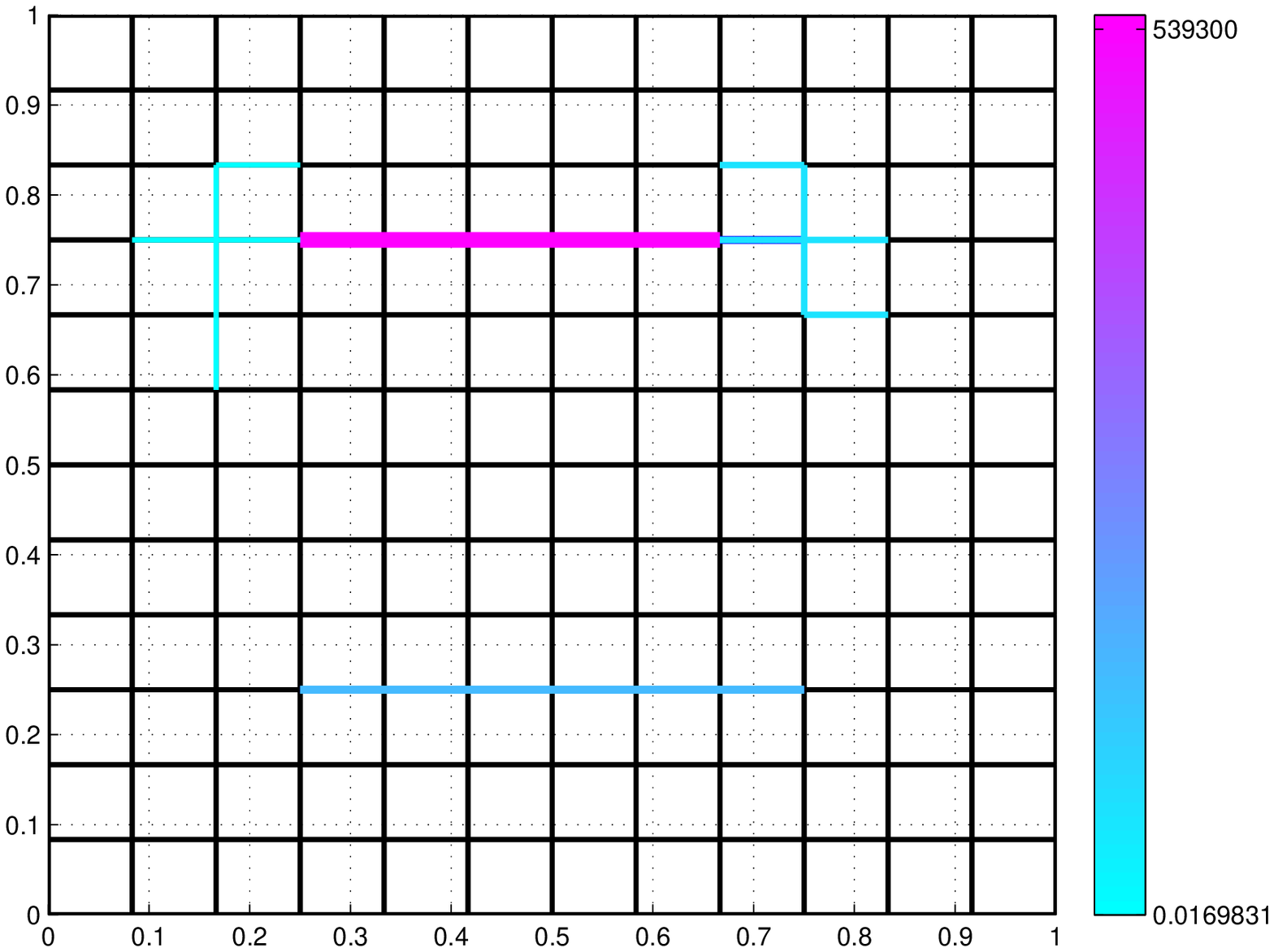} \\
    (e)
    \includegraphics[width=\twofigs]{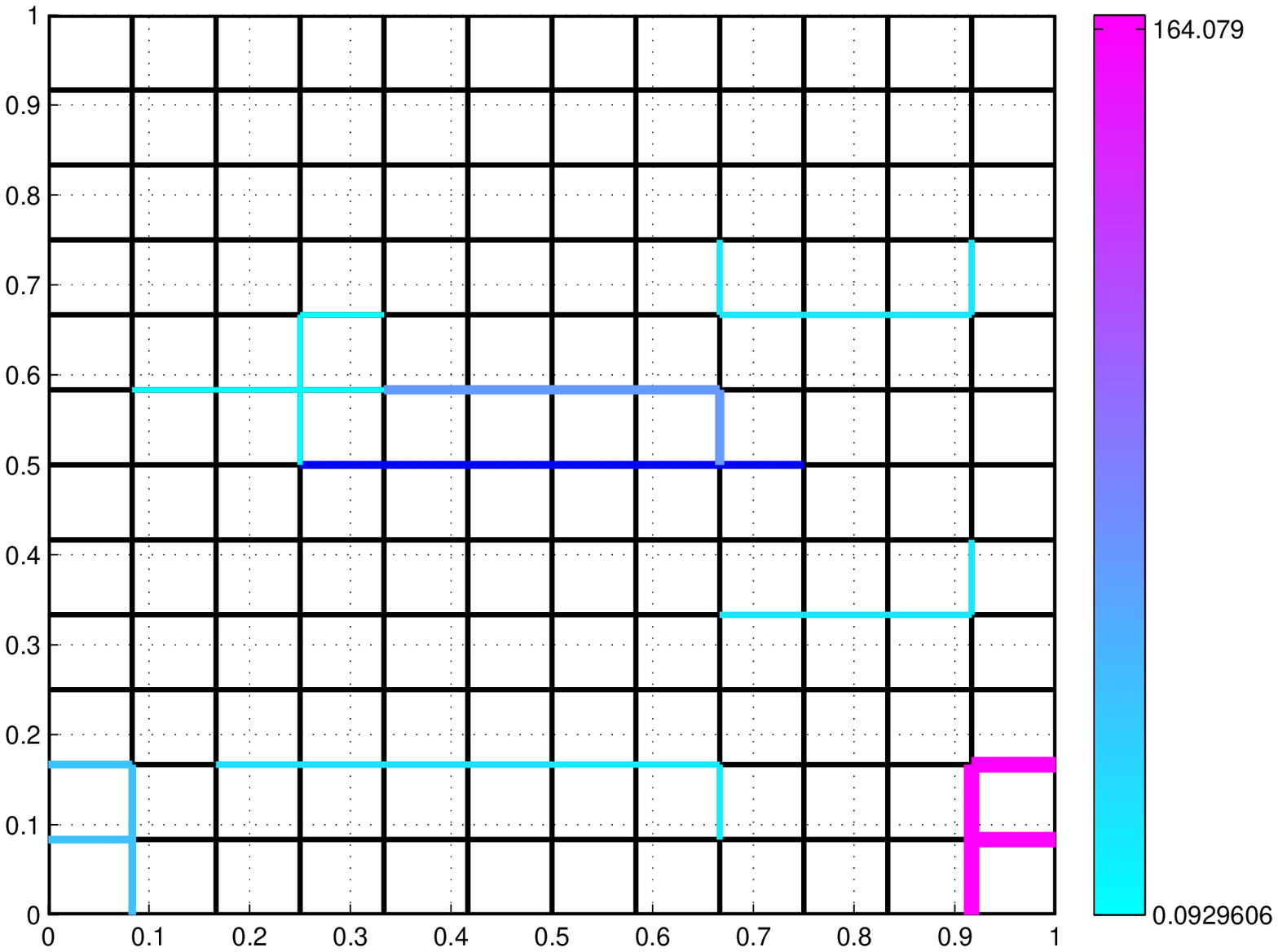}
    (f)
    \includegraphics[width=\twofigs]{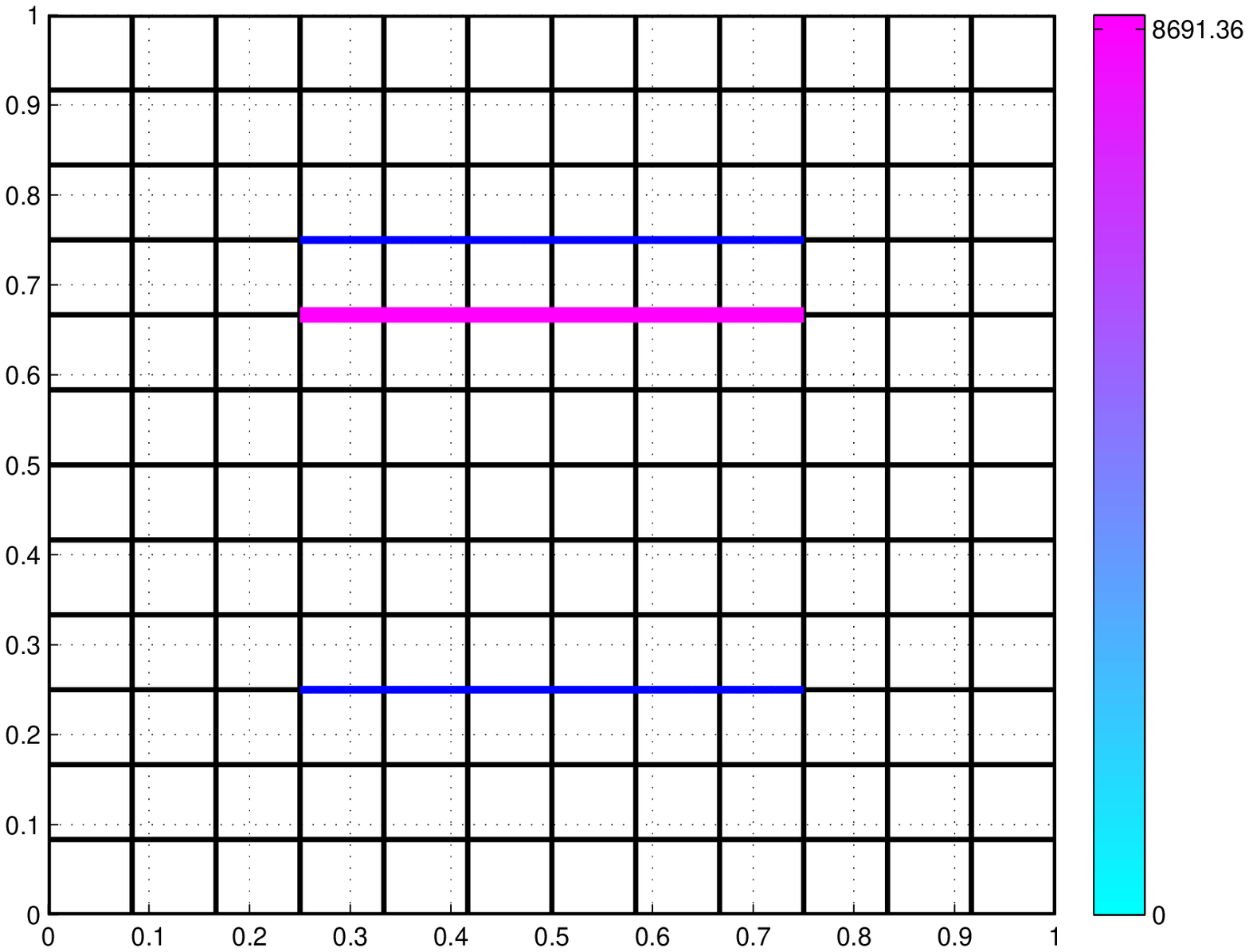}
    \caption{Estimation of faults with less measurement data.
      The inverted faults are represented with colors ranging from light blue
      (low permeabilities) to pink (high permeabilities).
      The target fault is represented in dark blue.
      \newline
      (a) Single target fault, no noise, $\Nm=6$.
      (b) Two target faults, no noise, $\Nm=8$.
      \newline
      (c) Single target fault, 4\% noise, $\Nm=8$.
      (d) Two target faults, 4\% noise, $\Nm=12$.
      \newline
      (e) Single target fault, 6\% noise, $\Nm=8$.
      (f) Two target faults, 2\% noise, $\Nm=8$.
    }
    \label{fig:fault:less:data}
  \end{center}
\end{figure}

\myparagraph{Influence of decimation of measurements}
When the measurement mesh~$\Thm=\Quad{\Nm}{\Nm}$ is reduced from $\Nm=72$
down to $\Nm=8$, numerical tests shown in~\cite{these:fatma:2016} indicate
that the algorithm still provides the target fault as the best {\candidate}
fault, and the minimization step still yields the target value of~$\alf$.
Moreover, this remains true for target values of~$\alf$ ranging from~0.2 up
to~200.
For values above~200, the direct model outputs almost the same pressure
field: the cost function becomes insensitive to the fault parameter~$\alf$.
It is the same for the case of two target faults: down to $\Nm=12$, the
algorithm still provides the most permeable target fault as the best
{\candidate} fault at the first iteration, the least permeable target fault
at the second iteration, and then, the minimization step yields the target
values of~$\alf$ for both faults.

Problems appear when the distribution of measurement points becomes too
loose.
At best, the target faults are still located quite accurately, but the
minimization step is unable to recover the target values of~$\alf$, and the
results may be misevaluated by several orders of magnitude as in the case of
two target faults with $\Nm=8$ (Figure~\ref{fig:fault:less:data}b).
Or the situation may degenerate progressively: the recovered locations become
quite diffuse, but the effective permeabilities can still be quite as,
accurate in the case of a single target fault with $\Nm=6$
(Figure~\ref{fig:fault:less:data}a).

\myparagraph{Influence of random noise}
When adding white Gaussian noise up to a relative noise level of~6\% (the
perturbation due to the presence of the fault becomes hardly visible by the
eye), the algorithm still exhibits the target fault as the best {\candidate}
fault, and the minimization step recovers a value of $\alf=2.01$, almost
identical to the target value $\alf=2$, see~\cite{these:fatma:2016}.
When also decimating the measurements down to~$\Nm=8$, the algorithm still
recovers the target fault as the best {\candidate} fault up to a relative
noise level of 4\%, and the minimization step yields promising values of~1.5
(with~2\% noise) and~1.2 (with~4\% noise) for the fault parameter~$\alf$
(Figure~\ref{fig:fault:less:data}c).
The algorithm is not as successful with~6\% noise and $\Nm=8$: a reasonable
estimated fault in the middle (in medium blue) is polluted by a very strong
fault in the lower right corner of the domain
(Figure~\ref{fig:fault:less:data}e).

When considering the case of two faults, the conditions are still favorable
up to a noise level of~4\%: the location of the most permeable fault is
recovered at the first iteration and that of the other at the second
iteration.
Moreover, the minimization step recovers the value $\alf=17$ for the upper
fault and $\alf=1.99$ for the lower fault, which is still close to the target
values~20 and~2.
When also decimating the measurements down to $\Nm=12$, the algorithm still
recovers correctly the least permeable (lower) fault with an excellent
inverted value of~$\alf=1.98$, but the most permeable (upper) fault is highly
overestimated (values over half a million, Figure~\ref{fig:fault:less:data}d)
and with~4\% noise, the location of this fault exhibits artifacts associated
with negligible permeabilities under~0.1.
For $\Nm=8$, only the strongest (upper) fault is estimated with~2\% noise, it
is located one cell under the target location and the permeability is highly
overvalued to about~8700 (Figure~\ref{fig:fault:less:data}f).

\myclearpage

\begin{figure}[\figloc]
  \begin{center}
    (a)
    \includegraphics[width=\twofigs]{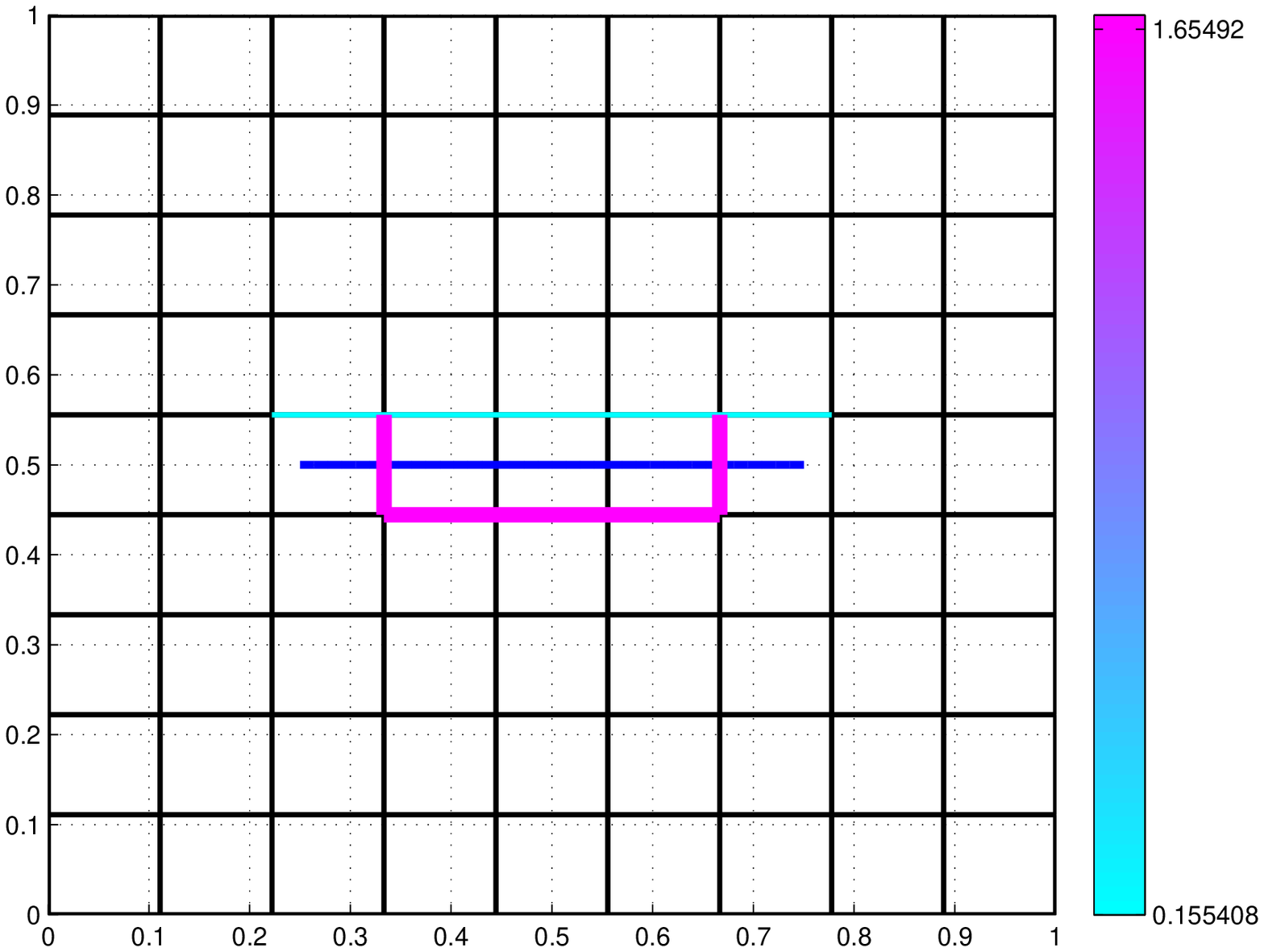}
    (b)
    \includegraphics[width=\twofigs]{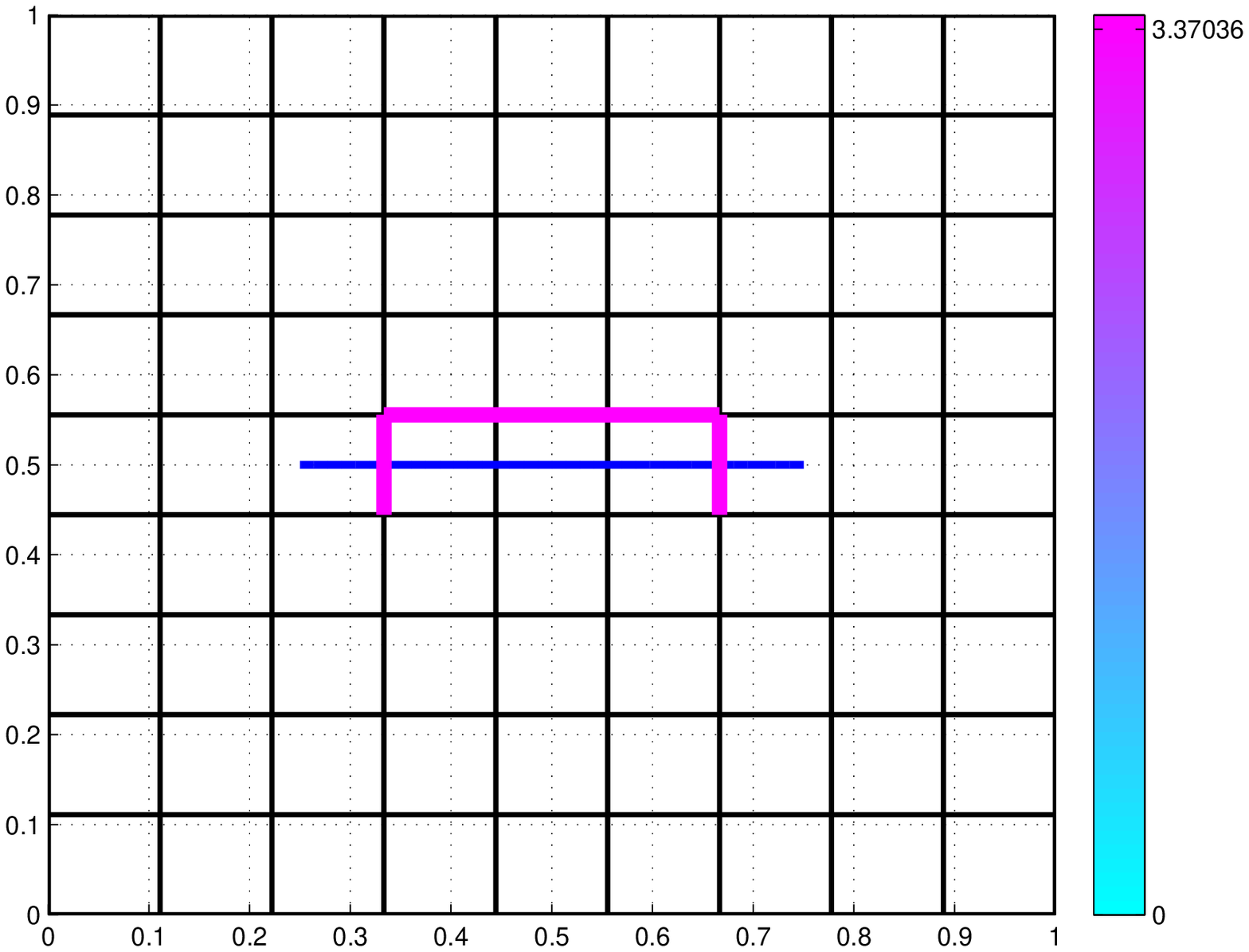} \\
    (c)
    \includegraphics[width=\twofigs]{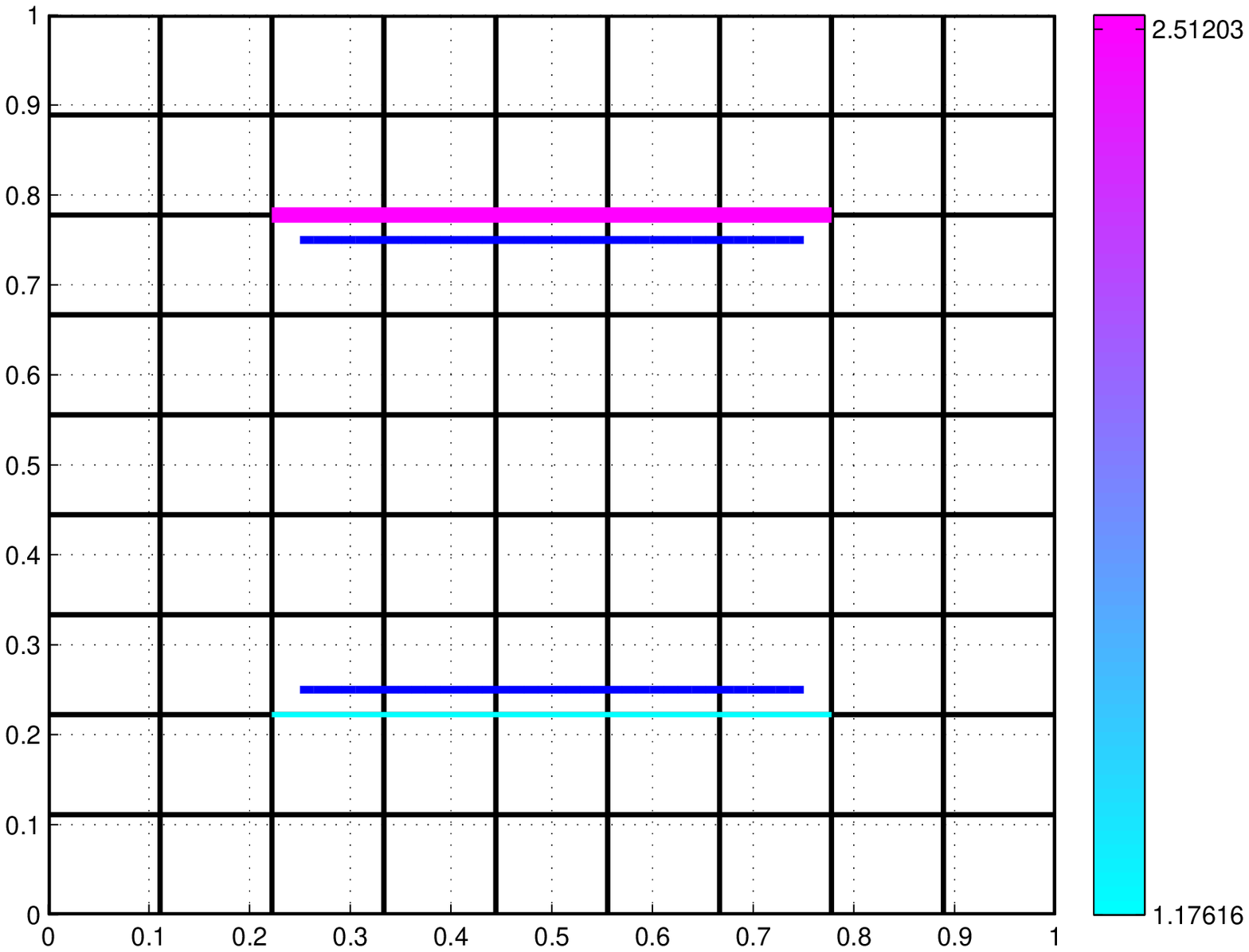}
    (d)
    \includegraphics[width=\twofigs]{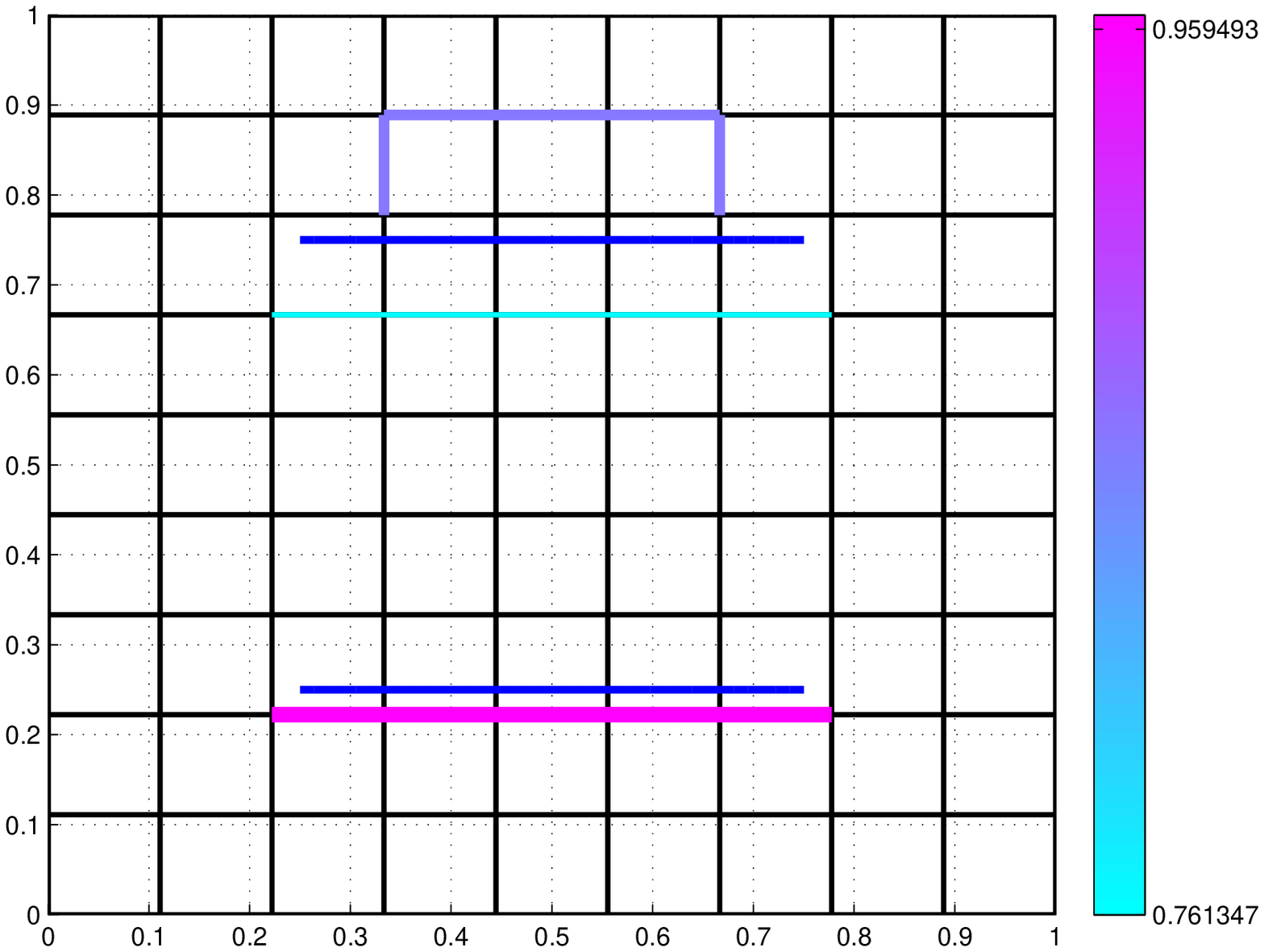}
    \caption{Estimation of faults with coherent noise:
      the target faults (dark blue) are not carried by the fracture research
      mesh ($\Nf=9$).
      The inverted faults are represented with colors ranging from light blue
      (low permeabilities) to pink (high permeabilities).
      \newline
      (a) Single target fault, $\Nm=72$.
      (b) Single target fault, $\Nm=8$
      \newline
      (c) Two target faults, $\Nm=72$.
      (d) Two target faults, $\Nm=8$.
    }
    \label{fig:fault:noise}
  \end{center}
\end{figure}

\myparagraph{Influence of coherent noise}
Finally, we investigate the much more difficult situation where the synthetic
pressure data do not belong to the range of the direct model used for
inversion, {\ie} the global minimum of the objective function is not~0.
Of course, we can no longer expect to perfectly retrieve the location and
permeability of the target faults, but at least we can hope to recover
groups of faults with similar hydrogeological signatures; in fact, values
of~$\alf$ are adding up for nearby parallel faults tangential to the flow.

We consider $\Nf=9$, hence the target faults located in the middle, or at one
quarter and three quarters of the domain, are no longer carried by the
fracture research mesh anymore.
We also decimate the pressure data from $\Nm=72$ down to $\Nm=8$.
The estimated locations of the faults contain most of the coarse edges of the
fracture research mesh that are closest to the target faults.
Moreover, with the single target fault, the recovered effective
permeabilities range from~1.8 to~3.4, for a target value of~2
(Figures~\ref{fig:fault:noise}a and~\ref{fig:fault:noise}b).
And with the two target faults, the recovered effective permeabilities are
higher for the upper fault than for the lower, but they are underestimated up
to an order of magnitude (Figures~\ref{fig:fault:noise}c
and~\ref{fig:fault:noise}d): from~1.6 to~2.5 for the most permeable fault
(the target value is~20), and from~1.0 to~1.6 for the least one (the target
value is~2).

\begin{figure}[\figloc]
  \begin{center}
    (a)
    \includegraphics[width=\threefigs]{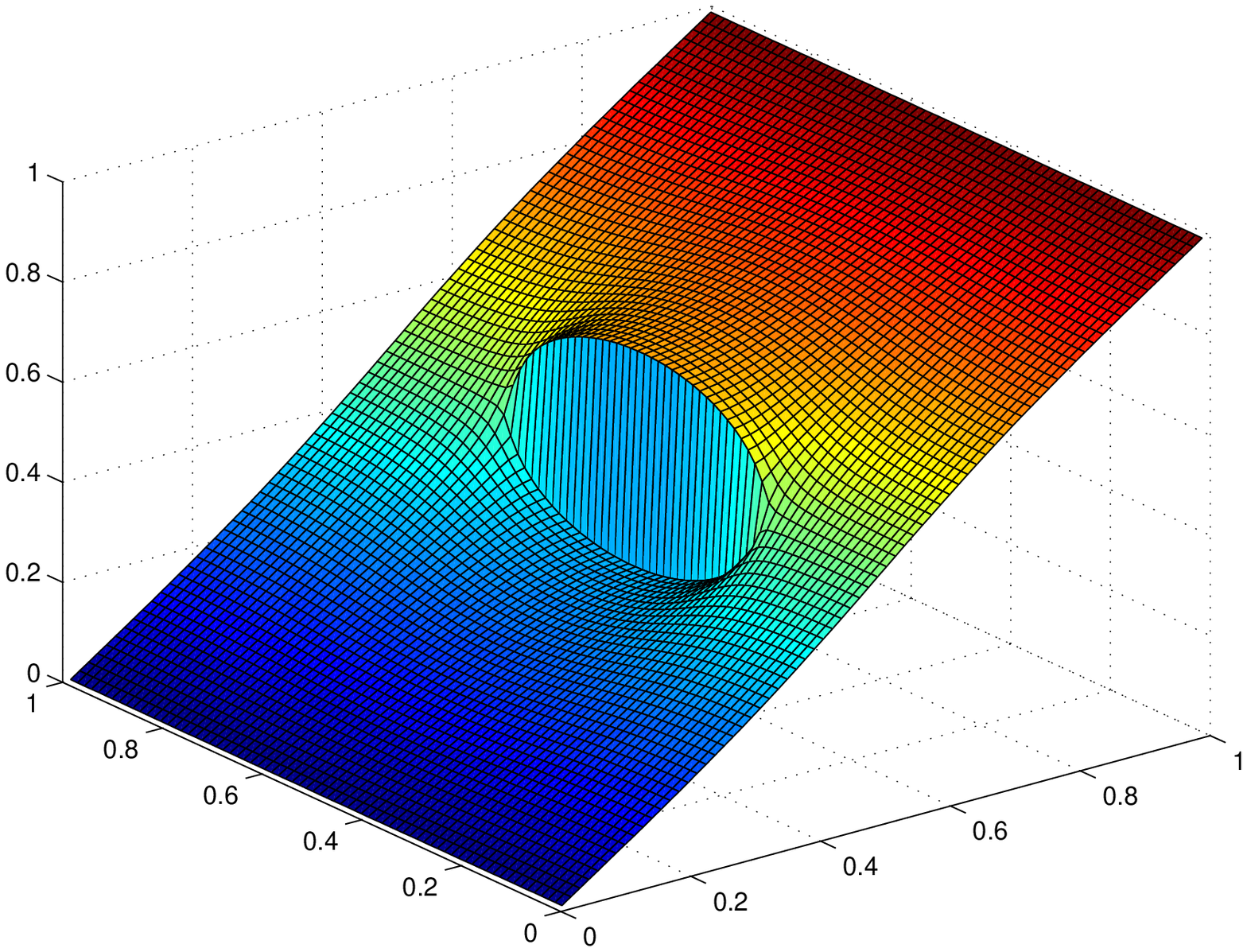}
    (b)
    \includegraphics[width=\threefigs]{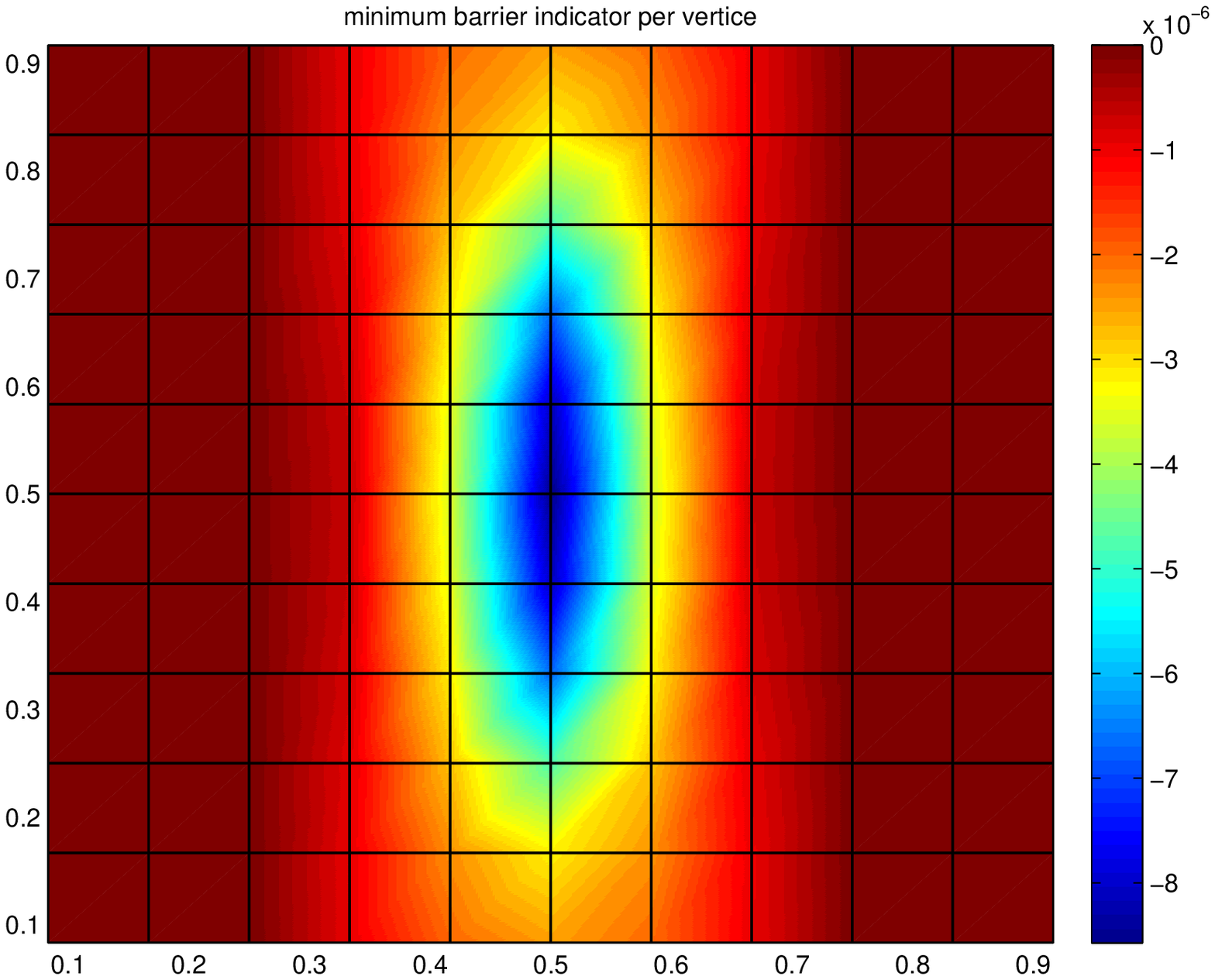}
    (c)
    \includegraphics[width=\threefigs]{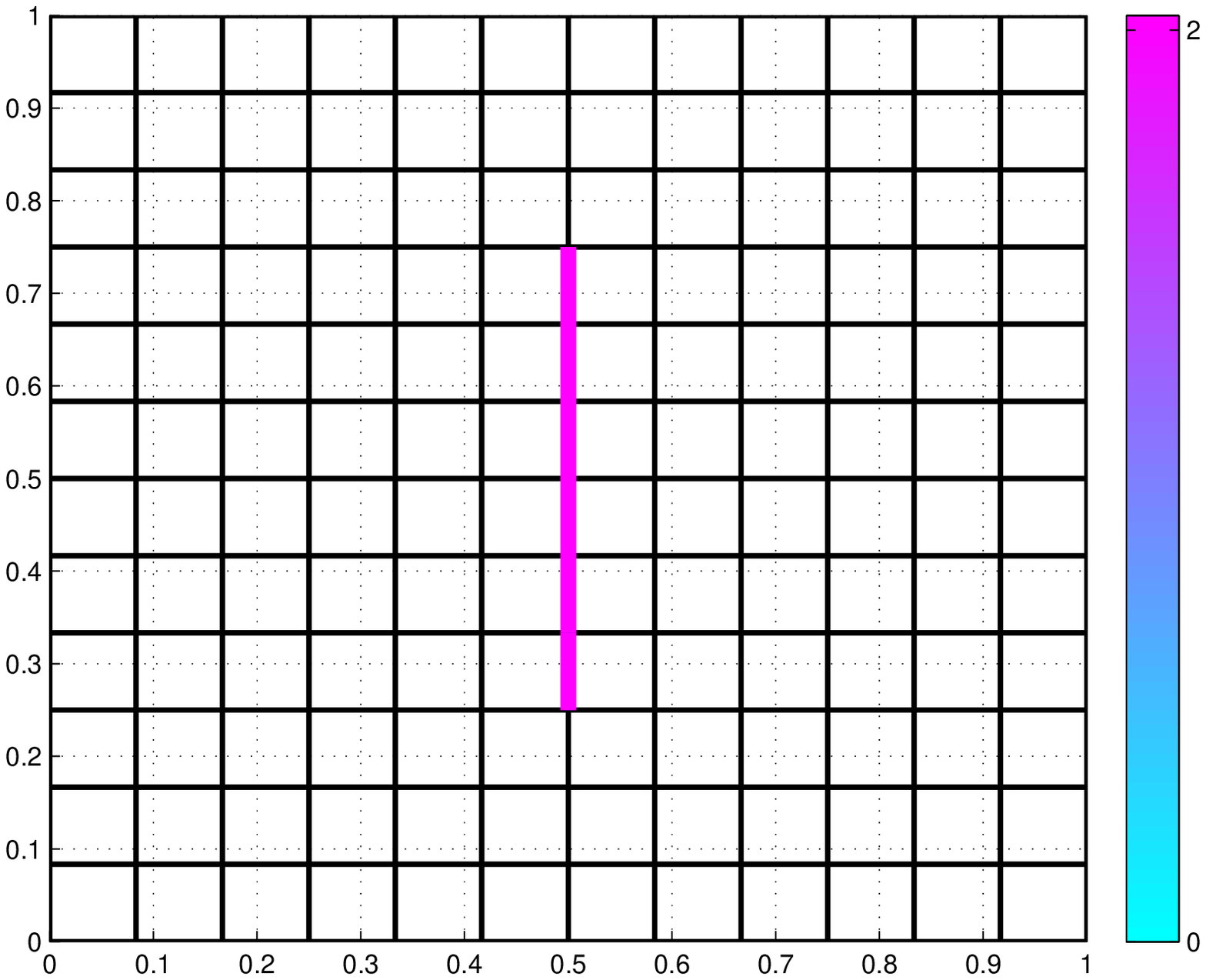} \\
    (d)
    \includegraphics[width=\twofigs]{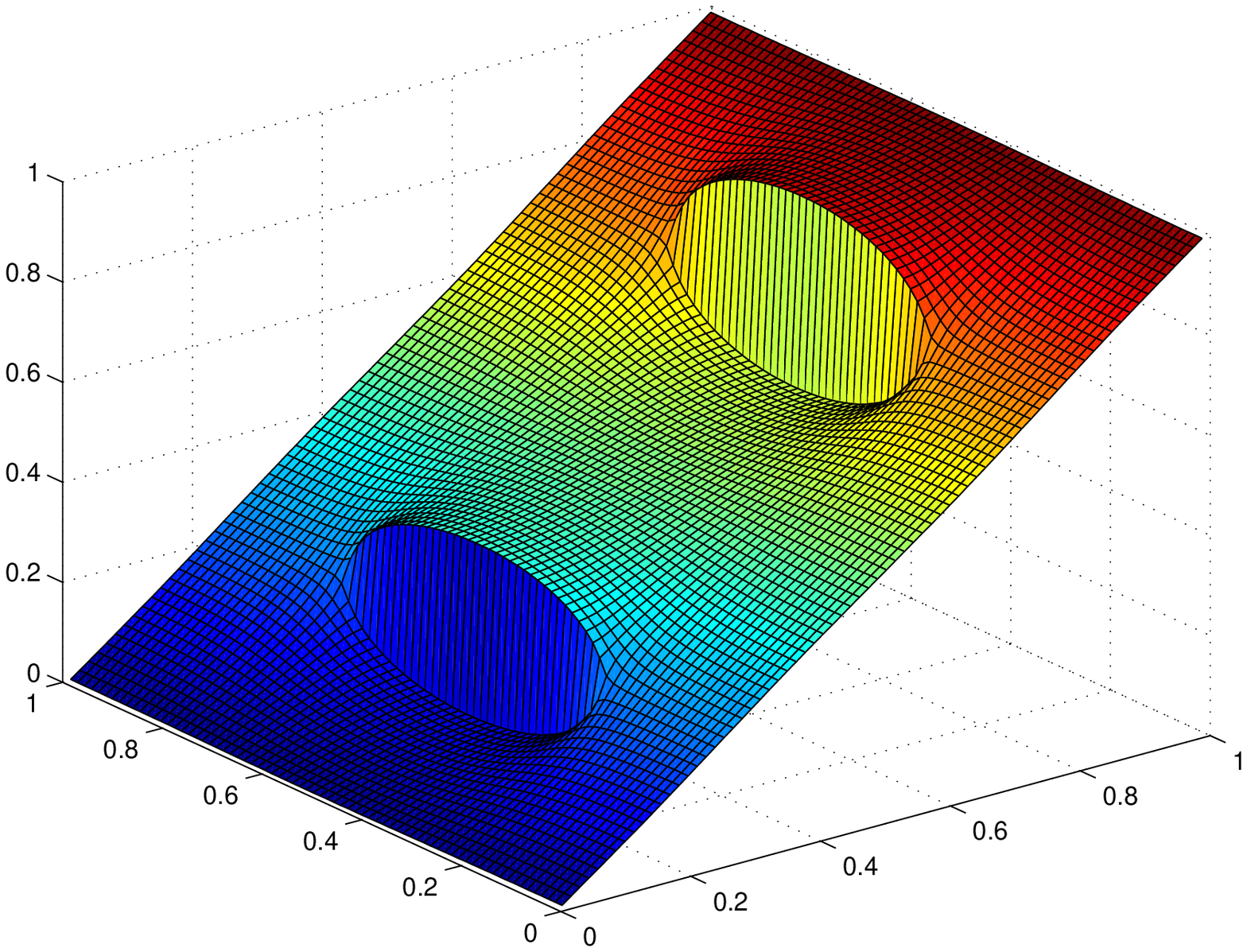}
    (e)
    \includegraphics[width=\twofigs]{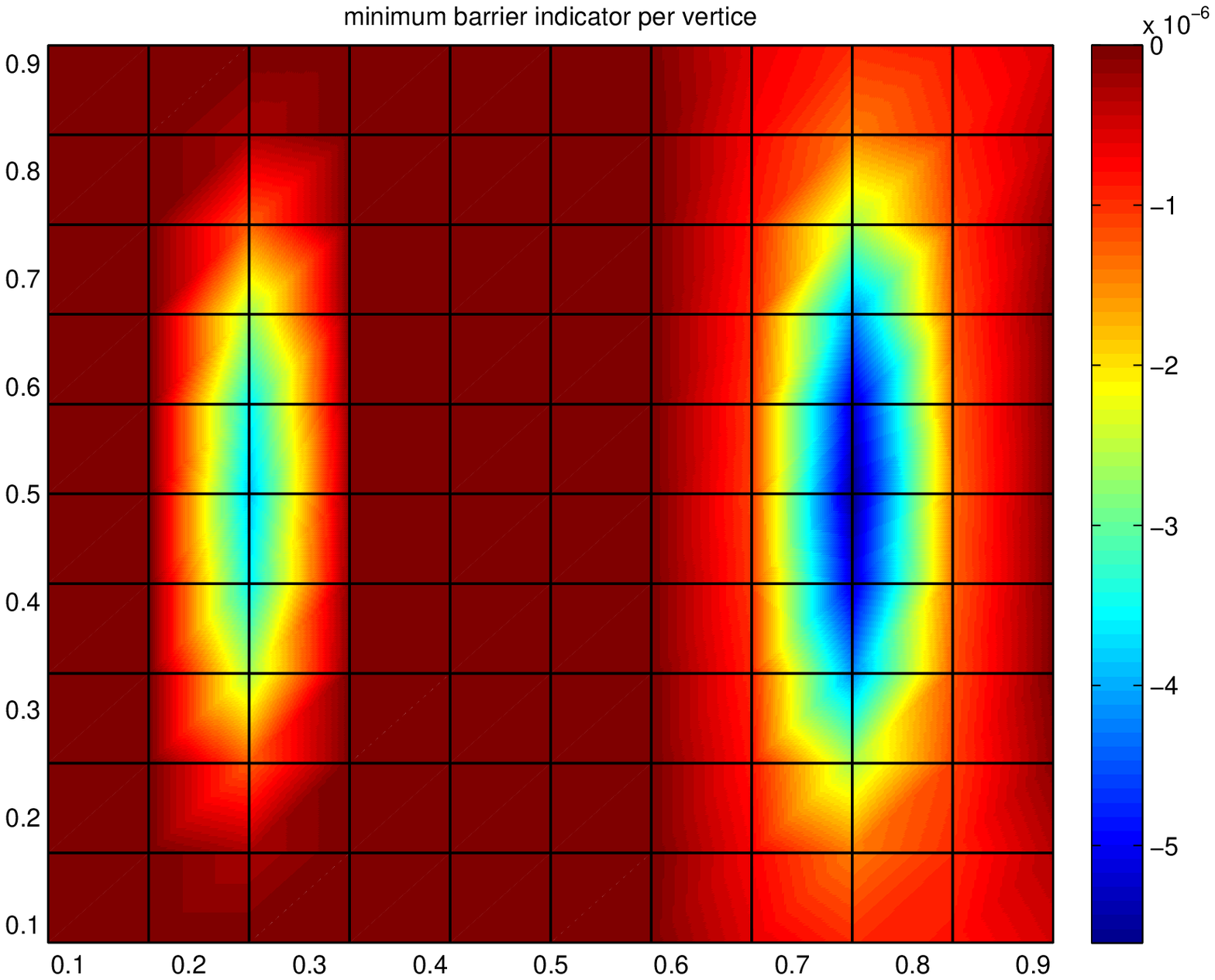} \\
    (f)
    \includegraphics[width=\twofigs]{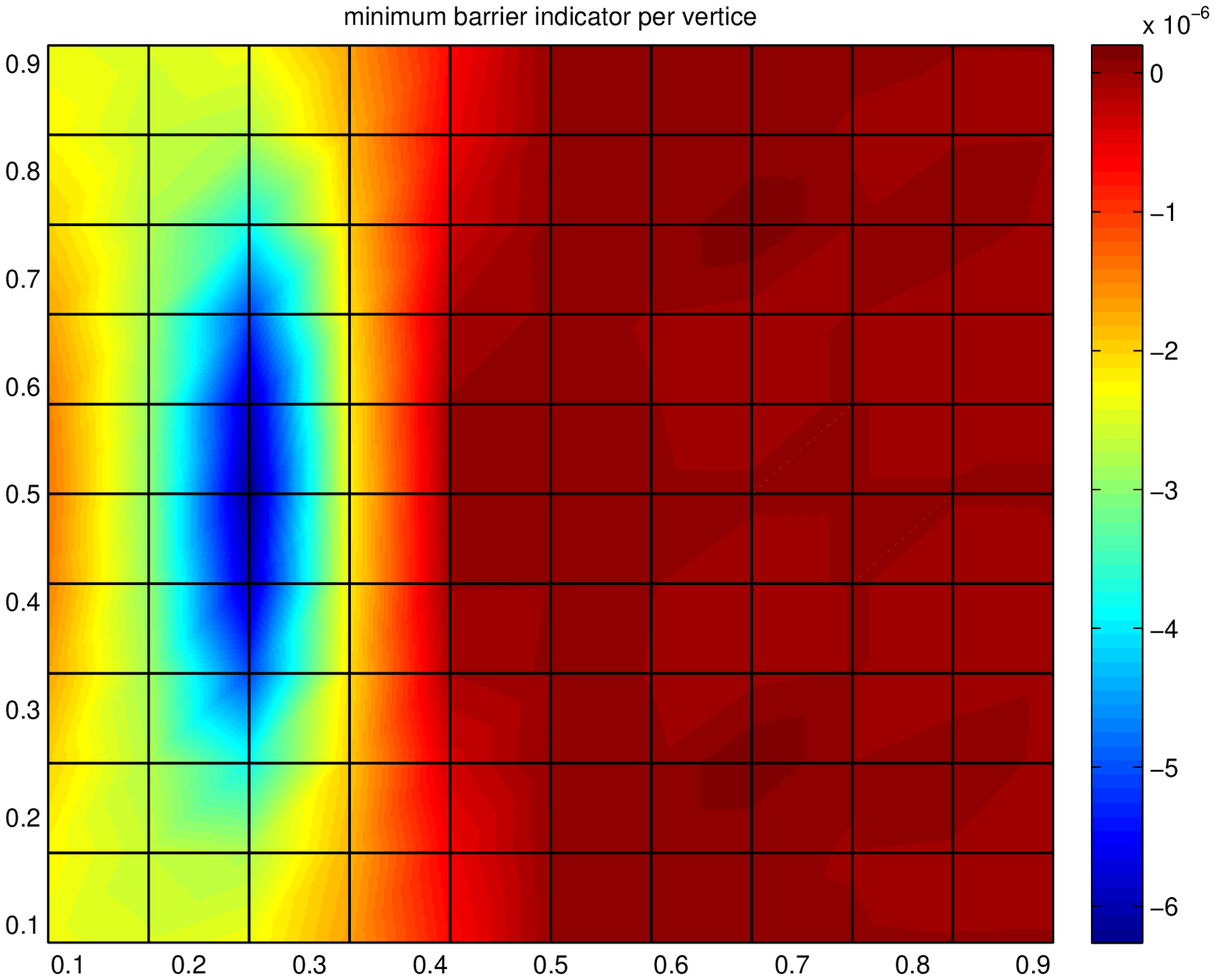}
    (g)
    \includegraphics[width=\twofigs]{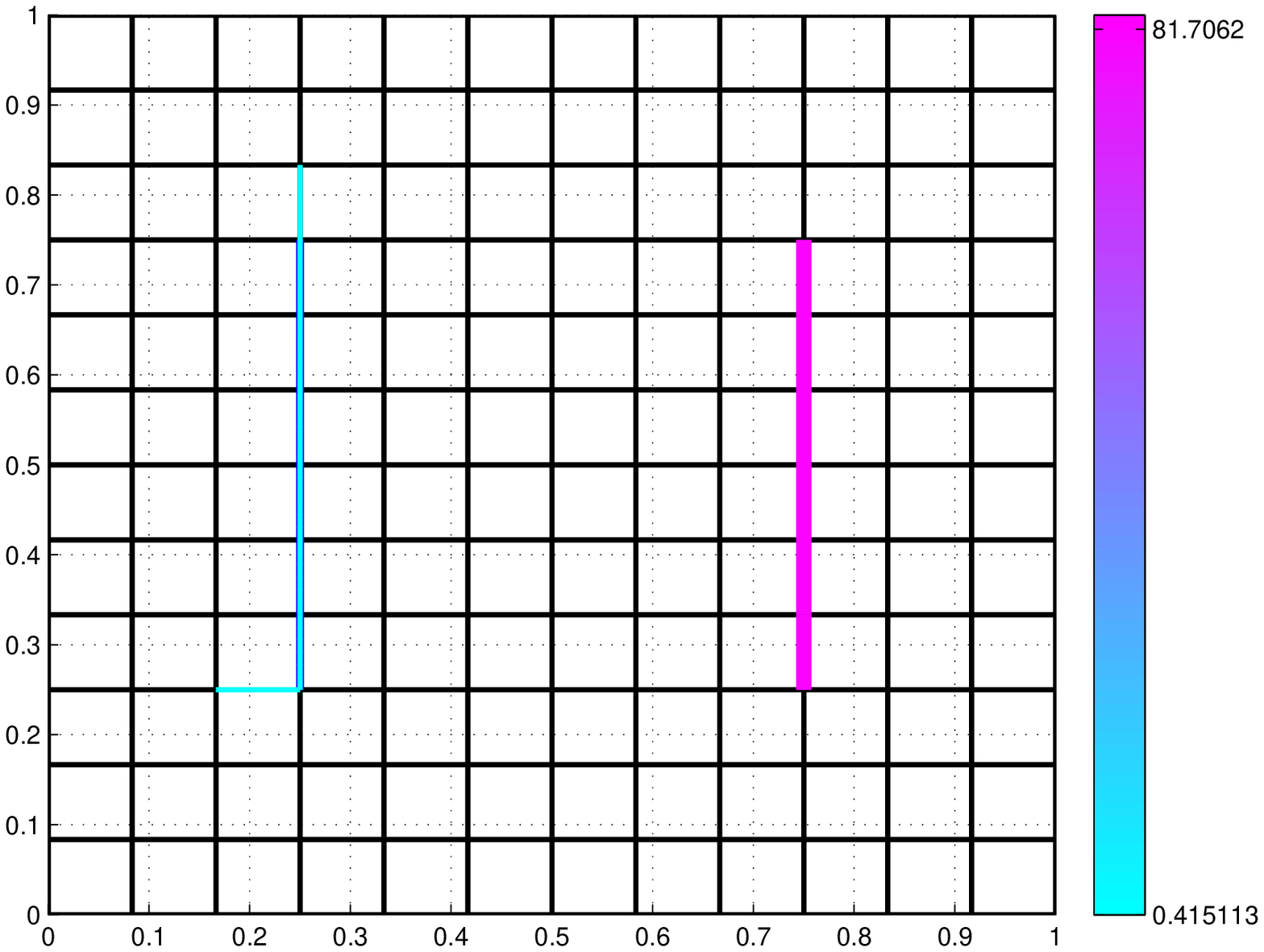}
    \caption{Cases of a single barrier (a-c), and of two barriers (d-g),
      normal to the flow.
      \newline
      (a) Data pressure, solution of the direct
      model~\eqref{eq:fracmodel:direct:discr:current} with a barrier at
      $x=0.5$ ($\bet=2$).
      \newline
      (b) Distribution of indicators for the long list of elementary
      {\candidate} barriers (see Figure~\ref{fig:one:fault}c).
      \newline
      (c) Best result after minimization for all {\candidate} barriers of the
      short list.
      \newline
      (d) Data pressure, solution of the direct
      model~\eqref{eq:fracmodel:direct:discr:current} with two barriers at
      $x=0.25$ ($\bet=2$) and $x=0.75$ ($\bet=20$).
      (e) Distribution of indicators at the first iteration.
      \newline
      (f) Distribution of indicators at the second iteration.
      \newline
      (g) Best result after minimization for all {\candidate} barriers of the
      short list.
    }
    \label{fig:barrier}
  \end{center}
\end{figure}

\myclearpage
\subsection{Estimation of barriers}
\label{ss:num:barrier}

The case of barriers is quite similar.
We present in Figure~\ref{fig:barrier} some tests for one or two barriers
normal to the flow.
The barrier indicators for the long list of elementary {\candidate} barriers
are even more precise: the target barriers are almost already drawn at the
beginning of the first iteration of the algorithm in the cartographies of
indicators in Figures~\ref{fig:barrier}b and~\ref{fig:barrier}e.
Indeed, the location of the most resistive target barrier is associated with
the strongest indicator, and the optimization step picks it up for it
produces a perfect fit to the data (Figures~\ref{fig:barrier}c
and~\ref{fig:barrier}g).
Nevertheless, in the case of two barriers, during the extension stage of the
second iteration, the location of the weakest target barrier is not
associated with the strongest indicator, and it is not even picked up by the
optimization step.
Consequently, although the best location for the second barrier contains all
edges of the least resistive target barrier, the optimization results are not
perfect.
The recovered values for the resistivity are $\bet=82$ and $\bet=0.41$,
instead of the target values~20 and~2;
they are quite far, but the hierarchy is correct.

\begin{figure}[\figloc]
  \begin{center}
    (a)
    \includegraphics[width=\twofigs]{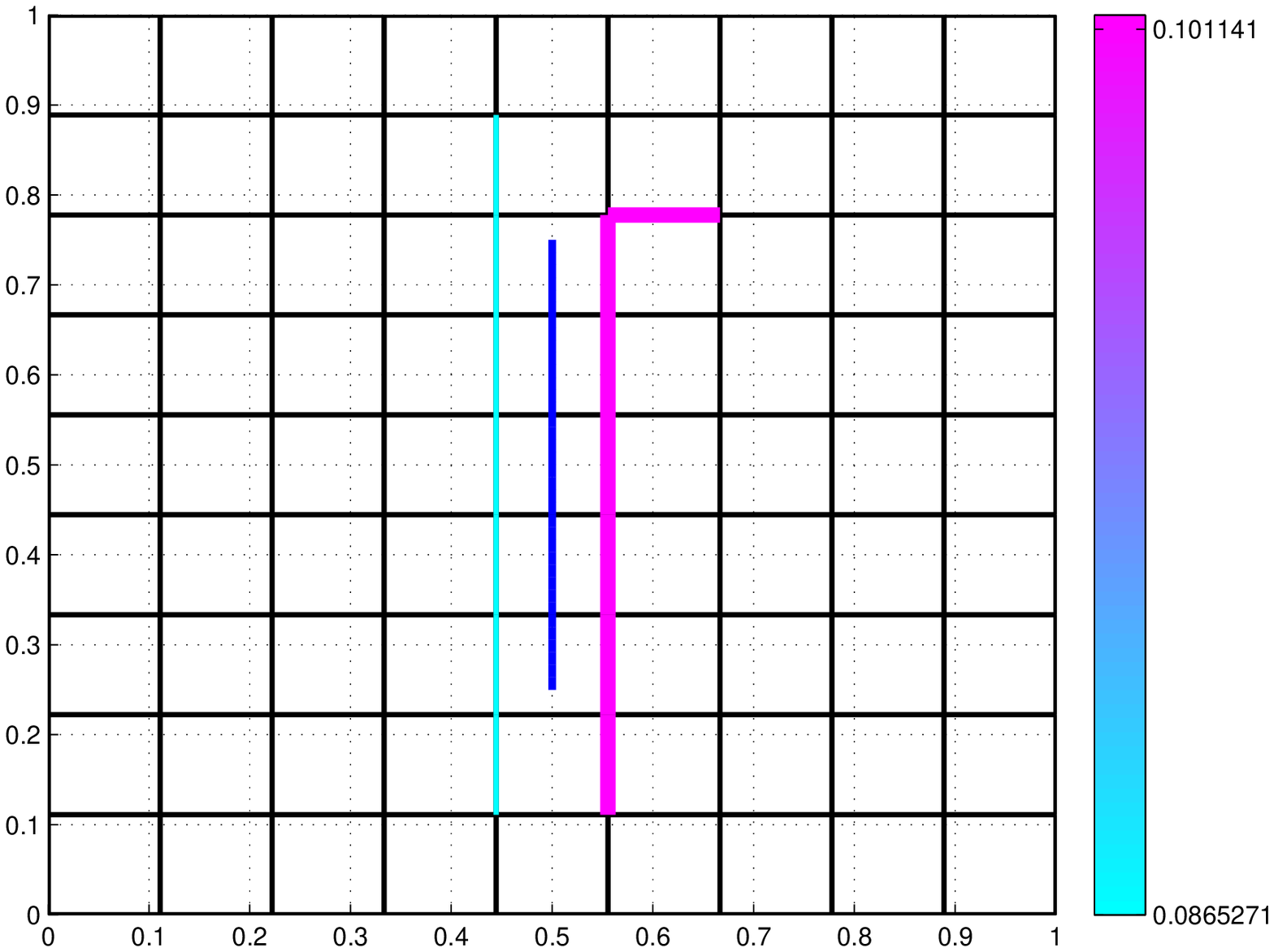}
    (b)
    \includegraphics[width=\twofigs]{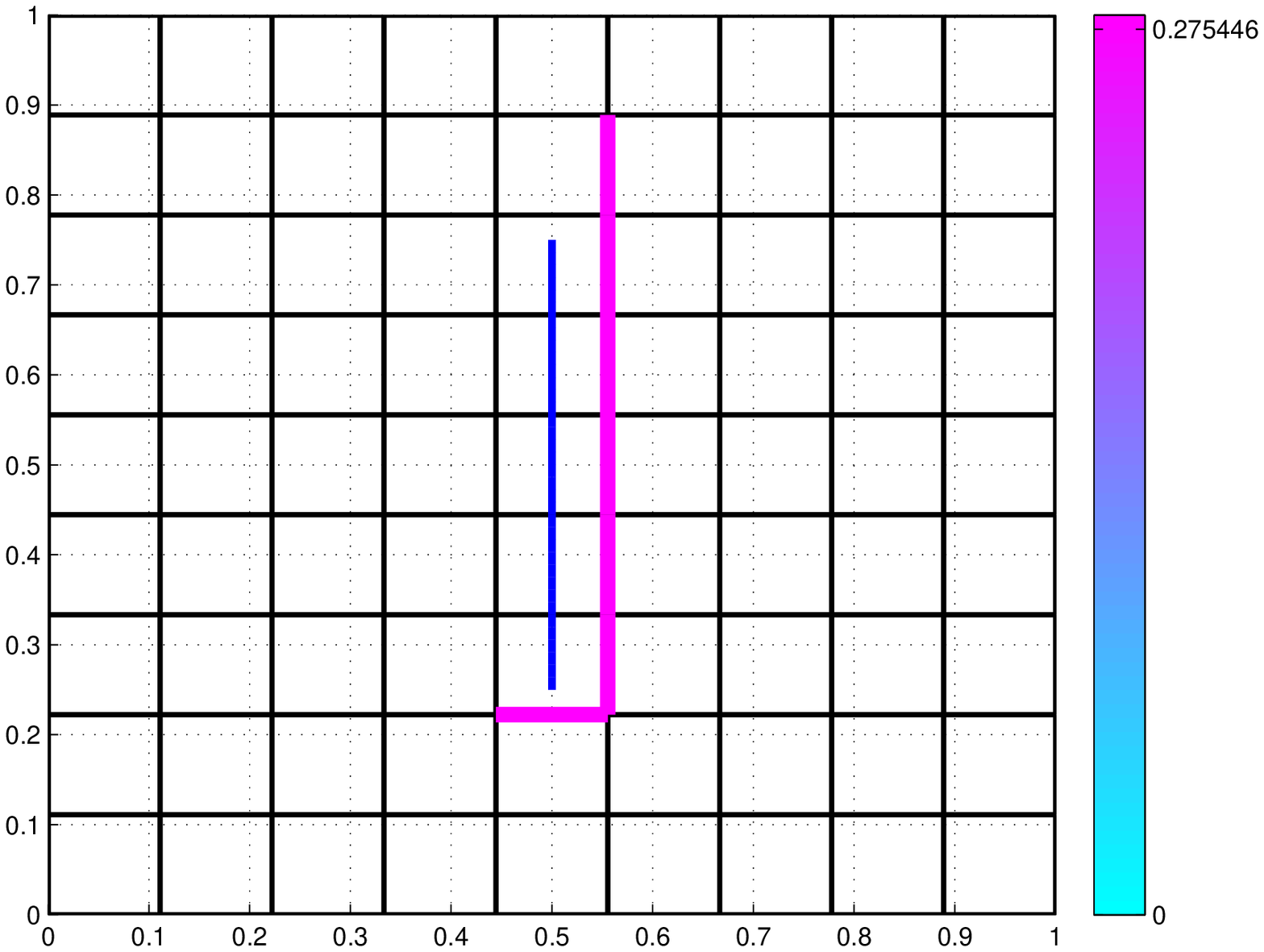} \\
    (c)
    \includegraphics[width=\twofigs]{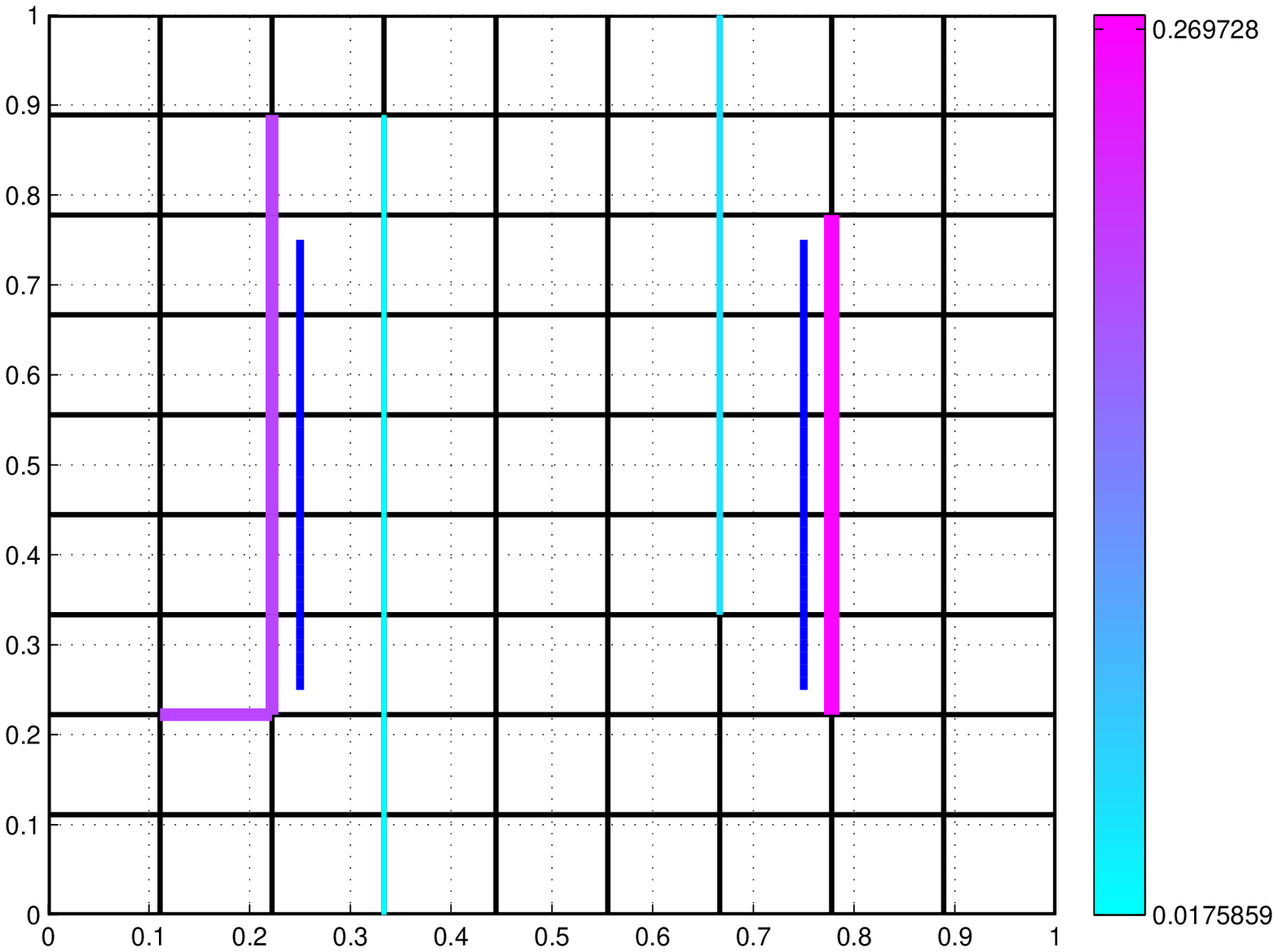}
    (d)
    \includegraphics[width=\twofigs]{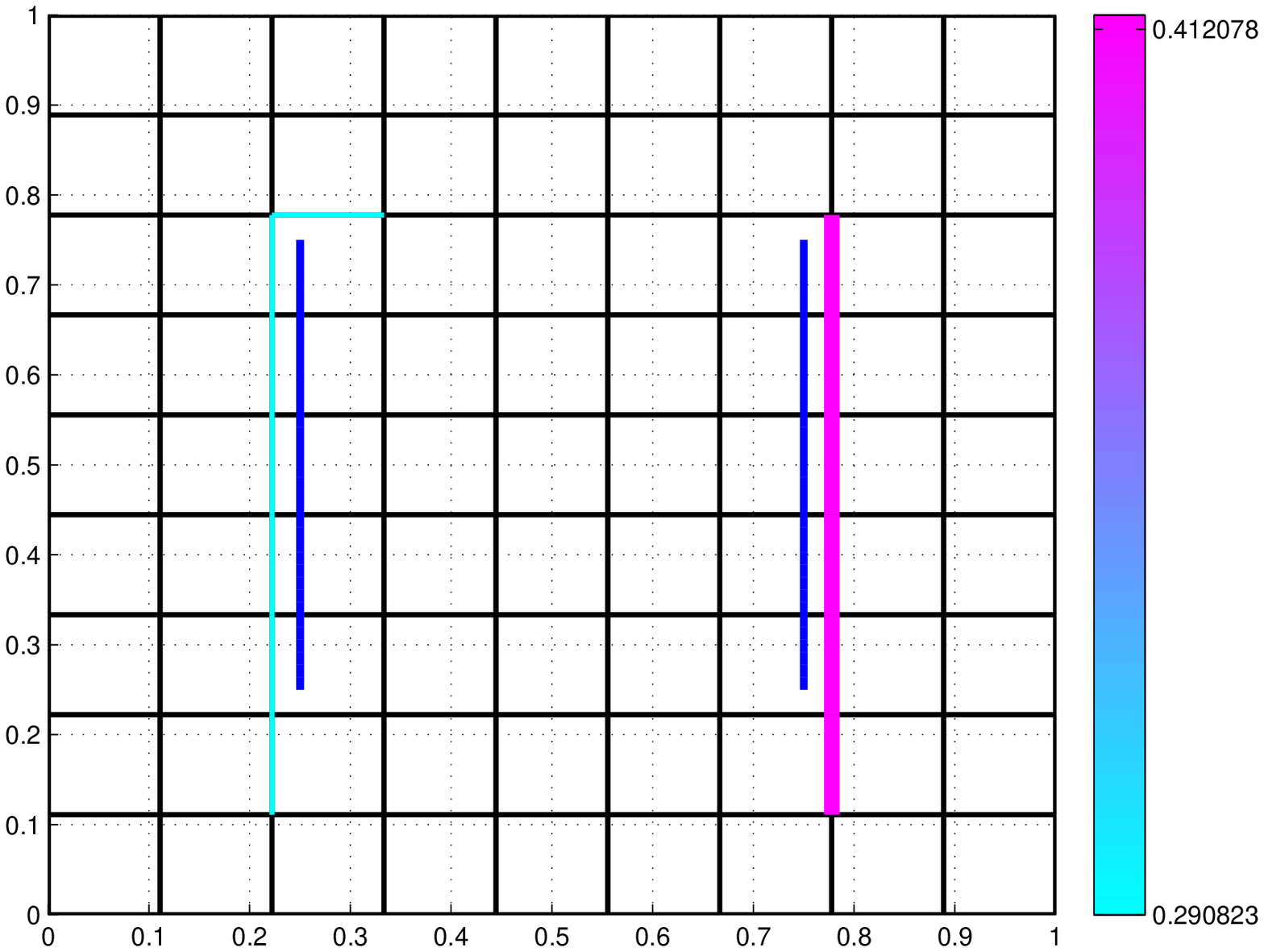}
    \caption{Estimation of barriers with coherent noise:
      the target barriers (dark blue) are not carried by the fracture
      research mesh ($\Nf=9$).
      The inverted barriers are represented with colors ranging from light
      blue (low resistivities) to pink (high resistivities).
      \newline
      (a) Single target barrier, $\Nm=72$.
      (b) Single target barrier, $\Nm=8$.
      \newline
      (c) Two target barriers, $\Nm=72$.
      (d) Two target barriers, $\Nm=8$.
    }
    \label{fig:barrier:noise}
  \end{center}
\end{figure}

When using a fracture research mesh that no longer carries the target
barriers ($\Nf=9$), the inversion results are still quite reasonable, see
Figure~\ref{fig:barrier:noise}.
We also decimate the pressure data measures from $\Nm=72$ down to $\Nm=8$.
The estimated locations of the barriers contains most of the coarse edges of
the fracture research mesh that are closest to the target barriers.
After minimization, the recovered effective resistivities (values of~$\bet$
are also adding up for nearby parallel barriers normal to the flow) are
underestimated by one or two orders of magnitude, but the hierarchy is still
preserved.

\myclearpage
\subsection{Discussion}
\label{ss:num:concl}

The first lesson learned from the numerical experiments is that, as hoped,
the inversion algorithm presented in Section~\ref{s:algorithm} performs
accurately in favorable situations.
Indeed, exact fit to synthetic data that are in the range of the inversion
model can be reached in a wide variety of situations: for one or two faults
tangential to the flow, for a barrier normal to the flow, for a wide range of
target parameter values (up to four orders of magnitude), and when loosening
the number of distributed measurements (up to a factor of about~80).

Much more interesting is the stability property shown by the algorithm.
The recovered fractures remain perfectly located, or at least superimposed
with the target location, when adding random noise up to a level
of~6\% for faults, and of~8\% for barriers, see~\cite{these:fatma:2016}.
And this remains true when loosening the number of measurements by a factor
of about~80 for up to~4\% noise (for faults) and up to~5\% noise for
barriers.
Moreover, when the situation deteriorates a bit more, for instance when the
measurement mesh becomes just too loose, when the noise level becomes just
too high, or when the target fractures no longer belong to the fracture
research mesh, then the location of recovered fractures remains within a
distance of no more than a coarse cell from the target location.
Furthermore, even when fractures are estimated by several smaller fractures
of close location, the value of the recovered effective parameter is of the
same order than the target value for single fractures, and the hierarchy in
values is kept in the case of two target fractures.

However, the heuristic strategy detailed in Section~\ref{ss:build:cand:frac}
was designed for the case of fault detection, and it may be suboptimal for
the estimation of barriers.
Indeed, even in the most advantageous case of the inverse crime with full
measurements, the indicator step of the algorithm is unable to nominate the
target location of the least resistive barrier in the case of two target
barriers normal to the flow.

Finally, the amount of measurements necessary to accurately recover the
location of fractures may seem prohibitive (64 in our tests).
This indicates that one needs a fine enough distribution of measurement to
hope solving this difficult inverse problem.
Some further experiments could deal with a study of the influence of the
location of measurements: {\apriori} information may help driving the
measurements in the vicinity of the target fracture in order to reduce their
number.

%% file: concl.tex

\section{Conclusion}

The estimation of fractures in porous media from pressure and/or flow data is
known to be an ill-posed inverse problem.
We propose a new approach for the recovering of both location and
hydrogeological properties of a small number of large fractures.
The method does not need any remeshing, nor shape derivation.

We have chosen a numerical flow model that treats fractures as interfaces
interacting with the surrounding porous matrix, and we have derived a forward
model for the simulation of both faults and barriers with vanishing
intensities, provided they are located on edges of the simulation mesh.
Then, the approach is based on the minimization of an error function that
measures the distance between data and simulated measurements.

We have defined specific fracture indicators, whose sign and modulus give the
effect on the data misfit of opening any {\candidate} fault or barrier.
These indicators give only first-order information, but are computationally
inexpensive.

Overparameterization is avoided by restricting the search for fractures to
edges of a coarser mesh than the computational grid, and by introducing
fractures one at a time.
When the data misfit is too high, a short list of new {\candidate} fracture
locations is built from a large number of elementary fracture locations that
span the entire domain of interest, with the assistance of fracture
indicators.
The actual enhancement brought by each candidate is computed by optimization,
and the best performing fracture is retained.

Numerical tests have been performed on synthetic data corresponding to simple
typical situations.
They demonstrate the ability of the proposed algorithm to automatically
retrieve one or two faults parallel to the flow, or one or two barriers
perpendicular to the flow.
The algorithm is shown to be fairly stable when noise is added to the data,
and when the fracture to be detected is not located on the fracture search
mesh.